\documentclass{aoscm}

\usepackage{apacite}

\usepackage{apacite,latexsym,graphics,epsfig,subfigure}

\def\sign{\mathop{\rm sign}\nolimits}
\def\half{{\scriptstyle\frac{1}{2}}}
\def\fourth{{\scriptstyle\frac{1}{4}}}
\def\threefourth{{\scriptstyle\frac{3}{4}}}

\newcommand{\cov}{{\mbox{\rm cov}}}
\newcommand{\var}{{\mbox{\rm var}}}

\newcommand{\tr}{\mbox{\rm tr}}
\newcommand{\erf}{\mbox{\rm erf}}
\newcommand{\erfi}{\mbox{\rm erfi}}

\newtheorem{theorem}{Theorem}
\newtheorem{lemma}{Lemma}

\newtheorem{corollary}{Corollary}

\NONUMBIB

\newcommand{\ind}{\perp\hspace*{-2.1mm}\perp}
\def\bR{{\bf R}}
\newtheorem{example}{Example}

\begin{document}

\SPECFNSYMBOL{1}{2}{}{}{}{}{}{}{}%
\AOSMAKETITLE

 \AOSyr{2006}
 \AOSvol{00}
 \AOSno{00}
 \AOSpp{000--000}
 \AOSReceived{Received }
 \AOSAMS{Primary 62H20; secondary 62E99}
 \AOSKeywords{correlation, covariance, test of independence,
spectral decomposition, eigenvalues, eigenfunctions, Hilbert-Schmidt
operator, Fr{\'e}chet bounds, contingency tables, phi-square,
canonical correlation, Cram{\'e}r-von Mises tests, rank tests,
Fredholm integral equation of the second kind.}
 \AOStitle{ON A NEW CORRELATION COEFFICIENT, ITS ORTHOGONAL DECOMPOSITION AND ASSOCIATED TESTS OF INDEPENDENCE}%
 \AOSauthor{Wicher Bergsma\footnote{Supported by The Netherlands
Organization for Scientific Research (NWO), Project Number
400-20-001. }}
 \AOSaffil{London School of Economics and Political Science}%
 \AOSlrh{Wicher Bergsma}%
 \AOSrrh{A new correlation coefficient}%
 \AOSAbstract{A possible drawback of the ordinary correlation coefficient $\rho$
for two real random variables $X$ and $Y$ is that zero correlation
does not imply independence. In this paper we introduce a new
correlation coefficient $\rho^*$ which assumes values between zero
and one, equalling zero iff the two variables are independent and
equalling one iff the two variables are linearly related. The
coefficients $\rho^*$ and $\rho^2$ are shown to be closely related
algebraically, and they coincide for distributions on a $2\times 2$
contingency table. We derive an orthogonal decomposition of $\rho^*$
as a positively weighted sum of squared ordinary correlations
between certain marginal eigenfunctions. Estimation of $\rho^*$ and
its component correlations and their asymptotic distributions are
discussed, and we develop visual tools for assessing the nature of a
possible association in a bivariate data set. The paper includes
consideration of grade (rank) versions of $\rho^*$ as well as the
use of $\rho^*$ for contingency table analysis. As a special case a
new generalization of the Cram{\'e}r-von Mises test to $K$ ordered
samples is obtained.
}%

\maketitle

\BACKTONORMALFOOTNOTE{3}




\noindent{\bf Contents:}

\contentsline {section}{\numberline {1}Introduction}{2}
\contentsline {section}{\numberline {2}Properties of the kernel
function $h_F$}{4} \contentsline {subsection}{\numberline {{\rm
2}.{\rm 1}}Key properties of $h_F$}{5} \contentsline
{subsection}{\numberline {{\rm 2}.{\rm 2}}Spectral decomposition of
$h_F$}{9} \contentsline {subsection}{\numberline {{\rm 2}.{\rm
3}}Obtaining the eigensystem in the discrete case}{13} \contentsline
{subsection}{\numberline {{\rm 2}.{\rm 4}}Obtaining the eigensystem
in the continuous case}{14} \contentsline {subsection}{\numberline
{{\rm 2}.{\rm 5}}Discrete approximation of the continuous case}{17}
\contentsline {subsection}{\numberline {{\rm 2}.{\rm 6}}Relation to
Anderson-Darling kernel}{18} \contentsline {section}{\numberline
{3}Properties of $\kappa $ and $\rho ^*$}{18} \contentsline
{subsection}{\numberline {{\rm 3}.{\rm 1}}Key properties of $\kappa
$ and $\rho ^*$}{19} \contentsline {subsection}{\numberline {{\rm
3}.{\rm 2}}Orthogonal decomposition}{22} \contentsline
{subsection}{\numberline {{\rm 3}.{\rm 3}}Parameterization of the
likelihood}{24} \contentsline {subsection}{\numberline {{\rm 3}.{\rm
4}}Fr{\'e}chet bounds for component correlations}{24} \contentsline
{section}{\numberline {4}Estimation and tests of independence}{27}
\contentsline {subsection}{\numberline {{\rm 4}.{\rm 1}}U and V
statistic estimators of $\kappa $}{27} \contentsline
{subsection}{\numberline {{\rm 4}.{\rm 2}}Permutation tests}{28}
\contentsline {subsection}{\numberline {{\rm 4}.{\rm 3}}Asymptotic
distribution of estimators under independence}{29} \contentsline
{subsection}{\numberline {{\rm 4}.{\rm 4}}Bonferroni corrections for
testing significance of component correlations}{30} \contentsline
{section}{\numberline {5}Grade versions of $\kappa $ and $\rho ^*$,
copulas, and rank tests}{31} \contentsline {subsection}{\numberline
{{\rm 5}.{\rm 1}}Rank tests}{32} \contentsline
{subsection}{\numberline {{\rm 5}.{\rm 2}}A new class of $K$-sample
Cram{\'e}r-von Mises tests as a special case}{32} \contentsline
{subsection}{\numberline {{\rm 5}.{\rm 3}}$\kappa $ as a weighted
$\phi $-coefficient}{33} \contentsline {section}{\numberline {6}Data
analysis: investigating the nature of the association}{34}
\contentsline {subsection}{\numberline {{\rm 6}.{\rm 1}}Some
artificial data sets}{35} \contentsline {subsection}{\numberline
{{\rm 6}.{\rm 2}}Mental health data}{41} \contentsline
{subsection}{\numberline {{\rm 6}.{\rm 3}}Norwegian stock
exchange}{43} \contentsline {subsection}{\numberline {{\rm 6}.{\rm
4}}Discussion}{43}

\section{Introduction}\label{sec1}

We introduce a correlation coefficient $\rho^*$ which has the
potential advantage compared to the ordinary correlation $\rho$ that
it detects arbitrary forms of association between two real random
variables $X$ and $Y$. In fact $\rho^*$, to be defined below, can be
viewed as a simple modification of $\rho^2$, as we now show. The
ordinary covariance is defined as
\begin{eqnarray*}
    \cov(X,Y) = E(X-EX)(Y-EY)
\end{eqnarray*}
and the ordinary correlation as
\begin{eqnarray*}
    \rho(X,Y) = \frac{\cov(X,Y)}{\sqrt{\cov(X,X)\cov(Y,Y)}}
\end{eqnarray*}
Now suppose that $Z$, $Z_1$ and $Z_2$ are iid with distribution
function $F$, that $X$ and $Y$ have marginal distributions $F_1$ and
$F_2$, and that $(X_1,Y_1),$ and $(X_2,Y_2)$ are iid replications of
$(X,Y)$. Then with
\begin{eqnarray}
    u_F(z_1,z_2) = (z_1-EZ)(z_2-EZ) = E(z_1-Z_1)(z_2-Z_2)
    \label{udef}
\end{eqnarray}
it is easy to verify that
\begin{eqnarray}
    \cov(X,Y)^2 = E
    u_{F_1}(X_1,X_2)u_{F_2}(Y_1,Y_2)
    \label{covdef}
\end{eqnarray}
Now straightforward algebra based on the left hand side
of~(\ref{udef}) shows that we can rewrite $u_F$ as
\begin{eqnarray*}
    u_F(z_1,z_2) =
   -
   \half E
   {\mbox{\Large $($}}\,
   |z_1-z_2|^2 - |z_1-Z_2|^2 - |Z_1-z_2|^2 + |Z_1-Z_2|^2
   \,{\mbox{\Large )}}
\end{eqnarray*}
Replacing the squares in $u_F$ by absolute values then gives
\begin{eqnarray}
   h_F(z_1,z_2) =
   -
   \half
   E{\mbox{\Large $($}}\,
   |z_1-z_2| - |z_1-Z_2| - |Z_1-z_2| + |Z_1-Z_2|
   \,{\mbox{\Large )}}
   \label{hdef1}
\end{eqnarray}
and we can define a new `covariance' $\kappa$ by replacing $u$ by
$h$ in~(\ref{covdef}):
\begin{eqnarray*}
   \kappa(X,Y) = Eh_{F_1}(X_1,X_2)h_{F_2}(Y_1,Y_2)
\end{eqnarray*}
Now we can also define
\begin{eqnarray*}
   \rho^*(X,Y)
   =
   \frac
   {\kappa(X,Y)}
   {\sqrt{\kappa(X,X)\kappa(Y,Y)}}
\end{eqnarray*}
Thus, whereas the squared covariance and $\rho^2$ are based on
squared differences, $\kappa$ and $\rho^*$ are based on absolute
differences. In this paper we demonstrate the perhaps surprising
result that $0\le \rho^*(X,Y)\le 1$, such that $\rho^*(X,Y)=0$ iff
$X$ and $Y$ are independent, and $\rho^*(X,Y)= 1$ iff $X$ and $Y$
are linearly related. A further main result we give is an orthogonal
decomposition of $\rho^*$ in terms of component correlations between
eigenfunctions of $h_{F_1}$ and $h_{F_2}$.

Based on their formulas, the following statistical interpretation of
$\rho^2$ and $\rho^*$ can be given: they measure how much two $X$
observations which are `far' apart tend to occur with $Y$
observations which are `far' apart, and similarly how much two $X$
observations which are `close' together tend to occur with $Y$
observations which are `close' together.

This paper is organized as follows. In Section~\ref{sec2}, the
properties of the kernel function $h_F$ are investigated in detail.
Some general properties are given, including conditions for its
existence and a proof that it is positive, and a large part of the
section is devoted to the spectral decomposition of $h_F$. We show
that if $h_F$ is square integrable, it has a mean square convergent
spectral decomposition in terms of the eigenvalues and vectors of
$h_F$. For discrete $F$, a set of difference equations is given
which has this eigensystem as its solution, and for continuous $F$
an analogous differential equation is given. The numerical solution
of these equations is treated in some detail. Closed form solutions
are only available in some special cases, for example, if $F$
belongs to the uniform distribution on $[0,1]$, the eigenfunctions
of $h_F$ are the Fourier cosine functions.

The results of Section~\ref{sec2} are used in Section~\ref{sec3} to
derive properties of $\rho^*$. We demonstrate the aforementioned
result that $0\le \rho^*(X,Y)\le 1$, such that $\rho^*(X,Y)=0$ iff
$X$ and $Y$ are independent, and $\rho^*(X,Y)= 1$ iff $X$ and $Y$
are linearly related. Furthermore, a decomposition of $\rho^*$ is
given in terms of a sum of squared correlations between marginal
eigenfunctions of $h_{F_1}$ and $h_{F_2}$ weighted with the product
of the corresponding marginal eigenvalues. We give a
parameterization of the likelihood in terms of component
correlations of $\rho^*$ and the marginal eigenfunctions, somewhat
analogous to the well-known canonical correlation decomposition.
Fr{\'e}chet bounds for the component correlations are discussed,
which gives some insight into the possible structure of the
dependence between two random variables. Finally in this section,
component correlations for the normal distribution are discussed as
an illustration.

In Section~\ref{sec4}, we derive sample and unbiased estimators of
$\kappa(X,Y)$ and related estimates of $\rho^*(X,Y)$, which can be
calculated in time $O(n^2)$. The asymptotic distributions of the
estimators under independence, which is a mixture of chi-squares, is
derived. Finally, small sample permutation tests and Bonferroni
corrections for testing the significance of component correlations
are discussed.

Section~\ref{sec5} concerns grade versions of $\kappa$ and $\rho^*$,
copulas, and rank tests. Rank statistics, obtained from the grade
versions of $\kappa$ and $\rho^*$, are discussed. It is shown that
the two sample Cram{\'e}r-von Mises statistic is obtained as a
special case, as well as a new generalization to the case of $K$
ordered samples. Furthermore, it is shown that $\rho^*$ is a
weighted mean of phi-square coefficients obtained from collapsing
the distribution onto a $2\times 2$ table with respect to cut-points
$(x,y)$.

In Section~\ref{sec6} we propose a methodology for gaining an
understanding of the association between two variables from a data
set. The methodology is based on combining hypothesis tests with
visual tools for displaying how much individual observations
contribute to the association.

Although many of the results of the present paper are new, we have,
naturally, also borrowed much from the literature, particularly
concerning the eigensystems and orthogonal decompositions. Some
important references here are \citeA{ad52,dk72}, \citeA{dwv73} and
\citeA{dewet87}, among others. However, the focus of much of the
literature we refer to is on studying power of hypothesis tests. The
aim of this paper, on the other hand, is on providing a meaningful
coefficient for describing association, which we hope leads to a
useful methodology for gaining an understanding of the association
between two variables, and, along the way, to tests with high power
against salient alternatives, the salience of the alternatives being
determined by the size of $\rho^*$.



Throughout this paper, we use the following conventions and
assumptions. We assume that $(X,Y)$, $(X_1,Y_1)$ and $(X_2,Y_2)$ are
iid with marginal distribution functions $F_1$ and $F_2$,
respectively, and joint distribution function $F_{12}$. We impose no
restrictions on the distributions, i.e., they may be continuous,
discrete, or mixed continuous-discrete. For simplification of some
of the derivations, we define distribution functions in the
following slightly non-standard way:
\begin{eqnarray*}
   F_1(x) = P(X<x) + \half P(X=x)\\
   F_2(y) = P(Y<y) + \half P(Y=y)
\end{eqnarray*}
and
\begin{eqnarray*}
   \lefteqn{F_{12}(x,y) = P(X<x,Y<y) + \half P(X<x,Y=y) + }\\
   && \half P(X=x,Y<y) + \fourth P(X=x,Y=y)
\end{eqnarray*}

\section{Properties of the kernel function $h_F$}\label{sec2}

In this section a detailed description is given of the kernel
function $h_F$ defined by~(\ref{hdef1}). Section~\ref{basichf}
concerns existence, continuity, positivity, square integrability,
existence of the trace and the shape of the graph of $h_F$. Methods
for verifying whether several of these properties hold are given. In
Section~\ref{spect}, under the assumption of square integrability of
$h_F$, its spectral decomposition is given, and some properties of
the associated eigenvalues and functions are derived. The
eigensystem is the solution to an integral equation which may be
difficult to solve. The problem is reformulated in terms of
difference equations for the discrete case in Section~\ref{obt eig
1} and in terms of a differential equation for the continuous case
in Section~\ref{obt eig 2}. Both rewrites appear much easier to
solve than the integral equation. In Section~\ref{contapprox},
efficient numerical approximation of the eigensystem of $h_F$ for
continuous $F$ is discussed. For several well-known distributions,
including the uniform and the normal, closed form solutions or
numerical approximations of (parts of) the eigensystem are given. A
new distribution $F$ is introduced which has the seemingly rare
property that $h_F$ is square integrable but has infinite trace. In
Section~\ref{sec ad} the relation between $h_F$ and a kernel
introduced by \citeA{ad52} is given. We are not aware of the kernel
$h_F$, depending on $F$, having been described previously.

\subsection{Key properties of $h_F$}\label{basichf}


The kernel $h_F$ {\em exists} if $h_F(z_1,z_2)$ is finite for some
$(z_1,z_2)\in\bR^2$. The kernel $h_F$ is {\em positive} if
\begin{eqnarray}
   Eg(Z_1)g(Z_2)h_F(Z_1,Z_2) \ge 0 \label{pd1}
\end{eqnarray}
for every function $g:\bR\rightarrow\bR$ for which the expectation
exists. The kernel $h_F$ is {\em square integrable} if
\begin{eqnarray}
   Eh_F(Z_1,Z_2)^2 = \int h_F(z_1,z_2)^2dF(z_1)dF(z_2)
   \label{hfsq}
\end{eqnarray}
is finite. The kernel $h_F$ is {\em trace class} if its trace
\begin{eqnarray*}
   \tr(h_F) = Eh_F(Z,Z) = \int h_F(z,z)dF(z)
\end{eqnarray*}
is finite. In several lemmas below we give some relatively easily
verifiable conditions for checking whether these properties hold for
$h_F$. The final Lemma~\ref{lem8} concerns the shape of the graph of
$h_F$.

The next lemma may simplify verification of the existence of $h_F$,
and asserts continuity and positivity as well as giving another
integral representation of $h_F$. First we need the following
notation:
\begin{eqnarray*}
   \gamma(x,y) =
   \left\{
   \begin{array}{cc}
      0     & x>y \\
      \half & x=y \\
      1     & x<y
   \end{array}
   \right.
\end{eqnarray*}
Note that
\begin{eqnarray}
    F(z)   &=& E\gamma(Z,z) \label{gamma1}\\
    1-F(z) &=& E\gamma(z,Z) \label{gamma2}
\end{eqnarray}
\begin{lemma}\label{lem1}
If $h_F$ exists it exists on $\bR^2$. It is then continuous and
positive, with equality in~(\ref{pd1}) only for the constant
function, and has the representation
\begin{eqnarray*}
   h_F(z_1,z_2) =
   \int_{-\infty}^{\infty}[\gamma(z_1, w)-F(w)][\gamma(z_2, w)-F(w)]dw
   \hspace{6mm}\forall z_1,z_2
\end{eqnarray*}
\end{lemma}
\proof{Proof} We first show continuity of $h_F$ on its domain. Let
$\delta>0$. Then if $|z_1-z_1'|<\delta$ and $|z_2-z_2'|<\delta$,
\begin{eqnarray*}
   \lefteqn{|h_F(z_1',z_2') - h_F(z_1,z_2)|}\\
   &=& \half |E\left[ (|z_1'-z_2'|-|z_1-z_2|)-(|z_1'-Z_2|-|z_1-Z_2|)-(|z_2'-Z_1|-|z_2-Z_1|) \right]| \\
   &\le& \half E\left[ 2\delta +\delta + \delta \right] \\
   &=& 2\delta
\end{eqnarray*}
Hence, $h_F$ is continuous and bounded on any finite domain.
Therefore, if $h_F$ exists in one point it exists on $\bR^2$.

We next derive the integral representation of $h_F$. We have
\begin{eqnarray}
   |z_1-z_2| = \int_{-\infty}^{\infty}[\gamma(z_1,w)-\gamma(z_2,w)]^2dw
   \label{difint}
\end{eqnarray}
and with $z_{i:4}$ the $i$th largest number in the set
$\{z_1,z_2,z_3,z_4\}$, we have
\begin{eqnarray}
   \lefteqn{
   \int_{-\infty}^{\infty}
   \left|
   [\gamma(z_1, w)-\gamma(z_3, w)]
   [\gamma(z_2, w)-\gamma(z_4, w)]
   \right|
   dw =}\nonumber\hspace*{30mm}\\
   &&
   \left\{
   \begin{array}{ll}
      0                & \mbox{if $z_1,z_3<z_2,z_4$ or $z_1,z_3>z_2,z_4$}  \\
      z_{3:4}-z_{2:4}  & \mbox{otherwise}
   \end{array}
   \right.
   \label{fubini1}
\end{eqnarray}
Now we can derive the desired result first using~(\ref{difint}),
then applying Fubini's theorem which is justified
because~(\ref{fubini1}) is finite, and finally using~(\ref{gamma1}):
\begin{eqnarray}
   h_F(z_1,z_2)
   &=& -\half E\int_{-\infty}^{\infty}
   \left(
   [\gamma(z_1, w)-\gamma(z_2, w)]^2
   -
   [\gamma(z_1, w)-\gamma(Z_2, w)]^2-\right.
   \nonumber\\
   &&
   \left.
   [\gamma(Z_1, w)-\gamma(z_2, w)]^2
   +
   [\gamma(Z_1, w)-\gamma(Z_2, w)]^2
   \right)
   dw
   \nonumber\\
   &=&
   E\int_{-\infty}^{\infty}
   [\gamma(z_1, w)-\gamma(Z_1, w)]
   [\gamma(z_2, w)-\gamma(Z_2, w)]
   dw
   \nonumber\\
   &=&
   \int_{-\infty}^{\infty}[\gamma(z_1, w)-F(w)][\gamma(z_2,
   w)-F(w)]dw\label{intrep}
\end{eqnarray}

Finally, we show positivity of $h_F$. Let $g$ be nonconstant. Then
using~(\ref{intrep}) and Fubini's theorem,
\begin{eqnarray*}
   \lefteqn{Eg(Z_1)g(Z_2)h_F(Z_1,Z_2)}\\
   &=&
   E\int_{-\infty}^{\infty}
   g(Z_1)[\gamma(Z_1, w)-F(w)]
   g(Z_2)[\gamma(Z_2, w)-F(w)]
   dw\\
   &=&
   \int_{-\infty}^{\infty}
   \left(
   Eg(Z)[\gamma(Z, w)-F(w)]
   \right)^2
   dw > 0
\end{eqnarray*}
Hence, $h_F$ is positive. If $g$ is constant it is easily verified
that the expression is zero.
 \hspace*{\fill}$\Box$\endproof

\noindent Note that from the lemma, it follows that for checking
existence of $h_F$, it suffices to check existence of $h_F(0,0)$.
Now $h_F(0,0)$ has the following convenient representations
\begin{eqnarray}
   h_F(0,0)
   &=&
   \half E
   \left(
      |Z_1|+|Z_2|-|Z_1-Z_2|
   \right) \nonumber\\
   &=&
   \int_{-\infty}^0F(z)^2dz + \int_0^{\infty}[1-F(z)]^2dz
   \label{hf00}
\end{eqnarray}
These representations are immediately verified from~(\ref{hdef1})
and from the representation of $h_F$ given in Lemma~\ref{lem1}.

An example of a random variable for which $h_F$ does not exist is
$Z=V^2$, where $V$ has a Cauchy distribution. This can be verified
by checking that~(\ref{hf00}) does not converge.

By giving some alternative representations of~(\ref{hfsq}), the next
lemma may be helpful in the verification of square integrability of
$h_F$:
\begin{lemma}\label{lem3}
We have:
\begin{eqnarray}
   Eh_F(Z_1,Z_2)^2
   =
   \frac16E(Z_{2:4}-Z_{3:4})^2
   =
   2\int_{z_1<z_2}F(z_1)^2[1-F(z_2)^2]dz_1dz_2
   \label{sqintrepr}
\end{eqnarray}
\end{lemma}
\proof{Proof} To prove the first equality, write $a_{ij}=|Z_i-Z_j|$.
Then, since for example $Ea_{12}=Ea_{13}$ and
$Ea_{34}a_{12}=Ea_{34}a_{15}$, we obtain
\begin{eqnarray*}
   Eh_F(Z_1,Z_2)^2
   &=&   \frac14 E
   \left(a_{12}-a_{13}-a_{24}+a_{34}\right)
   \left(a_{12}-a_{15}-a_{26}+a_{56}\right)
   \\
   &=&   \frac14 E
   \left(
   a_{12}^2
   -a_{12}a_{15}
   -a_{12}a_{26}
   +a_{12}a_{56}\right.
   \\ &&
   -a_{13}a_{12}
   +a_{13}a_{15}
   +a_{13}a_{26}
   -a_{13}a_{56}
   \\
   &&
   -a_{24}a_{12}
   +a_{24}a_{15}
   +a_{24}a_{26}
   -a_{24}a_{56}
   \\ &&
   \left.
   +a_{34}a_{12}
   -a_{34}a_{15}
   -a_{34}a_{26}
   +a_{34}a_{56}
   \right)
   \\
   &=&   \frac14 E
   \left(
   a_{12}^2
   -a_{12}a_{13}
   -a_{12}a_{24}
   +a_{12}a_{34}\right)
   \\
   &=& \frac{1}{16} E
   \left(
   a_{12}-a_{13}-a_{24}+a_{34}
   \right)^2
\end{eqnarray*}
It may now be verified that $a_{12}-a_{13}-a_{24}+a_{34}$ equals 0
if $Z_1,Z_4\le Z_2,Z_3$ or $Z_1,Z_4\ge Z_2,Z_3$, and equals $\pm
2(Z_{2:4}-Z_{3:4})$ otherwise. Hence, and because
$P(a_{12}-a_{13}-a_{24}+a_{34}\ne 0)=2/3$, we have
\begin{eqnarray*}
   Eh_F(Z_1,Z_2)^2
   &=&
   \frac{1}{16}
   \left(
   a_{12}-a_{13}-a_{24}+a_{34}
   \right)^2\\
   &=&
   \frac{1}{16}\times
   \frac{2}{3}\times
   4\,E(Z_{2:4}-Z_{3:4})^2
   = \frac16E(Z_{2:4}-Z_{3:4})^2
\end{eqnarray*}
which is the first part of the lemma.

To prove the second equality, first note that
\begin{eqnarray}
   E\gamma(Z,z_1)\gamma(Z,z_2)=\min\{F(z_1),F(z_2)\}
   \label{eqmin}
\end{eqnarray}
We now have using Lemma~\ref{lem1}, Fubini's theorem (for
justification see proof of Lemma~\ref{lem1}) and~(\ref{eqmin}),
\begin{eqnarray*}
   Eh_F(Z_1,Z_2)^2
   &=& E\int[\gamma(Z_1,z_1)-F(z_1)][\gamma(Z_2,z_1)-F(z_1)]dz_1 \times\nonumber\\
   &&   \int[\gamma(Z_1,z_2)-F(z_2)][\gamma(Z_2,z_2)-F(z_2)]dz_2\nonumber\\
   &=& \int E[\gamma(Z_1,z_1)-F(z_1)][\gamma(Z_1,z_2)-F(z_2)]\times \nonumber\\
   &&E[\gamma(Z_2,z_1)-F(z_1)][\gamma(Z_2,z_2)-F(z_2)]dz_1dz_2 \nonumber\\
   &=& \int [\min\{F(z_1),F(z_2)\}-F(z_1)F(z_2)]^2dz_1dz_2\nonumber\\
   &=& 2\int_{z_1<z_2} F(z_1)^2[1-F(z_2)]^2dz_1dz_2\nonumber
\end{eqnarray*}
 \hspace*{\fill}$\Box$\endproof

\noindent An example of a distribution for which $h_F$ exists but is
not square integrable is the Cauchy distribution. In particular,
$h_F(0,0)=2\pi^{-1}\log 2$, so $h_F$ exists. In this case,
nonexistence of $Eh_F(Z_1,Z_2)^2$ can most easily be verified using
the right hand side of~(\ref{sqintrepr}) and a computer algebra
package such as provided in Mathematica.

The next lemma may be helpful in verifying whether or not $h_F$ is
trace class:
\begin{lemma}\label{lemtr}
We have
\begin{eqnarray*}
   \tr(h_F) = \half E|Z_1-Z_2| = \int F(z)[1-F(z)]dz
\end{eqnarray*}
which is finite iff $Z$ has finite mean.
\end{lemma}
\proof{Proof} The first equality follows directly
from~(\ref{hdef1}), the second is well-known and can be found using
similar techniques as in the proof of the second equality in
Lemma~\ref{lem3}. As is well-known, integration by parts leads to
the representation of the mean as
\begin{eqnarray*}
   EZ = \int_0^{\infty}[1-F(z)]dz - \int_{-\infty}^0F(z)dz
\end{eqnarray*}
so the mean exists iff the terms on the right hand side exists. Now
since
\begin{eqnarray*}
   F(0)\int_0^{\infty}[1-F(z)]dz \le \int_0^{\infty}F(z)[1-F(z)]dz \le \int_0^{\infty}[1-F(z)]dz
\end{eqnarray*}
and
\begin{eqnarray*}
   [1-F(0)]\int_{-\infty}^0[1-F(z)]dz \le \int_{-\infty}^0F(z)[1-F(z)]dz \le \int_{-\infty}^0F(z)dz
\end{eqnarray*}
it follows that
\begin{eqnarray*}
   \tr(h_F) = \int_0^{\infty}F(z)[1-F(z)]dz + \int_{-\infty}^0F(z)[1-F(z)]dz
\end{eqnarray*}
exists iff $EZ$ exists.
 \hspace*{\fill}$\Box$\endproof

\noindent The quantity $E|Z_1-Z_2|$ is also called {\em Gini's mean
difference}. Note that, by Lemmas~\ref{lem3} and~\ref{lemtr}, both
$Eh_F(Z_1,Z_2)^2$ and $\tr(h_F)$ can be used as measures of
dispersion.

An example of a distribution function $F$ for which $h_F$ is square
integrable but not trace class is given in Example~\ref{ex1} in the
next subsection.

\begin{figure}
\begin{center}
\includegraphics[width=.6\linewidth,clip]{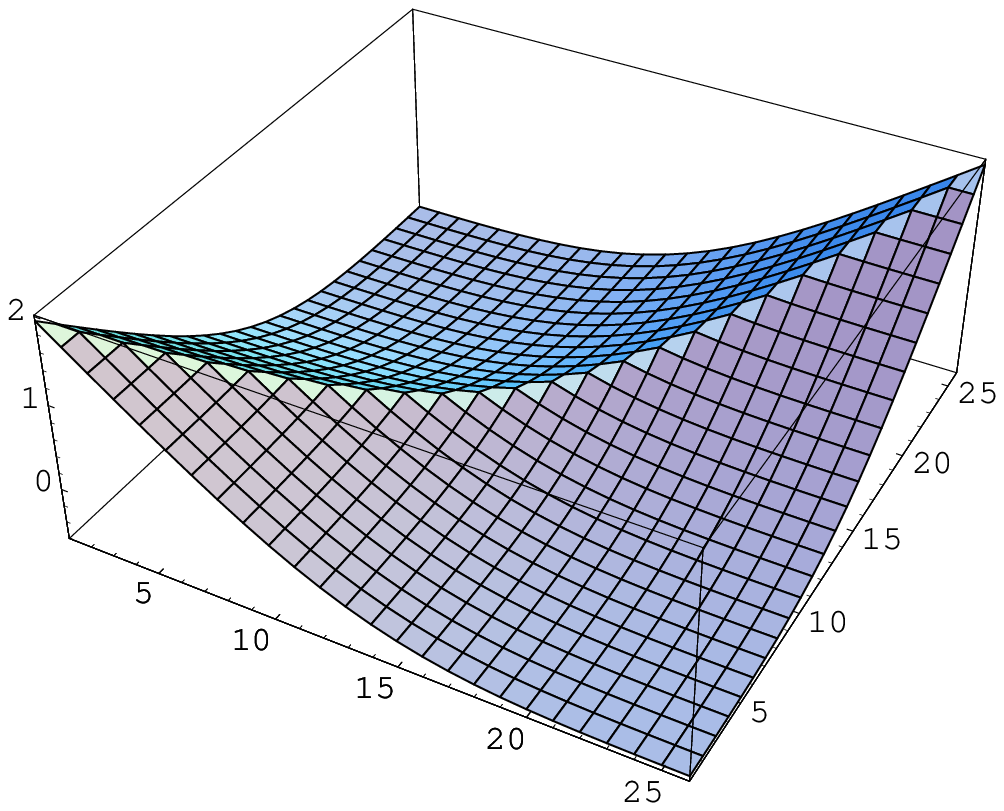}
\end{center}
{\caption{Graph of $h_F$ with $F$ the CDF for the standard normal
distribution{\label{hnorm}}}}
\end{figure}

We conclude this section with a lemma concerning the shape of the
graph of $h_F$. In Figure~\ref{hnorm} a representation of the graph
of $h_F$ with $F$ the CDF of the normal distribution is given. The
statements of Lemma~\ref{lem8} can be verified in the plot. Then:
\begin{lemma}\label{lem8}
Suppose $F$ is such that $h_F$ exists. Then:
\begin{enumerate}

\item For given $z_1$, $h_F(z_1,z_2)$ is strictly decreasing in $z_2$ on the domain $\{z:z\ge z_1,F(z)<1\}$ and strictly increasing in $z_2$
on the domain $\{z:z\le z_1,F(z)>0\}$.

\item $h_F(z,z)$ is strictly increasing in $z$ on the domain $\{z:F(z)>\half\}$ and strictly decreasing in $z$ on the domain $\{z:F(z)<\half\}$.

\end{enumerate}
\end{lemma}
\proof{Proof} Part 1: With $z_1,z_2,z_3$ such that $z_1<z_2\le z_3$
and $F(z_3)<1$, we obtain using Lemma~\ref{lem1} that
\begin{eqnarray*}
   h_F(z_1,z_3) - h_F(z_1,z_2)
   &=& \int[\gamma(z_1,w)-F(w)][\gamma(z_3,w)-\gamma(z_2,w)]dw \\
   &=& -\int_{z_2}^{z_3}[1-F(w)]dw\\
   &<& 0
\end{eqnarray*}
where the strict inequality holds because $F(z_3)<1$. This proves
the strict decreasingness part, the strict increasingness is proven
by appropriately reversing signs in the above.

Part 2: For $z<z'$ we have
\begin{eqnarray*}
   h_F(z,z) - h_F(z',z')
   &=&
   \int_{-\infty}^{\infty}[\gamma(z,w)-F(w)]^2-[\gamma(z',w)-F(w)]^2dw\\
   &=&
   \int_{z}^{z'}[1-2F(w)]dw
\end{eqnarray*}
which is positive if $F(z')\ge F(z)>\half$ and negative if $F(z)\le
F(z')<\half$, proving the monotonicity relations.
 \hspace*{\fill}$\Box$\endproof

\subsection{Spectral decomposition of $h_F$}\label{spect}

A sequence of random variables $Z^{(N)}$ is said to converge to $Z$
in mean square if
\begin{eqnarray*}
   E\left(Z - Z^{(N)}\right)^2 \rightarrow 0
   \mbox{ as } N\rightarrow\infty
\end{eqnarray*}
For a distribution function $F$, we define
\begin{eqnarray*}
   L_2(F) = \left\{g:\bR\rightarrow\bR\left|\int g(z)^2dF(x)<\infty\right.\right\}
\end{eqnarray*}
as the set of square integrable functions with respect to $F$. A set
of functions $\{g_k\}$ is said to be {\em orthonormal} if
\begin{eqnarray}
   Eg_k(Z)g_l(Z) = \delta_{kl} \label{g1}
\end{eqnarray}
(with $\delta$ the Kronecker delta) and
\begin{eqnarray*}
   Eg_k(Z)^2=1
\end{eqnarray*}
An orthonormal system $\{g_k\}$ is {\em complete} if for any
function $g\in L_2(F)$ there exist numbers $\{\alpha_k\}$ such that
\begin{eqnarray*}
   g(Z) = \sum_{k=1}^\infty\alpha_k g_k(Z)
\end{eqnarray*}
with convergence in mean square. Then the system $\{g_k\}$ is also
called a {\em basis} of $L_2(F)$. A number $\lambda$ is called an
eigenvalue of $h_F$ with corresponding eigenvector $g$ if
$Eg(Z)^2=1$ and
\begin{eqnarray}
   \lambda g(z) = Eh_F(z,Z)g(Z)  \label{g2}
\end{eqnarray}

Then we have:
\begin{theorem}\label{th2a}
Suppose $h_F$ is square integrable. Then there exists a complete
system of orthonormal functions $\{g_k\}$ of $L_2(F)$ consisting of
eigenfunctions of $h_F$, the corresponding eigenvalues
$\{\lambda_k\}$ being nonnegative. Each $g_k$ is continuous and
satisfies
\begin{eqnarray}
   Eg_k(Z) = 0 \label{degen}
\end{eqnarray}
if $\lambda_k\ne 0$. Furthermore, $h_F$ has the spectral
decomposition
\begin{eqnarray}
   h_F(Z_1,Z_2) = \sum_{k=1}^{\infty}\lambda_kg_k(Z_1)g_k(Z_2)
   \label{spectr1}
\end{eqnarray}
where convergence is in mean square.
\end{theorem}
\proof{Proof} A continuous square integrable kernel on a
sigma-finite measure space is a Hilbert-Schmidt kernel which by a
generalization of Mercer's theorem has the desired spectral
decomposition \cite{zaanen60}. It is easy to verify that $(\bR,B,F)$
is a $\sigma$-finite measure space, so the operator mapping the
function $g$ to $\int h_F(x,y)g(y)dF(y)$ is a Hilbert-Schmidt
operator. Nonnegativity of the eigenvalues follows from positivity
of $h_F$ (see Lemma~\ref{lem1}). Furthermore, from~(\ref{g2}),
$\lambda Eg_k(Z)=0$ so $Eg_k(Z)=0$ for nonzero $\lambda_k$.
 \hspace*{\fill}$\Box$\endproof

Note that if $g$ is a solution to~(\ref{g2}), then so is $-g$. To
identify the solutions, we proceed as follows. A function $g$ is
{\em initially positive} if there exists a $z$ such that $g(z)>0$
and $g(z')<0$ for all $z'<z$. Without loss of generality, we may
assume that the $g_k$ are initially positive. We also assume the
eigenvalues are ordered: $\lambda_1>\lambda_2>\ldots$. An
interesting property of eigenfunctions of homogeneous positive
Fredholm integral equations of the second kind is that they are
oscillating, in the sense that for every $k$, $g_k$ has $k$ distinct
zeroes and no more than $k-1$ local extremes.

Some further results are as follows:
\begin{lemma}\label{lem6}
For square integrable $h_F$:
\begin{enumerate}\label{lemtr2}

\item if finite, $\tr(h_F)=\sum_{k=0}^{\infty}\lambda_k$

\item $Eh_F(Z_1,Z_2)^2=\sum_{k=0}^{\infty}\lambda_k^2$.

\end{enumerate}
\end{lemma}
\proof{Proof} From the spectral decomposition~(\ref{spectr1}) and
Fubini's theorem,
\begin{eqnarray*}
   Eh_F(Z,Z) = E\sum\lambda_kg_k(Z)^2 = \sum\lambda_k
\end{eqnarray*}
For the second part, the desired result is obtained by Parseval's
theorem.
 \hspace*{\fill}$\Box$\endproof

As an example, we consider the dichotomous case which is the
simplest possible and has closed form solutions:
\begin{example}
Consider the dichotomous case that $P(Z=0)=1-P(Z=1)=p$. Then
\begin{eqnarray*}
   h_F(0,0) &=& (1-p)^2 \\
   h_F(0,1) &=& -p(1-p) \\
   h_F(1,0) &=& -p(1-p) \\
   h_F(1,1) &=& p^2
\end{eqnarray*}
The eigenvalues are $(\lambda_0,\lambda_1) = (0,p(1-p))$ and the
eigenfunctions are $g_0(z)=1$ and
\begin{eqnarray*}
   g_1(z) = \frac{1}{\sqrt{p(1-p)}}[z-(1-p)]
\end{eqnarray*}
which has mean zero and variance one. Now $h_F(z_1,z_2) =
\lambda_1g_1(z_1)g_1(z_2)$.
\end{example}
Solutions to~(\ref{g2}) for various other $F$ are given in
Table~\ref{eigtable}, see Section~\ref{obt eig 2} for an
explanation.

Equation~(\ref{g2}) is a homogeneous Fredholm integral equation of
the second kind based on the degenerate kernel $h_F$ (see, for
example, \citeNP{tricomi85}). In general these equations are
difficult to solve. In Sections~\ref{obt eig 1} and~\ref{obt eig
2},~(\ref{g2}) is reduced to a simpler problem for the discrete and
continuous case, respectively. We conclude this section with an
alternative formula to~(\ref{g2}) for finding the eigensystem which
we use in the next two subsections. With
\begin{eqnarray*}
   G(z) = \int_{-\infty}^{z}g(y)dF(y),
\end{eqnarray*}
we have:
\begin{lemma}\label{lem5}
The nonzero eigenvalues and eigenvectors of $h_F$ are solutions to
the equation
\begin{eqnarray*}
   \lambda
   \left[
   g(z) - g(z')
   \right]
   +
   \int_{z}^{z'}G(w)dw = 0
\end{eqnarray*}
for all $z<z'$.
\end{lemma}
\proof{Proof} Substitution of the expression for $h_F$ of
Lemma~\ref{lem1} into~(\ref{g2}) yields
\begin{eqnarray*}
   \lambda g (z)
   &=&
   \int_{-\infty}^{\infty}
   \int_{-\infty}^{\infty}
   [\gamma(z, w)-F(w)][\gamma(z_2, w)-F(w)]g (z_2)dF(z_2)dw\\
   &=&
   \int_{-\infty}^{\infty}
   [\gamma(z, w)-F(w)]G(w)dw\\
   &=&
   -\int_{-\infty}^{z}F(w)G (w)dw +
   \int_{z}^{\infty}[1-F(w)]G (w)dw
\end{eqnarray*}
The Lemma follows from this and since by~(\ref{degen})
$G(z)\rightarrow 0$ as $|z|\rightarrow\infty$.
 \hspace*{\fill}$\Box$\endproof

\noindent We have the following interesting implication:
\begin{corollary}\label{lindens}
If $\langle a,b\rangle$ is an interval with zero probability mass,
i.e., $F(a)=F(b)$, then a solution $g(z)$ to~(\ref{g2}) is linear on
$\langle a,b\rangle$.
\end{corollary}
\proof{Proof} If $F(a)=F(b)$ then from its definition it follows
that $G$ is constant on $\langle a,b\rangle$. Therefore, for any
$z,z'\in\langle a,b\rangle$, we obtain from Lemma~\ref{lem5} that
\begin{eqnarray*}
   \lambda
   \left[
   g(z) - g(z')
   \right]
   = c(z-z')
\end{eqnarray*}
for some constant $c$. Hence, $g$ is linear on $\langle a,b\rangle$.
 \hspace*{\fill}$\Box$\endproof

\noindent Note that for discrete distributions, it follows that the
eigenfunctions $g$ are piecewise linear.


\subsection{Obtaining the eigensystem in the discrete case}\label{obt eig 1}

We consider the case that $Z$ is a.s.\ in a finite set, say
$P(Z\in\{z_1,z_2,\ldots,z_K\})=1$. We use the shorthand
$P(Z=z_i)=p_i$ and assume without loss of generality that $p_i>0$
and $z_i<z_{i+1}$ for all $i$. Then we obtain:
\begin{lemma}\label{cor2}
With $c_i=(z_i-z_{i-1})^{-1}$ the nonzero eigenvalues and
eigenvectors of $h_F$ are solutions to the equations
\begin{eqnarray*}
   &&p_1g (z_1) = \lambda c_2[g (z_1)-g (z_2)] \\
   &&p_Kg (z_K) = \lambda c_K[g (z_{K-1})-g (z_{K})]\\
   && p_ig (z_i) =
   \lambda
   \left[c_ig (z_{i-1})-(c_{i}+c_{i+1})g (z_{i})+c_{i+1}g (z_{i+1})
   \right]
   \hspace{6mm}2\le i\le K-1
\end{eqnarray*}
\end{lemma}
\proof{Proof} First note that from the definition
for $z\in\langle z_i,z_{i+1}\rangle$,
\begin{eqnarray*}
   G (z) = \sum_{j=1}^{i}p_jg (z_j)
\end{eqnarray*}
It follows that
\begin{eqnarray}
   \int_{z_i}^{z_{i+1}}G (w)dw =
   (z_{i+1}-z_i)\sum_{j=1}^{i}p_jg (z_j)
   \label{int2}
\end{eqnarray}
The first two displayed equations of the lemma now follow
from~(\ref{int2}), Lemma~\ref{lem5} and from $\sum p_ig (z_i)=0$.
From~(\ref{int2}) we further obtain
\begin{eqnarray}
   c_{i+1}\int_{z_i}^{z_{i+1}}G (w)dw -c_{i} \int_{z_{i-1}}^{z_i}G (w)dw =
   p_ig (z_i)
   \label{dif1}
\end{eqnarray}
Again from Lemma~\ref{lem5}, we have
\begin{eqnarray*}
   \lambda
   \left[
   g (z_{i}) - g (z_{i+1})
   \right]
   +
   \int_{z_{i}}^{z_{i+1}}G (w)dw
   &=& 0
   \\
   \lambda
   \left[
   g (z_{i-1}) - g (z_{i})
   \right]
   +
   \int_{z_{i-1}}^{z_{i}}G (w)dw
   &=& 0
\end{eqnarray*}
Multiplying these equations by $c_{i+1}$ and $c_{i}$ respectively,
taking differences and the use of~\ref{dif1} yields
\begin{eqnarray*}
   p_ig_k(z_i) =
   \lambda
   \left[c_ig_k(z_{i-1})-(c_{i}+c_{i+1})g_k(z_{i})+c_{i+1}g_k(z_{i+1})
   \right]
\end{eqnarray*}
which completes the proof.
 \hspace*{\fill}$\Box$\endproof

More details on difference equations of the form given in
Lemma~\ref{cor2} can be found in \citeA{agarwal92}. Lemma~\ref{cor2}
allows fast and memory efficient computation of the eigenvalues and
vectors. In matrix notation, we must solve the generalized
eigenvalue problem
\begin{eqnarray*}
   D_pg = \lambda C g
\end{eqnarray*}
where $g$ is the eigenvector with corresponding eigenvalue
$\lambda$, $D_p$ is a diagonal matrix with the coordinates of $p$
on the main diagonal, and
\begin{eqnarray*}
   C =
   \left(
   \begin{array}{ccccccc}
      c_2 & c_2        & 0          & 0 \\
      c_2 & -(c_2+c_3) & c_3        & 0   &  & \ldots\\
      0   & c_3        & -(c_3+c_4) & c_4 & \\
      0   & 0          & c_4        & -(c_4+c_5) \\
          &            &            &        & \ddots \\
          & \vdots     &            &        & & -(c_{K-1}+c_K) & c_K      \\
          &            &            &        & & c_K         & c_K      \\
   \end{array}
   \right)
\end{eqnarray*}
This method can also be used for fast approximation of continuous
systems (see Section~\ref{contapprox}). The exact solution for
continuous systems is given in the next subsection.

\subsection{Obtaining the eigensystem in the continuous case}
\label{obt eig 2}

If $F$ is differentiable, the problem of finding the eigenvalues
and vectors can be reduced to a differential equation which is
sometimes easier to solve:
\begin{lemma}
Suppose $f$ is the derivative of $F$. Then the nonzero eigenvalues
and eigenvectors of $h_F$ are solutions to the equation
\begin{eqnarray}
   \lambda g''(z) + f(z)g (z) = 0
   \label{de1}
\end{eqnarray}
subject to the side condition
\begin{eqnarray*}
   g'(z) \rightarrow 0 \hspace{4mm}\mbox{ as }|z|\rightarrow\infty
\end{eqnarray*}
\end{lemma}
\proof{Proof} Letting $z'\rightarrow z$ in Lemma~\ref{lem5} we
obtain
\begin{eqnarray*}
   \lambda g'(z) + G (z) = 0
\end{eqnarray*}
Differentiating both sides with respect to $z$ yields~(\ref{de1}).
Finally, since by~(\ref{degen}) $G (z)\rightarrow 0$ as
$|z|\rightarrow\infty$, the side condition follows.
 \hspace*{\fill}$\Box$\endproof

Equation~(\ref{de1}) leads to an interesting observation on the
behavior of the eigenfunctions $g_k$: the second derivative
$g_k''(x)$ is proportional to the local density $f(x)$ times
$g_k(x)$. As mentioned in Section~\ref{spect} and as follows from
the theory of Sturm-Liouville differential equations, the
eigenfunctions oscillate. Now if the local density is low, the
oscillatory behavior will be slower. In particular, on intervals
where the local density is zero the second derivative of an
eigenfunction is zero so the eigenfunction is linear.  (NB: this
also holds for non-continuous distributions, see
Corollary~\ref{lindens}).

If $F$ is both differentiable and invertible, we can reformulate
differential equation~(\ref{de1}) in the standard Sturm-Liouville
form:
\begin{lemma}
Suppose $F$ is invertible. Let $q (u)=g [F^{-1}(u)]$ and suppose $q
$ is twice differentiable. Then the eigenvalues and eigenvectors are
solutions to the equation
\begin{eqnarray}
   \frac{d}{du}f[F^{-1}(u)]q '(u) + \lambda ^{-1}q (u) = 0
   \label{de2}
\end{eqnarray}
subject to the side condition
\begin{eqnarray*}
   f[F^{-1}(u)]q '(u) \rightarrow 0 \mbox{ as } u\downarrow 0 \mbox{ or } u\uparrow 1
\end{eqnarray*}
\end{lemma}
\proof{Proof} Note that
\begin{eqnarray*}
   \frac{d}{du}F^{-1}(u) =
   \frac{1}{f[F^{-1}(u)]}
\end{eqnarray*}
so that
\begin{eqnarray*}
   f[F^{-1}(u)]q '(u) = f[F^{-1}(u)]\frac{d}{du}g [F^{-1}(u)] =
   g '[F^{-1}(u)]
\end{eqnarray*}
From this the side condition follows. Substituting $z=F^{-1}(u)$
into the left hand side of~(\ref{de1}) yields
\begin{eqnarray*}
   \lefteqn{\lambda g ''[F^{-1}(u)] -  f(F^{-1}(u))g [F^{-1}(u)]}\\
   &=&
   \lambda \frac{du}{dF^{-1}(u)}\,\frac{d}{du}\,g '[F^{-1}(u)] -
   f(F^{-1}(u))g [F^{-1}(u)]\\
   &=&
   \lambda f[F^{-1}(u)]\,\frac{d}{du}\,f[F^{-1}(u)]q '(u) -  f[F^{-1}(u)]q (u)
\end{eqnarray*}
Hence dividing both sides of~(\ref{de1}) by $\lambda f[F^{-1}(u)]$
yields the desired result.
 \hspace*{\fill}$\Box$\endproof
Note that for the $q_k$ the orthonormality condition~(\ref{g1})
reduces to
\begin{eqnarray*}
   \int q_k(u)q_l(u)du = \delta_{kl}
\end{eqnarray*}

In general, the differential equations~(\ref{de1}) and~(\ref{de2})
do not have closed form solutions (i.e., solutions in terms of
well-known functions). We obtained the solutions for the uniform,
the logistic and the exponential distributions which are given in
Table~\ref{eigtable}, where $P_k$ is the $k$th Legendre polynomial,
$J_i$ is the $i$th Bessel function of the first kind, $\alpha_{k}$
the $k$th zero of $J_1$ and
\begin{eqnarray*}
   \beta_k=\left(J_0(\alpha_k)^2+J_1(\alpha_k)^2\right)^{-1/2}
\end{eqnarray*}
Numerical solutions were obtained for the normal, Laplace and
chi-square distribution with one degree of freedom, see the next
subsection for details on obtaining these solutions. For the normal
and Laplace distributions we obtained the exact value of
$\sum\lambda_k$ and for the normal distribution of $\sum\lambda_k^2$
using Lemmas~\ref{lem3} and~\ref{lemtr} combined with results for
order statistics by \citeA{bg59} and \citeA{govindarajulu63}
summarized in \citeA{jkb94} and using Mathematica.

We also obtained a closed form expression for the eigensystem for
the distribution introduced in the next example. It is an example of
a distribution $F$ for which $h_F$ is square integrable but not
trace class, i.e., by Lemma~\ref{lemtr2}, $\sum\lambda_k^2$ is
finite but $\sum\lambda_k=\infty$. We are not aware of any previous
studies of this distribution.
\begin{example}\label{ex1}
With $V$ standard normally distributed, let $Z$ be defined as the
following function of $V$:
\begin{eqnarray*}
   Z = \sqrt{\frac{\pi}2}\int_0^V \exp(t^2/2)dt
\end{eqnarray*}
where the convention is used that for $a>b$,
\begin{eqnarray*}
   \int_a^b f(t)dt=-\int_b^a f(t)dt
\end{eqnarray*}
Close to zero, $Z$ has approximately a normal density, but for large
values the density is much lower, that is, $Z$ has much heavier
tails than a standard normal.

With $F$ the distribution function of $Z$, we now show that $h_F$ is
square integrable but not trace class. To show this, we shall derive
an expression for $f[F^{-1}(u)]$ to be plugged into
Equation~(\ref{de2}) so that it can be solved. The derivation
involves the so-called complex error function. The {\em error
function} is defined as
\begin{eqnarray*}
  \erf(z) = \sqrt{2\pi}\int_0^z\exp(-t^2)dt
\end{eqnarray*}
Note that the CDF of the standard normal distribution is
$\Phi(v)=(1+\erf(v/\sqrt{2}))/2$. The {\em imaginary error function}
is defined as
\begin{eqnarray*}
  \erfi(z) = -i\,\erf(i\,z)
\end{eqnarray*}
(Here $i=\sqrt{-1}$. See \citeNP{weisstein:erf-erfi}, for some of
the properties of $\erf$ and $\erfi$.) We define the inverse
$\erfi^{-1}(z)$ as the unique real $y$ satisfying $z=\erfi(y)$.

We see that
\begin{eqnarray*}
  Z = \pi\,\erfi\left(\frac{V}{\sqrt{2}}\right)
\end{eqnarray*}
Now since $V$ has a standard normal distribution we obtain for the
CDF of $Z$:
\begin{eqnarray*}
  F(z)
  = P(Z<z)
  = P\left(V < \sqrt{2}\,\erfi^{-1}(z/\pi)\right)
  = \Phi\left(\sqrt{2}\,\erfi^{-1}(z/\pi)\right)
\end{eqnarray*}
From this,
\begin{eqnarray*}
  F^{-1}(u) = \pi\,\erfi\left(\Phi^{-1}(u)/\sqrt{2}\right)
\end{eqnarray*}
Some tedious but straightforward algebra then gives
\begin{eqnarray*}
  f(F^{-1}(u)) = 1/\phi[\Phi^{-1}(u)]^{2}
\end{eqnarray*}
Plugging this expression into~(\ref{de2}) leads to the solution
\begin{eqnarray*}
   \lambda_k &=& 1/k \\
   q_k(u) &\equiv& H_k\left[\Phi^{-1}(u)\right]\sqrt{\phi\left[\Phi^{-1}(u)\right]}
\end{eqnarray*}
where $H_k$ is the $k$th Hermite polynomial. The `$\equiv$' symbol
is used to indicate that the expression for $q_k$ needs to be
suitably normalized. This solution of~(\ref{de2}) was derived by
\citeA{dwv73}, who provided a method for solving differential
equations of the form $\frac{d}{du}w(u)q_k'(u) +
\lambda_k^{-1}q_k(u)=0$ for certain types of weights $w(u)$,
including $w(u)=1/\phi[\Phi^{-1}(u)]^{2}$. Note that for $g_k$ we
obtain
\begin{eqnarray*}
   g_k(x) \equiv H_k\left[\sqrt{2}\,\erfi^{-1}(\pi\,z)\right]\sqrt{\phi\left[\sqrt{2}\,\erfi^{-1}(\pi\,z)\right]}
\end{eqnarray*}
It is well-known that here $\sum\lambda_k$ is divergent and
$\sum\lambda_k^2=\pi^2/6$.
\end{example}

\begin{table}
\begin{center}
\begin{minipage}{120mm}
\begin{tabular}{l|llll}
 Distribution      & Logistic              & Uniform             & Normal                &  \\ \hline
$1/f[F^{-1}(u)]$  & $u(1-u)$                   & $1$             & $\phi[\Phi^{-1}(u)]$  &   \\
$\lambda_k$       & $\frac{1}{k(k+1)}$  & $\frac{1}{k^2\pi^2}$   & $^{\mbox{\it\scriptsize a}}$  \\
$\sum\lambda_k$   & $1$         & $\frac{1}{6}$                  & $\frac{1}{\sqrt\pi}$            \\
$\sum\lambda_k^2$ & $\frac13(\pi^2-9)$        & $\frac{1}{90}$   & $\frac13 - \frac{\sqrt{3}-1}{\pi}$           \\
$q_k(u)$          & $\sqrt{2k+1}\,P_k(2u-1)$       & $2\cos(k\pi u)$ & $^{\mbox{\it\scriptsize a}}$               \\
$\lambda_1/\sum\lambda_k$ & 0.5000        & 0.6079               & 0.5269 \\
$\lambda_2/\sum\lambda_k$ & 0.1667        & 0.1520               & 0.1635 \\
$\lambda_3/\sum\lambda_k$ & 0.0833        & 0.0675               & 0.0795 \\
$\lambda_4/\sum\lambda_k$ & 0.0500        & 0.0380               & 0.0470 \\
\\
 Distribution      & Exponential           & Laplace              & Chi-square  & Example~\ref{ex1} \\ \hline
$1/f[F^{-1}(u)]$  & $1-u$                 & $\min(u,1-u)$        & $\frac{\phi[\Phi^{-1}(\frac{u+1}{2})]}{\Phi^{-1}\left(\frac{u+1}{2}\right)}$   & $\phi[\Phi^{-1}(u)]^2$     \\
$\lambda_k$       & $\frac{4}{\alpha_k^2}$& $^{\mbox{\it\scriptsize a}}$                 & $^{\mbox{\it\scriptsize a}}$     & $\frac{1}{k}$  \\
$\sum\lambda_k$   & $\frac12$             & $\frac34$            & 0.6360  & $\infty$ \\
$\sum\lambda_k^2$ & $\frac{1}{12}$        & 0.1458               & 0.1399  & $\frac{6}{\pi^2}$  \\
$q_k(u)$          & $\beta_kJ_0(\alpha_{k}\sqrt{1-u})$ & $^{\mbox{\it\scriptsize a}}$     & $^{\mbox{\it\scriptsize a}}$     & $H_k[\Phi^{-1}(u)]\sqrt{\phi[\Phi^{-1}(u)]}^{\mbox{\it\scriptsize b}}$ \\
$\lambda_1/\sum\lambda_k$ & 0.5453        & 0.4611               & 0.5567  & 0\\
$\lambda_2/\sum\lambda_k$ & 0.1627        & 0.1816               & 0.1615  & 0\\
$\lambda_3/\sum\lambda_k$ & 0.0774        & 0.0875               & 0.0758  & 0\\
$\lambda_4/\sum\lambda_k$ & 0.0451        & 0.0542               & 0.0438  & 0\\
\end{tabular}
\end{minipage}
\caption{Eigenvalues and eigenvectors of kernel function $h_F$ for
various $F$. \protect{\newline} $\,^{\mbox{\it\scriptsize a}}$ No
closed form expression is available. \protect{\newline}
$\,^{\mbox{\it\scriptsize b}}$ Expression not normalized.
 \label{eigtable}}
\end{center}
\end{table}

The differential equation~(\ref{de2}) with $f[F^{-1}(z)]$ replaced
by a weight function $w(z)$ is given in \citeA{ad52}, see also
\citeA{dwv73} and \citeA{dewet87}. We are not aware of
equation~(\ref{de1}) having been studied previously.



\subsection{Discrete approximation of the continuous
case}\label{contapprox}

For many continuous distribution functions $F$ the differential
equations~(\ref{de1}) and~(\ref{de2}) do not have a closed form
solution and the first $t$ eigenvalues and eigenfunctions can be
approximated by using a discrete approximation of $F$ and solving
the difference equations of Lemma~\ref{cor2}. For $i=1,\ldots,t$,
let $c_i=F^{-1}\left(\frac{i-1/2}{t}\right)$ and let $Z^{(t)}$ be a
discrete random variable with
$P(Z^{(t)}=c_i)=p_i=F(i/t)-F((i-1)/t)$. The eigenvalues and
eigenvectors of $h_{F^{(t)}}$, with $F^{(t)}$ the distribution
function of $Z^{(t)}$, can then be calculated using the method of
Section~\ref{obt eig 1}. For large $t$, this method seems to give
good approximations of the eigenvalues and eigenvectors of $h_F$. An
idea of the quality of the approximations can be gained from
Table~\ref{eigtable2}. The numerical results in Table~\ref{eigtable}
were obtained using this method. For further details on discrete
approximations of eigenvalues and vectors of kernels see
\citeA{tricomi85}.

To obtain a good approximation, $t$ should of course be chosen as
large as possible. Using Mathematica~5.2 on a Pentium~IV computer at
3.0MHz, using no special routines for calculating the eigensytem of
tridiagonal matrices, calculation of a complete solution for
$t=1000$ took 29 seconds. We expect that using software with such
special routines, it is possible to obtain solutions of the first
few eigenvalues and eigenvectors much more quickly and for much
larger $t$.

\begin{table}
\begin{center}
\begin{tabular}{l|l|l|l}
                  & True value                & Estimate ($t=101$)        & Estimate ($t=1001$)         \\
$\sum\lambda_k$   & 1                         & 0.99303                   & 0.99931           \\
$\sum\lambda_k^2$ & 0.28987                   & 0.29027                   & 0.28988           \\
$\lambda_1^*$     & 0.50000                   & 0.50370                   & 0.50035   \\
$\lambda_{10}^*$  & 9.0909$\times 10^{-3}$    & 9.3093$\times 10^{-3}$    & 9.1056$\times 10^{-3}$\\
$\lambda_{100}^*$ & 9.9010$\times 10^{-5}$    & 9.9708$\times 10^{-5}$    & 9.9145$\times 10^{-5}$     \\
$\lambda_{1000}^*$& 9.9900$\times 10^{-7}$    & -                         & 9.9970$\times 10^{-7}$    \\
\end{tabular}
\caption{Eigenvalues and their estimates based on discrete
approximations for the kernel $h_F$ with $F$ the logistic
distribution function.
 \label{eigtable2}}
\end{center}
\end{table}

For calculation of the eigensystem from a sample, see
Section~\ref{uv}.

\subsection{Relation to Anderson-Darling kernel}\label{sec ad}

A related kernel was studied by \citeA{ad52} and \citeA{dwv73},
among others, in the context of Cram{\'e}r-von Mises tests. With $w$
a nonnegative weight function, they considered the kernel
\begin{eqnarray*}
   r_w(u,v) = \int_0^1
   [\gamma(u, t)-t]
   [\gamma(v, t)-t]
   w(t)
   dt
\end{eqnarray*}
The kernel $r_w$ is closely related to the kernel $h_F$: with
$w(t)=\frac{1}{f[F^{-1}(t)]}$, we obtain
\begin{eqnarray*}
   r_w(u,v) = h_F[F^{-1}(u),F^{-1}(v)]
\end{eqnarray*}
Note that this conversion does not work for discrete $F$. The
eigensystem for $r_w$ is given by the set of solutions
to~(\ref{de2}) with $f[F^{-1}(z)]$ replaced by the weight function
$w(z)$ \cite{ad52}.



Other closely related kernels have been given in the context of
two-sample tests by \citeA{za03} and \citeA{bf04}. (See
Section~\ref{cvm} for the relation between two-sample tests and
tests of independence.)

\section{Properties of $\kappa$ and $\rho^*$}\label{sec3}

We now apply the results of Section~\ref{sec2} in order to derive
properties of $\kappa$ and $\rho^*$. In Section~\ref{sec3.1}, some
key properties are derived, including that $0\le\rho^*(X,Y)\le 1$,
with $\rho^*(X,Y)=0$ iff $X\ind Y$ and $\rho^*(X,Y)=1$ iff $X$ and
$Y$ are linearly related. In Section~\ref{orthog} we give a
decomposition of $\kappa(X,Y)$ and $\rho^*(X,Y)$ as weighted sums of
squared correlations between the marginal eigenfunctions of
$h_{F_1}$ and $h_{F_2}$, weighted by functions of the eigenvalues.
In Section~\ref{parlik}, we describe a decomposition of the
likelihood in terms of marginal eigenfunctions and component
correlations of $\rho^*$. In Section~\ref{frechet}, Fr{\'e}chet
bound for the component correlations are given, which gives some
insight into their meaning.

\subsection{Key properties of $\kappa$ and $\rho^*$}\label{sec3.1}

Some key properties of $\kappa$ and of $\rho^*$, are given in the
following theorem:
\begin{theorem}\label{th1}
Suppose $X$ and $Y$ and $Z$ are real random variables for which the
marginal kernels $h_{F_1}$ and $h_{F_2}$ exist. Then:
\begin{enumerate}
\item If $a$, $b$, $c$ and $d$ are constants, then
$\kappa(aX+b,cY+d)=a\,c\,\kappa(X,Y)$ and
$\rho^*(aX+b,cY+d)=\rho^*(X,Y)$.\label{th1p1}

\item $\kappa(X,Y)\ge 0$ with equality iff $X\ind Y$.\label{th1p2}

\item If $\kappa(X,X)<\infty$ and $\kappa(Y,Y)<\infty$ then
$\kappa(X,Y)\le\sqrt{\kappa(X,X)\kappa(Y,Y)}$ with equality iff $X$
and $Y$ are a.s.\ linearly related.\label{th1p3}

\item If both $X$ and $Y$ are dichotomous then
$\kappa(X,Y)=\cov(X,Y)^2$ and $\rho^*(X,Y)=\rho(X,Y)^2$\label{th1p4}

\item With $Z_{i:n}$ the $i$th order statistic in a sample of size $n$, $\kappa(Z,Z)=\frac{1}{6}E(Z_{2:4}-Z_{3:4})^2$.\label{th1p5}

\end{enumerate}
\end{theorem}
The proof of the theorem is given at the end of this section. Note
that $\kappa(X,X)=Eh_{F_1}(X_1,X_2)^2$ and
$\kappa(Y,Y)=Eh_{F_2}(Y_1,Y_2)^2$ so the condition in Part~3 is
equivalent to square integrability of the marginal kernels. From
Theorem~\ref{th1}, Parts~2 and~3, we immediately have the following:
\begin{corollary}\label{cor1}
Suppose $\kappa(X,X)<\infty$ and $\kappa(Y,Y)<\infty$. Then
$0\le\rho^*(X,Y)\le 1$, with $\rho^*(X,Y)=0$ iff $X\ind Y$ and
$\rho^*(X,Y)=1$ iff $X$ and $Y$ are a.s.\ linearly related.
\end{corollary}
We may compare Corollary~\ref{cor1} to the related well-known result
for the ordinary correlation $\rho$: If $\var(X)<\infty$ and
$\var(Y)<\infty$ then $0\le\rho^2\le 1$ with $\rho^2=0$ if $X\ind Y$
and $\rho^2=1$ iff $X$ and $Y$ are a.s.\ linearly related. The
important difference is that $X\ind Y$ is equivalent to $\rho^*=0$
but $X\ind Y$ only implies $\rho=0$, not vice versa. By
Part~\ref{th1p5} of Theorem~\ref{th1}, $\kappa(Z,Z)$ can be used as
measures of dispersion for a real random variable $Z$. Note the
relation with the variance, which can be defined as
$E(Z_{1:2}-Z_{2:2})^2$.

Note that, even though $\rho^*(X,Y)=1$ iff $X$ and $Y$ are linearly
related, $\rho^*$ is not a measure of linear association in the
following sense: if the slope of the linear regression line of $Y$
given $X$ is zero, $\rho^*(X,Y)$ need not be equal to zero.


From Lemma~\ref{lem3} and Lemma~\ref{lemtr}, a sufficient condition
for $\rho^*(X,Y)$ to exist is that $EX$ and $EY$ exist. An example
showing that this condition is not necessary is Example~\ref{ex1}.
Note that the existence of the ordinary correlation $\rho$ has the
much stronger requirement of finite marginal variances. Summarizing,
we have
\begin{eqnarray*}
    \lefteqn{
    \rho(X,Y)\mbox{ exists }
    \Leftrightarrow
    \left\{\sigma^2(X) \mbox{ and }\sigma^2(Y)\mbox{ exist }\right\}
    \Rightarrow
    \left\{EX\mbox{ and }EY\mbox{ exist }\right\}
    \Rightarrow
    }\\
    &&
    \Rightarrow\left\{E(X_{2:4}-X_{3:4})^2\mbox{ and }E(Y_{2:4}-Y_{3:4})^2\mbox{ exist }\right\}
    \Leftrightarrow \rho^*(X,Y)\mbox{ exists }
\end{eqnarray*}
where the one-way implications are strict.

We now proceed to the proof of Theorem~\ref{th1}:\vspace{3mm}\\
\proof{Proof of Theorem~\ref{th1}}

Part 1 follows directly from the definition.

Part 2: From Lemma~\ref{lem1} we obtain
\begin{eqnarray*}
   h_{F_i}(a,b)
   &=& \int
   [\gamma(a,w)-F_i(w)][\gamma(b,w)-F_i(w)]
   dw
\end{eqnarray*}
Furthermore, note that
\begin{eqnarray*}
   {F_{12}(x,y)-F_1(x)F_2(y)}
   = E[\gamma(X,x)-F_1(x)][\gamma(Y,y)-F_2(y)]
\end{eqnarray*}
which is easy to verify. Using these results and the finiteness of
each side of~(\ref{fubini1}) which allows us to apply Fubini's
theorem, we obtain
\begin{eqnarray}
   \kappa(X,Y)
   &=& Eh_{F_1}(X_1,X_2)h_{F_2}(Y_1,Y_2)\nonumber\\
   &=& E\int[\gamma(X_1,x)-F_1(x)][\gamma(X_2,x)-F_1(x)]dx \times\nonumber\\
   &&   \int[\gamma(Y_1,y)-F_2(y)][\gamma(Y_2,y)-F_2(y)]dy\nonumber\\
   &=& \int E[\gamma(X_1,x)-F_1(x)][\gamma(Y_1,y)-F_2(y)]\times \nonumber\\
   &&E[\gamma(X_2,x)-F_1(x)][\gamma(Y_2,y)-F_2(y)]dxdy \nonumber\\
   &=& \int [F_{12}(x,y)-F_1(x)F_2(y)]^2dxdy\label{fxy}
   \label{intf12}
\end{eqnarray}
If $X\ind Y$, the integrand is zero so $\kappa(X,Y)=0$. It remains
to be shown that $X\not\!\!\!\ind Y$ implies $\kappa(X,Y)\ne 0$.  We
next sketch the proof.

If $X\not\!\!\!\ind Y$ then there is an $(a,b)$ such that
$D(a,b)=F_{12}(a,b)-F_1(a)F_2(b)\ne 0$. We now show that if
$D(a,b)\ne 0$, then there is an open interval, which has $(a,b)$ as
a limit point, such that $D\ne 0$ on that interval. It then follows
that~(\ref{intf12}) is nonzero. Let $\varepsilon_{12}\ge 0$ be the
probability mass in $(a,b)$, $\varepsilon_{1}\ge 0$ be the marginal
probability mass in $a$ and $\varepsilon_{2}\ge 0$ be the marginal
probability mass in $b$. Denote by $D(a^\pm,b^\pm)$ the limit
approaching from anywhere in one of four open `quadrants' defined by
$(a,b)$. Then from the definition of $F_{12}$, $F_1$ and $F_2$ given
in the introduction,
\begin{eqnarray*}
   D(a^-,b^-) &=& D(a,b) - \fourth\varepsilon_{12} + \half\varepsilon_{1}F_2(b) + \half\varepsilon_{2}F_1(a)
   - \fourth\varepsilon_{1}\varepsilon_{2} \\
   D(a^-,b^+) &=& D(a,b) + \fourth\varepsilon_{12} + \half\varepsilon_{1}F_2(b) - \half\varepsilon_{2}F_1(a)
   + \fourth\varepsilon_{1}\varepsilon_{2} \\
   D(a^+,b^-) &=& D(a,b) + \fourth\varepsilon_{12} - \half\varepsilon_{1}F_2(b) + \half\varepsilon_{2}F_1(a)
   + \fourth\varepsilon_{1}\varepsilon_{2} \\
   D(a^+,b^+) &=& D(a,b) + \threefourth\varepsilon_{12} - \half\varepsilon_{1}F_2(b) - \half\varepsilon_{2}F_1(a)
   - \fourth\varepsilon_{1}\varepsilon_{2}
\end{eqnarray*}
Now if $D(a,b)\ne 0$ these four expressions cannot all be zero, so
there must be an open set, in at least one of the four `quadrants'
and with $(a,b)$ as a limit point, where $D$ is nonzero.
Hence,~(\ref{intf12}) cannot be zero.

\newcommand{\aseq}{\stackrel{\mbox{\scriptsize a.s.}}{=}}

Part 3: By definition $\kappa(X,Y)^2\le\kappa(X,X)\kappa(Y,Y)$ is
equivalent to
\begin{eqnarray*}
   \left[Eh_{F_1}(X_1,X_2)h_{F_2}(Y_1,Y_2)\right]^2\le
   Eh_{F_1}(X_1,X_2)^2 Eh_{F_2}(Y_1,Y_2)^2
\end{eqnarray*}
This is a Cauchy-Schwartz inequality so it holds, and equality holds
iff
\begin{eqnarray}
   h_{F_1}(X_1,X_2) \aseq c\,h_{F_2}(Y_1,Y_2)\label{aseq1}
\end{eqnarray}
for some constant $c$. If $Y\aseq aX+b$ for certain constants $a$
and $b$ then it is immediately verified that~(\ref{aseq1}) holds
with $c=|a|$.

The reverse implication that~(\ref{aseq1}) implies linearity remains
to be shown. Suppose that~(\ref{aseq1}) holds. Then
\begin{eqnarray*}
   \lefteqn{h_{F_1}(X_1,X_2) - h_{F_1}(X_1,X_3) - h_{F_1}(X_2,X_4) - h_{F_1}(X_3,X_4)
   \aseq} \\
   && c\left[h_{F_2}(Y_1,Y_2) - h_{F_2}(Y_1,Y_3) - h_{F_2}(Y_2,Y_4) -
   h_{F_2}(X_3,X_4)\right]
\end{eqnarray*}
which reduces to
\begin{eqnarray*}
   \lefteqn{|X_1-X_2| - |X_1-X_3| - |X_2-X_4| + |X_3-X_4|
   \aseq} \\
   && c\left(|Y_1-Y_2| - |Y_1-Y_3| - |Y_2-Y_4| + |Y_3-Y_4|\right)
\end{eqnarray*}
But this is equivalent to
\begin{eqnarray}
   X_{3:4} - X_{2:4}
   \aseq
   c\left(Y_{3:4} - Y_{2:4}\right)
   \label{osdif}
\end{eqnarray}
Now without loss of generality suppose $Y_{3:4}=c X_{3:4}+b$ and
$Y_{2:4}=c X_{2:4}+b'$ for some $b$ and $b'$. Substitution
into~(\ref{osdif}) then yields $b=b'$, so the second and third order
statistics for $X$ and $Y$ are linearly related. Now the
distribution function of the second order statistic for $X$ is
\begin{eqnarray*}
   F_{1;2:4}(x) = \int_{-\infty}^x F_1(t)[1-F_1(t)]^2dF_1(t)
\end{eqnarray*}
and for $Y$
\begin{eqnarray*}
   F_{2;2:4}(y) = \int_{-\infty}^y F_2(t)[1-F_2(t)]^2dF_2(t)
\end{eqnarray*}
It is now straightforward to show that the equation
$F_{2;2:4}(y)=F_{1;2:4}(cx+b)$ leads to $F_2(y)=F_1(cx+b)$, so $X$
and $Y$ are linearly related.

Part 4: Without loss of generality assume $X\in\{0,1\}$ and
$Y\in\{0,1\}$ with probability one. Then $h_{F_1}=u_{F_1}$ and
$h_{F_2}=u_{F_21}$ (both $u$ and $h$ defined in Section~\ref{sec1}),
so
$\kappa(X,Y)=Eh_{F_1}(X_1,X_2)h_{F_2}(Y_1,Y_2)=Eu_{F_1}(X_1,X_2)u_{F_2}(Y_1,Y_2)=\cov(X,Y)$.

Part 5: Since $\kappa(Z,Z)=Eh_F(Z_1,Z_2)^2$, this follows directly
from Lemma~\ref{lem3}
 \hspace*{\fill}$\Box$\endproof

We conclude this section by giving some representations of $\kappa$
in terms of $h_F$ and the (conditional) distribution functions of
$X$ and $Y$. Let
\begin{eqnarray*}
   F_{2|1}(y|x)=P(Y<y|X=x)+\half P(Y=y|X=x)
\end{eqnarray*}
be the conditional distribution function of $Y$ given $X=x$.
\begin{lemma}\label{lem4}
The following equalities hold for $\kappa$:
\begin{enumerate}
\item  $\kappa(X,Y) = \int[F_{12}(x,y)-F_1(x)F_2(y)]^2dx dy$

\item
$\kappa(X,Y)=Eh_{F_1}(X_1,X_2)\int[F_{2|1}(y|X_1)-F_2(y)][F_{2|1}(y|X_2)-F_2(y)]dy$


\end{enumerate}
\end{lemma}
\proof{Proof} Part 1: This follows from the proof of
Theorem~\ref{th1}, Part~\ref{th1p2}

Part 2: First note that $d_{xy}F_{12}(x,y)=d_xF_1(x)d_yF_{2|1}(y|x)$
which we write in shorthand $dF_{12}(x,y)=dF_1(x)dF_{2|1}(y|x)$.
Hence,
\begin{eqnarray*}
   \kappa(X,Y)
   &=& Eh_{F_1}(X_1,X_2)h_{F_2}(Y_1,Y_2)\\
   &=& \half\int h_{F_1}(x_1,x_2)\int[\gamma(y_1,y)-F_2(y)][\gamma(y_2,y)-F_2(y)]dydF_{12}(x_1,y_1)dF_{12}(x_2,y_2)\\
   &=& \int h_{F_1}(x_1,x_2)\int[\gamma(y_1,y)-F_2(y)]dF_{2|1}(y_1|x_1)[\gamma(y_2,y)-F_2(y)]dF_{2|1}(y_2|x_2)dydF_1(x_1)dF_1(x_2)\\
   &=& \int h_{F_1}(x_1,x_2)\int[F_{2|1}(y|x_1)-F_2(y)][F_{2|1}(y|x_2)-F_2(y)]dydF_1(x_1)dF_1(x_2)\\
   &=& Eh_{F_1}(X_1,X_2)\int[F_{2|1}(y|X_1)-F_2(y)][F_{2|1}(y|X_2)-F_2(y)]dy
\end{eqnarray*}

 \hspace*{\fill}$\Box$\endproof

Note the similarity of Part~1 of Lemma~\ref{lem4} with the formula
for the covariance given by \citeA{hoeffding40}:
\begin{eqnarray*}
   \cov(X,Y) =
   \int
   [F_{12}(x,y)-F_1(x)F_2(y)]
    dx dy
\end{eqnarray*}

\subsection{Orthogonal decomposition}\label{orthog}

Let us assume $h_{F_1}$ and $h_{F_2}$ are square integrable and have
the spectral decompositions
\begin{eqnarray}
   h_{F_1}(x_1,x_2) &=& \sum_{k=0}^{\infty}\lambda_kg_{1k}(x_1)g_{1k}(x_2)\label{spec1}\\
   h_{F_2}(y_1,y_2) &=&
   \sum_{k=0}^{\infty}\mu_kg_{2k}(y_1)g_{2k}(y_2)\label{spec2}
\end{eqnarray}
See Lemma~\ref{lem3} on how to check for square integrability. For
ease of notation, we write the correlations between marginal
eigenfunctions as
\begin{eqnarray*}
   \rho_{kl}(X,Y) = \rho[g_{1k}(X),g_{2l}(Y)]
\end{eqnarray*}
We now have the orthogonal decomposition given as follows:
\begin{theorem}\label{th2}
Suppose $h_{F_1}$ and $h_{F_2}$ are square integrable with
spectral decompositions as above. Then with convergence in mean
square,
\begin{eqnarray*}
   \kappa(X,Y) = {\sum_{k=1}^\infty\sum_{l=1}^\infty\lambda_k\,\mu_l}\,\rho_{kl}(X,Y)^2
\end{eqnarray*}
and
\begin{eqnarray*}
   \rho^*(X,Y) =
   \frac{1}{\sqrt{\sum\lambda_k^2}\,\sqrt{\sum\mu_l^2}}
   {\sum_{k=1}^\infty\sum_{l=1}^\infty\lambda_k\,\mu_l}\,\rho_{kl}(X,Y)^2
\end{eqnarray*}
\end{theorem}
\proof{Proof} Write
\begin{eqnarray*}
   \kappa^{(N,N)}(X,Y)
   &=&
   E\left[\sum_{k=1}^{N}\lambda_kg_{1k}(X_1)g_{1k}(X_2)\,\sum_{l=1}^{N}\mu_lg_{2l}(Y_1)g_{2l}(Y_2)
   \right]\\
   \kappa^{(.,N)}(X,Y)
   &=&
   E\left[h_{F_1}(X_1,X_2)\,\sum_{l=1}^{N}\mu_lg_{2l}(Y_1)g_{2l}(Y_2)
   \right]\\
   \kappa^{(N,.)}(X,Y)
   &=&
   E\left[\sum_{k=1}^{N}\lambda_kg_{1k}(X_1)g_{1k}(X_2)\,h_{F_2}(Y_1,Y_2)
   \right]
\end{eqnarray*}
Then straightforward algebra gives
\begin{eqnarray*}
   \kappa^{(N,N)}(X,Y)
   =
   {\sum_{k=1}^N\sum_{l=1}^N\lambda_k\,\mu_l}\,\rho[g_{1k}(X),g_{2l}(Y)]^2
\end{eqnarray*}
By the Cauchy-Schwartz inequality we obtain
\begin{eqnarray*}
   \lefteqn{
   \left(\kappa(X,Y) - \kappa^{(N,.)}(X,Y) - \kappa^{(.,N)}(X,Y) + \kappa^{(N,N)}(X,Y)\right)^2
   }
   \\
   &=&
   E
   \left(
   \left[h_{F_1}(X_1,X_2)-
   \sum_{k=1}^{N}\lambda_kg_{1k}(X_1)g_{1k}(X_2)
   \right]
   \left[h_{F_2}(Y_1,Y_2)-
   \sum_{l=1}^{N}\mu_kg_{2l}(Y_1)g_{2l}(Y_2)
   \right]
   \right)^2
   \\
   &\le&
   E
   \left[h_{F_1}(X_1,X_2)-
   \sum_{k=1}^{N}\lambda_kg_{1k}(X_1)g_{1k}(X_2)
   \right]^2
   E
   \left[h_{F_2}(Y_1,Y_2)-
   \sum_{l=1}^{N}\mu_kg_{2l}(Y_1)g_{2l}(Y_2)
   \right]^2
\end{eqnarray*}
By mean square convergence of the spectral decomposition the latter
goes to zero as $N\rightarrow\infty$ so
\begin{eqnarray*}
   \kappa(X,Y) - \kappa^{(N,.)}(X,Y) - \kappa^{(.,N)}(X,Y) + \kappa^{(N,N)}(X,Y) \rightarrow 0
\end{eqnarray*}
as $N\rightarrow\infty$. Similarly we find
\begin{eqnarray*}
   \kappa(X,Y) - \kappa^{(N,.)}(X,Y) \rightarrow 0\\
   \kappa(X,Y) - \kappa^{(.,N)}(X,Y) \rightarrow 0
\end{eqnarray*}
as $N\rightarrow\infty$. It follows that
\begin{eqnarray*}
   \kappa^{(N,N)} \rightarrow \kappa(X,Y)
\end{eqnarray*}
as $N\rightarrow\infty$, which is the desired result.
 \hspace*{\fill}$\Box$\endproof

The simplest example of a decomposition is if both variables are
dichotomous, say $P(X=0)=1-P(X=1)=p$ and $P(Y=0)=1-P(Y=1)=q$, we
obtain $\lambda_1=2p(1-p)$, $\mu_1=2q(1-q)$, and $\lambda_k=\mu_k=0$
for $k>1$ so that $\rho^*(X,Y) = \rho(X,Y)^2$, see
Theorem~\ref{th1}, Part~4. In this special case the decomposition
consists of just one component.

\subsection{Parameterization of the likelihood}\label{parlik}

Let $f_{12}$ be the joint density of $(X,Y)$ with corresponding
marginal densities $f_1$ and $f_2$. Since
\begin{eqnarray*}
   \rho_{kl} =
   \int
   \frac{f_{12}(x,y)}{f_1(x)f_2(y)}
   \,g_{1k}(x)g_{2l}(y)dF_1(x)dF_2(y)
\end{eqnarray*}
we can decompose the joint density as:
\begin{eqnarray}
   f_{12}(x,y)
   =
   f_1(x)f_2(y)
   \left(
   1+
   \sum_{i=1}^\infty
   \sum_{j=1}^\infty
   g_{1k}(x)g_{2l}(y)\rho_{kl}
   \right)
   \label{mod1}
\end{eqnarray}
A similar equation can be given for discrete distributions, and a
general treatment can be given using the Radon-Nikodym derivative.

Decomposition~(\ref{mod1}) may be compared to the well-known
canonical correlation decomposition
\begin{eqnarray*}
   f_{12}(x,y)
   =
   f_1(x)f_2(y)
   \left(
   1+
   \sum_{k=1}^\infty
   a_{1k}(x)a_{2k}(y)\rho_{k}
   \right)
\end{eqnarray*}
Here, $a_{1k}$ and $a_{2k}$ are those functions maximizing the
correlation between $X$ and $Y$, subject to the restraint (for
$k>1$) that they are orthogonal to $a_{11},\ldots,a_{1,k-1}$ and
$a_{21},\ldots,a_{2,k-1}$, respectively, and $\rho_k$ is the
correlation between $a_{1k}(X)$ and $b_{2k}(Y)$.

\subsection{Fr{\'e}chet bounds for component
correlations}\label{frechet}

Below, we discuss the  interpretation and properties of the
component correlations. In particular, we look at bounds for the
component correlations.

For two random variables $X$ and $Y$ with joint distribution
function $F_{12}$ and marginal distribution functions $F_1$ and
$F_2$, the well-known Fr{\'e}chet upper bound $F_{12}^+$ is defined
by
\begin{eqnarray*}
   F_{12}^+(x,y) =\min\{F_1(x),F_2(y)\}
\end{eqnarray*}
and the Fr{\'e}chet lower bound $F_{12}^-$ is defined by
\begin{eqnarray*}
   F_{12}^-(x,y) = \max\{0,1-F_1(x)-F_2(y)\}
\end{eqnarray*}
Then $\rho(X,Y)=1$ if and only if $F_{12}=F_{12}^+$ and
$\rho(X,Y)=-1$ if and only if $F_{12}=F_{12}^-$. A more general
question is, for functions $g$ and $h$, for which $F_{12}$ the
correlation between $g(X)$ and $h(Y)$ is maximal or minimal. Let
\begin{eqnarray*}
   S^+_{g,h}
   &=&
   \left\{
   (x,y)\in\bR^2|g(x)=h(y)
   \right\}
\end{eqnarray*}
and
\begin{eqnarray*}
   S^-_{g,h}
   &=&
   \left\{
   (x,y)\in\bR^2|g(x)=- h(y)
   \right\}
\end{eqnarray*}
Then we have:
\begin{lemma}
For functions $g$ and $h$, $\rho[g(X),h(Y)]=1$ iff the support of
the distribution of $(X,Y)$ is a subset of $S^+_{g,h}$ and and
$\rho[g(X),h(Y)]=-1$ iff the support is a subset of $S^-_{g,h}$.
\end{lemma}
\proof{Proof} Suppose for simplicity that $g$ and $h$ are
standardized. Then $\rho[g(X),h(Y)]=1$ iff $P[g(X)=h(Y)]=1$ and
$\rho[g(X),h(Y)]=-1$ iff $P[g(X)=-h(Y)]=1$, and the lemma
immediately follows.
 \hspace*{\fill}$\Box$\endproof

For the dichotomous case we obtain the following:
\begin{example}
If both variables are dichotomous, say $P(X=0)=1-P(X=1)=p$ and
$P(Y=0)=1-P(Y=1)=q$, we obtain $S_{11}^+=\{(0,0),(1,1)\}$ and
$S_{11}^-=\{(0,1),(1,0)\}$.
\end{example}

Note that the bounds need not be attainable since it may be the case
that, for example, for certain $x$, there is no $y$ such that
$(x,y)\in S^+_{g,h}$.

Here we are interested in the bounds for the component correlations
$\rho_{ij}$ of $\rho^*$. For simplicity, we write
$S^\pm_{ij}=S^\pm_{g_i,g_j}$. The next example shows that if both
$X$ and $Y$ have uniform distributions on $[0,1]$, then the
component correlations of $\rho^*$ can attain the bounds $1$ and
$-1$:
\begin{example}
Suppose $X$ and $Y$ are uniformly distributed on $[0,1]$. Then the
eigenfunctions are the Fourier cosine functions (see
Table~\ref{eigtable}). The set $S^+_{kl}$ is formed by the solutions
$(x,y)$ to the equation
\begin{eqnarray*}
   \cos (k\pi x) = \cos (l\pi y)
\end{eqnarray*}
and $S^-_{kl}$ by the solutions of
\begin{eqnarray*}
   \cos (k\pi x) = -\cos (l\pi y)
\end{eqnarray*}
The solutions are plotted in Figure~\ref{univ fr}, for
$k=1,\ldots,4$ and $l=1,\ldots,4$. The bounds for the $\rho_{kl}$
are attainable since $\rho_{kl}(X,Y)=1$ for the uniform distribution
on $S^+_{kl}$ and $\rho_{kl}(X,Y)=-1$ for the uniform distribution
on $S^-_{kl}$.

\begin{figure}[htbp]
\begin{center}
 \subfigure[Positive bounds $S^+_{ij}$: $\rho_{ij}=1$ if support of $(X,Y)$ is subset of $S^+_{ij}$]
 {\includegraphics[width=70mm]{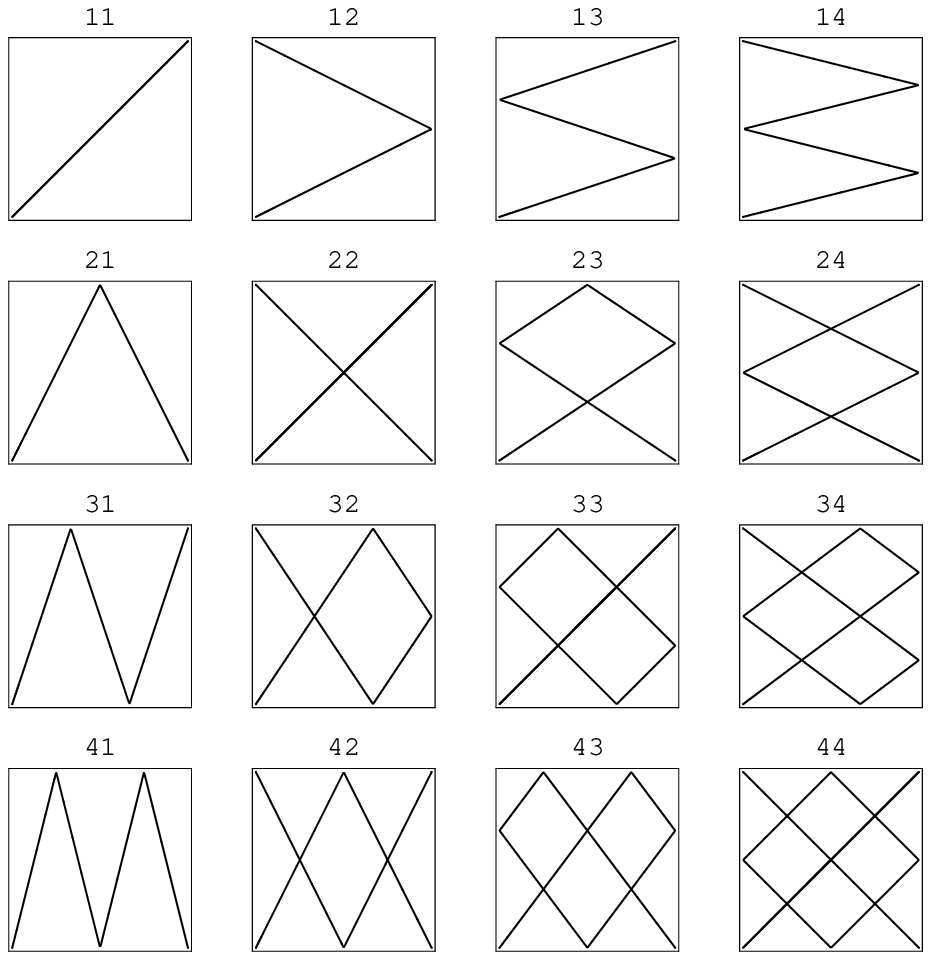}}
\end{center}
\begin{center}
 \subfigure[Negative bounds $S^-_{ij}$: $\rho_{ij}=1$ if support of $(X,Y)$ is subset of $S^-_{ij}$]
 {\includegraphics[width=70mm]{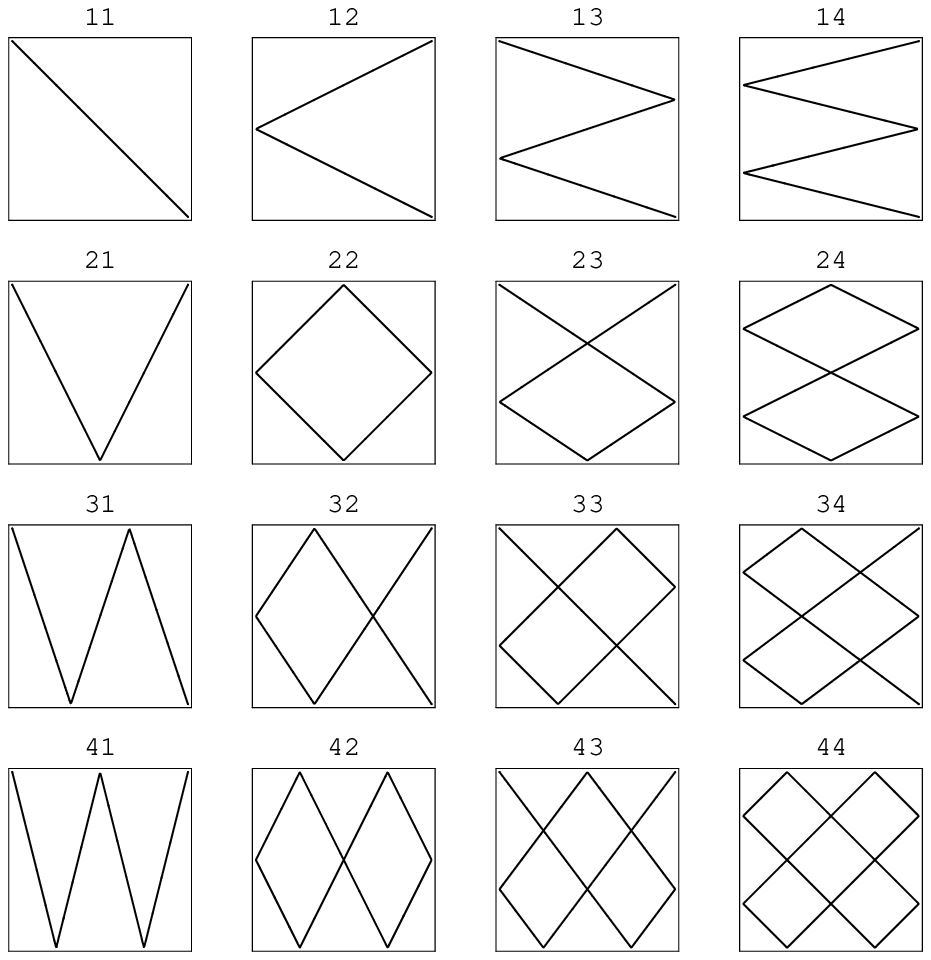}}
\end{center}
\caption{Supports of the Fr{\'e}chet bounds for the component
correlations $\rho_{ij}$ of $\rho^*(X,Y)$ when $X$ and $Y$ are
uniformly distributed on an interval.}\label{univ fr}
\end{figure}

Note that, if $\rho_{11}=1$, then $\rho_{22}=1$, since $S_{11}^+
\subset S_{22}^+$. More generally, by the same reasoning, we have
for all $i>1$,
\begin{eqnarray*}
   \rho_{kl}(X,Y)=1 \Rightarrow
      \rho_{i\times k,i\times l}(X,Y)=1
\end{eqnarray*}
and
\begin{eqnarray*}
   \rho_{kl}(X,Y)=-1 \Rightarrow
      \rho_{i\times k,i\times l}(X,Y)=(-1)^i
\end{eqnarray*}
By a symmetry argument, we also have
$\rho_{11}=1\Rightarrow\rho_{12}=0$, and there are various other
similar implications.
\end{example}

An overview of a large amount of literature on Fr{\'e}chet bounds is
given in \citeA{ruschendorf91}

\section{Estimation and tests of independence}\label{sec4}

In this section we discuss estimation of $\kappa$ and $\rho^*$ by U-
and V-statistics. Roughly speaking, the U-statistic estimator of a
parameter is an unbiased estimator based on taking averages
\cite{hoeffding48}, and the V-statistic estimator is the estimator
based on the distribution obtained by assigning a probability weight
$1/n$ to each sample point. For $\kappa$, both the U- and
V-statistic estimators are available, but for $\rho^*$ only the
latter is. However, we can estimate $\rho^*$ by a function of
U-statistic estimators.

In Section~\ref{uv}, it is shown how estimators of $\kappa$ and
$\rho^*$ by U- and V-statistics are obtained. In Section~\ref{perm}
permutation tests, useful for small samples, are described. In
Section~\ref{asymp}, the asymptotic distribution of these estimators
is derived under the null hypothesis of independence. In
Section~\ref{bonf}, Bonferroni corrections for tests of significance
of the component correlations are described.


\subsection{U and V statistic estimators of $\kappa$}\label{uv}

We first give a method for calculating the U- and V-statistic
estimators of $\kappa$ based on a sample
$(X_1,Y_1),\ldots,(X_n,Y_n)$, then we give the related estimates for
$\rho^*$.

The V-statistic estimator $\hat\kappa$ is the value of $\kappa$
based on the sample distribution functions $\hat F_1$ and $\hat
F_2$, and is obtained as follows. Let
\begin{eqnarray*}
   A_{1k} &=& \frac{1}{n}\sum_{i=1}^n|X_k-X_i|\\
   A_{2k} &=& \frac{1}{n}\sum_{i=1}^n|Y_k-Y_i|
\end{eqnarray*}
and
\begin{eqnarray*}
   B_{1} &=& \frac{1}{n^2}\sum_{i=1}^n\sum_{j=1}^n|X_i-X_j|\\
   B_{2} &=& \frac{1}{n^2}\sum_{i=1}^n\sum_{j=1}^n|Y_i-Y_j|
\end{eqnarray*}
Then we have for $k,l=1,\ldots,n$,
\begin{eqnarray*}
   h_{\hat F_1}(x_1,x_2) &=& -\half
   \left(|x_1-x_2| - A_{1k} - A_{1l} + B_1\right)\\
   h_{\hat F_2}(y_1,y_2) &=& -\half
   \left(|y_1-y_2| - A_{2k} - A_{2l} + B_2\right)\\
\end{eqnarray*}
and the sample or V-statistic estimator of $\kappa$ is given as
\begin{eqnarray*}
   \hat\kappa = \frac{1}{n^2}\sum_{i,j=1}^nh_{\hat F_1}(X_i,X_j)h_{\hat F_2}(Y_i,Y_j)
\end{eqnarray*}
Now with
\begin{eqnarray*}
   \tilde h_{\hat F_1}(x_1,x_2) &=& -\half
   \left(|x_1-x_2| - \frac{n}{n-1}A_{1k} - \frac{n}{n-1}A_{1l} + \frac{n}{n-1}B_1\right)\\
   \tilde h_{\hat F_2}(y_1,y_2) &=& -\half
   \left(|y_1-y_2| - \frac{n}{n-1}A_{2k} - \frac{n}{n-1}A_{2l} + \frac{n}{n-1}B_2\right)\\
\end{eqnarray*}
for $k,l=1,\ldots,n$, the unbiased or U-statistic estimator of
$\kappa$ is given as
\begin{eqnarray*}
   \tilde\kappa = \frac{2}{n(n-1)}\sum_{i=1}^{n-1}\sum_{j=i+1}^n\tilde h_{\hat F_1}(X_i,X_j)\tilde h_{\hat F_2}(Y_i,Y_j)
\end{eqnarray*}
By Hoeffding's theory of U-statistics we have that $\tilde\kappa$ is
an unbiased estimator of $\kappa$ \cite{hoeffding48,rw79}. Note that
$\hat\kappa\ge 0$ but $\tilde \kappa$ may be negative.

The related estimators of $\rho^*$ are the following:
\begin{eqnarray*}
   \hat\rho^*(X,Y)   &=& \frac{\hat\kappa(X,Y)}{\sqrt{\hat\kappa(X,X)\hat\kappa(Y,Y)}}\\
   \tilde\rho^*(X,Y) &=& \frac{\tilde\kappa(X,Y)}{\sqrt{\tilde\kappa(X,X)\tilde\kappa(Y,Y)}}
\end{eqnarray*}

For both types of estimators, the computational complexity of the
above method is $O(n^2)$.

The marginal eigenvalues and functions can be computed numerically
from $h_{\hat F_1}$ and $h_{\hat F_2}$ or from  $\tilde h_{\hat
F_1}$ and $\tilde h_{\hat F_2}$. See also Section~\ref{contapprox}
for computational aspects.

\subsection{Permutation tests}\label{perm}

Under independence, the sample marginal distributions of $X$ and $Y$
are ancillary statistics for $\hat\rho^*$ and $\tilde\rho^*$, so by
Fisher's theory of fiducial inference we should condition on the
sample marginals when testing for independence using $\hat\rho^*$
and $\tilde\rho^*$. If $X\ind Y$, conditioning on the marginals
ensures that $\hat\rho^*$ and $\tilde\rho^*$ are distribution free,
and exact conditional $p$-values can be calculated using the
permutation method. Evaluating all permutations quickly becomes
computationally prohibitive even for moderately large sample sizes,
and we recommend using a set of random permutations. Note that
permutation tests may also be applied to the component correlations
$\hat\rho_{ij}$ and $\tilde\rho_{ij}$.

Permutation tests may be computationally intensive. Using
non-optimized software, our experience shows that (bootstrap)
permutation tests for up to a several hundred observations are
feasible: for $n=100$, the permutation test based on 1000 random
permutations took less than four minutes, and for $n=500$ it took a
bit more than one hour. Techniques for the fast {\em exact}
evaluation of permutation tests using generating functions are
described by, among others, \citeA{bps96} and \citeA{vdwdbvdl99},
but it is not clear whether these techniques extend to statistics
such as $\hat\rho^*$ which are not based on ranks.


For categorical data the permutation test is better known as an {\em
exact conditional test} (where the conditioning is, again, on the
marginal distributions), the Fisher exact test being the best known
example. There is a large body of literature on fast evaluation of
exact conditional $p$-values for contingency tables, for an overview
see \citeA{agresti92} and for more recent developments see
\citeA{fms96,ds98,bb99}.

If the permutation test is too computationally intensive, an
asymptotic test may be done using the results of the next section.

\subsection{Asymptotic distribution of estimators under independence}\label{asymp}

For the asymptotic distribution of the estimators we obtain the
following:
\begin{theorem}\label{th asymp}
Suppose $h_{F_1}$ and $h_{F_2}$ are square integrable with spectral
decompositions~(\ref{spec1}) and~(\ref{spec2}). Then if $X\ind Y$
and with $Z_{ij}$ iid standard normal variables, we obtain
\begin{eqnarray*}
   n\tilde\kappa(X,Y) \rightarrow_D
   \sum_{i,j=0}^\infty \lambda_i\mu_j(Z_{ij}^2-1)
\end{eqnarray*}
If additionally $h_{F_1}$ and $h_{F_2}$ are trace class, we obtain
\begin{eqnarray*}
   n\hat\kappa(X,Y) \rightarrow_D
   \sum_{i,j=0}^\infty \lambda_i\mu_jZ_{ij}^2
\end{eqnarray*}
\end{theorem}
\proof{Proof} By the \citeA{hoeffding61} decomposition we can write
with $R_n=O(n^{-3})$
\begin{eqnarray*}
   \tilde\kappa
   &=&
   \left(\begin{array}{c}n \\ 2\end{array}\right)^{-1}
   \sum_{1\le i<j\le n} h_{F_1}(x_i,x_j)h_{F_2}(y_i,y_j) + R_n \\
   &=&
   \left(\begin{array}{c}n \\ 2\end{array}\right)^{-1}
   \sum_{1\le i<j\le n}
   \left( \sum_{k=1}^\infty \lambda_kg_{1k}(x_i)g_{1k}(x_j) \right)
   \left( \sum_{l=1}^\infty \mu_l g_{2l}(y_i)g_{2l}(y_j) \right) + R_n \\
   &=&
   \left(\begin{array}{c}n \\ 2\end{array}\right)^{-1}
   \sum_{k,l=1}^\infty \lambda_k\mu_l
   \left[
   \left( \sum_{i=1}^n g_{1k}(x_i)g_{2l}(y_i) \right)^2
   -
   \sum_{i=1}^n g_{1k}(x_i)^2g_{2l}(y_i)^2
   \right] + R_n
\end{eqnarray*}
Since $n^{-1}\sum_{i=1}^n g_{1k}(x_i)^2g_{2l}(y_i)^2\rightarrow 1$
a.s., and $n^{-1}\left( \sum_{i=1}^n g_{1k}(x_i)g_{2l}(y_i)
\right)^2\rightarrow_D Z_{kl}^2$ we obtain using the Cram{\'e}r-Wold
device that
\begin{eqnarray*}
   n\tilde\kappa \rightarrow_D
   \sum_{i,j=0}^\infty \lambda_i\mu_j(Z_{ij}^2-1)
\end{eqnarray*}

The proof for $\hat\kappa$ is similar; we have
\begin{eqnarray*}
   \hat\kappa
   &=&
   \frac1{n^2}
   \sum_{i,j=1}^n h_{F_1}(x_i,x_j)h_{F_2}(y_i,y_j) + R_n \\
   &=&
   \frac1{n^2}
   \sum_{i,j=1}^n
   \left( \sum_{k=1}^\infty \lambda_kg_{1k}(x_i)g_{1k}(x_j) \right)
   \left( \sum_{l=1}^\infty \mu_l g_{2l}(y_i)g_{2l}(y_j) \right) + R_n \\
   &=&
   \frac1{n^2}
   \sum_{k,l=1}^\infty \lambda_k\mu_l
   \left( \sum_{i=1}^n g_{1k}(x_i)g_{2l}(y_i) \right)^2
    + R_n
\end{eqnarray*}
Since $n^{-1}\left( \sum_{i=1}^n g_{1k}(x_i)g_{2l}(y_i)
\right)^2\rightarrow_D Z_{kl}^2$ we obtain using the Cram{\'e}r-Wold
device that
\begin{eqnarray*}
   n\hat\kappa \rightarrow_D
   \sum_{i,j=0}^\infty \lambda_i\,\mu_j\,Z_{ij}^2
\end{eqnarray*}
under the condition that
\begin{eqnarray*}
   \lim_{n\rightarrow\infty}E(n\hat\kappa)=\sum_{i,j=0}^\infty \lambda_i\,\mu_j = \sum_{i=1}^\infty \lambda_i\sum_{j=0}^\infty\mu_j
\end{eqnarray*}
is finite. Now by Lemma~\ref{lem6}, the two factors on the right
hand side are finite iff $h_{F_1}$ and $h_{F_2}$ are trace class,
completing the proof.
 \hspace*{\fill}$\Box$\endproof
The proof is similar to an adaptation by \citeA{dewet87} of a proof
by \citeA{eagleson79}. See also \citeA{gregory77} and
\citeA{hall79}.

Note that by Lemma~\ref{lemtr}, $h_{F_1}$ and $h_{F_2}$ are trace
class iff $EX$ and $EY$ exist. As follows from the theorem and noted
earlier by \citeA{dewet87} for related statistics, the U-statistic
estimator has an asymptotic distribution in more cases than the
V-statistic estimator. For example, if the marginal distribution of
at least one of $X$ and $Y$ is the distribution of
Example~\ref{ex1}, $n\hat\kappa$ does not have an asymptotic
distribution but $n\tilde\kappa$ does have one.

\subsection{Bonferroni corrections for testing significance of component correlations}\label{bonf}

As well as testing the significance of $\hat\rho^*$ directly, we can
test for the significance of the empirical component correlations
$\hat\rho_{ij}$. We recommend using $\hat\rho^*$ rather than
$\tilde\rho^*$ for calculating component correlations, since
$\tilde\rho^*$ may be negative in which case no component
correlations with nonnegative eigenvalues exist.

The proof of Theorem~\ref{th asymp} suggests that the component
correlations $\hat\rho_{kl}$ are asymptotically normal and
independent. Since there are many component correlations, a
simultaneous test of their significance needs a Bonferroni
correction. The ordinary Bonferroni correction, i.e., multiplying
the exceedance probabilities by the number of tests done, which in
this case is $n^2$, would be unreasonable since the multiplication
factor increases rapidly with $n$. Instead we propose dividing the
exceedance probability for $\hat\rho_{ij}$ by
\begin{eqnarray*}
   \frac{\hat\lambda_i\hat\mu_j}{\sum_{i=1}^n\hat\lambda_i\sum_{j=1}^n\hat\mu_j}
\end{eqnarray*}
Note that these numbers will converge to zero in probability if at
least one of $h_{F_1}$ and $h_{F_2}$ is not of trace class, i.e., by
Lemma~\ref{lemtr}, if at least one of $EX$ or $EY$ does not exist,
in which case the correction may not be the most appropriate one.

The idea of looking at components of a test seems to have first
appeared in \citeA{dk72}, where components of the Cram{\'e}r-von
Mises test were investigated. This test is a special case of the
tests based on $\rho^*$ described above (see Section~\ref{cvm}).

Other related work is by \citeA{kl99}, who looked at correlations
between orthogonal functions of the marginal cumulative distribution
functions, in particular, the Legendre polynomials. This work is an
extension of the so-called smooth tests of fit of \citeA{neyman37}.
Rather than looking at all correlations between the orthogonal
functions, they considered just the first few, and developed a
selection method based on Schwartz's rule for determining how many
correlations to base the overall test on.

\section{Grade versions of $\kappa$ and $\rho^*$, copulas, and rank
tests}\label{sec5}

For ordinal random variables $X$ and $Y$, any given scale is
arbitrary and it may be desirable to use scales based on the grades
$F_1(X)$ and $F_2(Y)$ of $X$ and $Y$. A general way to base $\kappa$
and $\rho^*$ on grades is as follows. For given invertible
distribution functions $K_1$ and $K_2$, we can define
\begin{eqnarray*}
   \kappa_{K_1,K_2}(X,Y) =
   \kappa[K_1^{-1}\circ F_1(X),K_2^{-1}\circ F_2(Y)]
\end{eqnarray*}
and
\begin{eqnarray*}
   \rho^*_{K_1,K_2}(X,Y) =
   \rho^*[K_1^{-1}\circ F_1(X),K_2^{-1}\circ F_2(Y)]
\end{eqnarray*}
Note that
\begin{eqnarray*}
   \kappa_{F_1,F_2}(X,Y) = \kappa(X,Y)
\end{eqnarray*}
and
\begin{eqnarray*}
   \rho^*_{F_1,F_2}(X,Y) = \rho^*(X,Y)
\end{eqnarray*}
With $K_1$ and $K_2$ uniform distribution functions,
$\rho^*_{K_1,K_2}(X,Y)$ is to $\rho^*$ what Spearman's rho is to the
ordinary correlation $\rho$.

We can use the results of Section~4 to obtain an orthogonal
decomposition of $\kappa_{K_1,K_2}$ and $\rho^*_{K_1,K_2}$ in terms
of component correlations. These component correlations then
parameterize the {\em copula}, which is defined as the joint
distribution of $(F_1(X),F_2(Y))$. From~(\ref{mod1}), and since the
marginal eigenfunctions of $h_F$ with $F$ the uniform distribution
are the cosine functions given in Table~\ref{eigtable}, we obtain
the following decomposition of $c_{12}$, the density function of the
copula:
\begin{eqnarray*}
   c_{12}(u,v)
   =
   1+
   \sum_{i=1}^\infty
   \sum_{j=1}^\infty
   \cos(k\pi u)\cos(l\pi v)\rho_{kl}
\end{eqnarray*}
where $\rho_{kl}=\int \cos(k\pi u)\cos(l\pi v)c_{12}(u,v)dudv$. This
decomposition was earlier given in \citeA{dewet80} and
\citeA{deheuvels81}. An overview of copula theory is given in
\citeA{nelsen06}. Possible drawbacks of using $\rho^*_{K_1,K_2}$ for
some given $K_1$ and $K_2$ is the arbitrariness of any choice of
$K_1$ and $K_2$ and the loss of scale information, but these issues
are hotly debated \cite{mikosch06}.

In Section~\ref{secrank}, a brief description of rank tests based on
$\kappa$ is given. In Section~\ref{cvm} a generalization of the
Cram{\'e}r-von Mises test to the case of $K$ ordered samples is
shown to be a special case, and a convenient representation is
given. In Section~\ref{secphi} we write $\kappa_{K_1,K_2}$ as a
weighted average of $\phi$-coefficients.

\subsection{Rank tests}\label{secrank}

Rank statistics which are distribution free under independence in
the continuous case are obtained as follows. For invertible
distribution functions $K_1$ and $K_2$ let
\begin{eqnarray*}
   \hat\rho^*_{K_1,K_2}(X,Y)   &=& \hat\rho^*[K_1^{-1}\circ \hat F_1(X),K_2^{-1}\circ \hat F_2(Y)] \\
   \tilde\rho^*_{K_1,K_2}(X,Y) &=& \tilde\rho^*[K_1^{-1}\circ \hat F_1(X),K_2^{-1}\circ \hat F_2(Y)]
\end{eqnarray*}
The derivation of the asymptotic distribution of these statistics is
slightly more involved than that of the asymptotic distribution of
$\hat\rho^*(X,Y)$. De Wet \citeyear{dewet80} has done this
derivation for statistics related to $\tilde\rho^*_{K_1,K_2}(X,Y)$.
He gave the weights for optimal tests in the Bahadur sense for
certain classes of alternatives, such as the bivariate normal.

With $K_1$ and $K_2$ the uniform distribution functions,
$n\hat\rho^*_{K_1,K_2}(X,Y)$ is a statistic discussed by
\citeA{bkr61}, see also \citeA{deheuvels81}. It can be viewed as a
generalization of the ordinary Cram{\'e}r-von Mises test (see next
subsection). Similarly, with $K_1$ and $K_2$ the logistic
distribution functions, $n\hat\rho_{K_1,K_2}^*(X,Y)$ can be viewed
as a generalization of the Anderson-Darling test.

\citeA{hoeffding48ind} described a related test, namely based on the
U-statistic estimator of
\begin{eqnarray*}
   \int \left[F_{12}(x,y)-F_1(x)F_2(y)\right]^2dF_{12}(x,y)
\end{eqnarray*}
which can be obtained from the representation of $\kappa$ in
Lemma~\ref{lem4}, Part~1, by replacing $dxdy$ by $dF_{12}(x,y)$.
Hoeffding's coefficient does not fall in the framework of the
present paper.

\subsection{A new class of $K$-sample Cram{\'e}r-von Mises tests as a special case}\label{cvm}

Suppose we have $K$ samples, the $k$th sample having $n_k$ iid
observations, say $\{U_{k1},\ldots,U_{kn_k}\}$. Then a test whether
the distributions of the observations in the different samples are
equal is called a {\em $K$-sample test}. A $K$-sample test can, in
fact, be viewed as a test of independence, namely, whether
`response' depends on `group membership,' the groups referring to
the different samples. Let us consider the case that the score
$c_k\in\bR$ is assigned to sample $k$ ($k=1,\ldots,K$) . With
$N_0=0$ and $N_k=\sum_{i=1}^kn_i$ let
$(X_{N_{k-1}+i_k},Y_{N_{k-1}+i_k})=(c_i,U_{k,i_k})$ for
$k=1,\ldots,K$ and $i_k=1,\ldots,n_k$. Then it can be seen that the
$K$ sample test is a test of independence of the $X$ observations
and the $Y$ observations. (Note that here the $X$ observations are
not random). A $K$-sample test can then be based on $\hat\rho^*$ or
$\tilde\rho^*$.

If samples are ordered but have no numerical scores assigned to
them, rank scores can be assigned, for example $c_k=N_k$.

In order to arrive at the Cram{\'e}r-von Mises test, we now use
Lemma~\ref{lem4}, Part~2 to give a representation of $\kappa$ in
terms of the conditional distribution functions. Let $G_k$ be the
distribution function of $U_k$, the response for sample $k$, and let
$p_k=n_k/N_K$ be the proportion of observations in sample $k$. Then
we obtain
\begin{eqnarray*}
   \kappa(X,Y) =\sum_{i,j}
   p_ip_jh_F(c_i,c_j)\int\left[G_i(y)-F_2(y)\right]\left[G_j(y)-F_2(y)\right]dy
\end{eqnarray*}
Some straightforward algebra shows that for the two-sample case this
reduces to
\begin{eqnarray*}
   \kappa(X,Y) =
   p_1^2p_2^2\int\left[G_1(y)-G_2(y)\right]^2dy
\end{eqnarray*}
A grade version of $\kappa$ is
\begin{eqnarray*}
   \kappa_{F_1,K}(X,Y) =\sum_{i,j}
   p_ip_jh_F(c_i,c_j)\int\left[G_i(y)-F_2(y)\right]\left[G_j(y)-F_2(y)\right]w[F_2(y)]dF_2(y)
\end{eqnarray*}
where
\begin{eqnarray*}
   w(u) = \frac{1}{k[K^{-1}(u)]}
\end{eqnarray*}
In the two-sample case, the sample version of $\kappa_{F_1,K}(X,Y)$
with $K$ the uniform distribution function reduces (essentially) to
the ordinary Cram{\'e}r-von Mises statistic, so we have a
generalization to the case of $K$ ordered samples. With $K$ the
logistic distribution, $\hat\kappa_{F_1,K}(X,Y)$ reduces to the
Anderson-Darling statistic.


A different generalization of the two-sample Cram{\'e}r-von Mises
test was given by \citeA{kiefer59}, namely to the case of $K$ {\em
unordered} samples.

\subsection{$\kappa$ as a weighted $\phi$-coefficient}\label{secphi}

\begin{figure}
\begin{center}
\includegraphics[width=.98\linewidth,clip]{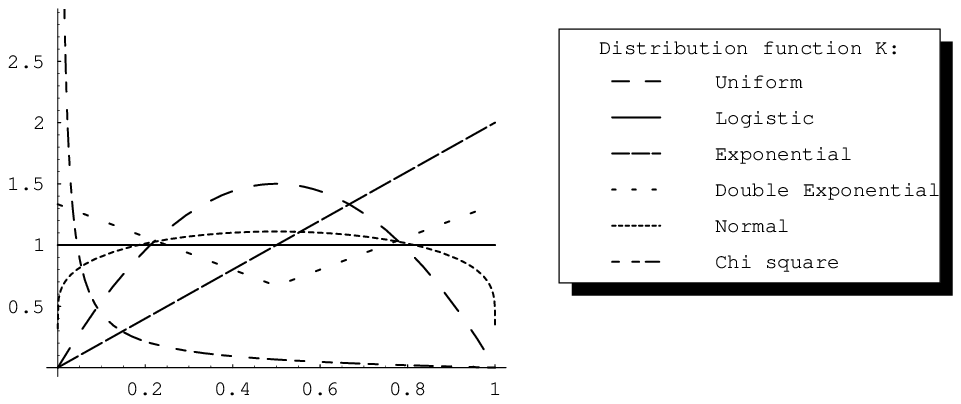}
\end{center}
{\caption{Plots of $\bar w_K(u)$ for several distributions
$K${\label{ffinv}}}}
\end{figure}

From Lemma~\ref{lem4}, Part~1, we directly obtain
\begin{eqnarray}
   \kappa_{K_1,K_2}(X,Y) =
   \int[F_{12}(x,y)-F_1(x)F_2(y)]^2\,dK_1^{-1}\circ
   F_1(x)\,dK_2^{-1}\circ F_2(y) \label{repr3}
\end{eqnarray}
This result leads to an interesting interpretation of
$\kappa_{K_1,K_2}$. The $\phi$ coefficient for measuring the
dependence in the $2\times 2$ table obtained by collapsing the
distribution with respect to the cut point $(x,y)$ is given as
\begin{eqnarray}
   \phi(x,y) = \frac{\left|F_{12}(x,y)-F_1(x)F_2(y)\right|}{\sqrt{F_1(x)[1-F_1(x)]F_2(y)[1-F_2(y)]}}
   \label{phidef1}
\end{eqnarray}
Now suppose $\psi$ is such that
\begin{eqnarray*}
   \phi(x,y) = \psi[F_1(x),F_2(y)]
\end{eqnarray*}
Then from~(\ref{repr3}) we obtain that $\kappa_{K_1,K_2}$ can be
written as a weighted average $\psi$-square:
\begin{eqnarray*}
   \kappa_{K_1,K_2}(X,Y) = \int\psi(u,v)^2\,w_{K_1}(u)\,w_{K_2}(v)dudv
\end{eqnarray*}
where the weight function $w$ is defined by
\begin{eqnarray*}
   w_K(u) = \frac{u(1-u)}{k[K^{-1}(u)]}
\end{eqnarray*}
The normalized weight function (integrating to one) is
\begin{eqnarray*}
  {\bar w}_K(u) = \frac{w_K(u)}{\int_0^1 w_K(u)du}
\end{eqnarray*}
where
\begin{eqnarray*}
   \int_0^1 w_K(u)du
   = \int_0^1{u(1-u)}dK^{-1}(u)
   = \int_{-\infty}^{\infty}K(x)[1-K(x)]dx
\end{eqnarray*}
In Figure~\ref{ffinv}, $\bar w_K$ is plotted for the distribution
functions $K$ given in Table~\ref{eigtable}.

From Figure~\ref{ffinv} we can deduce which $\phi$-coefficients are
given most weight for different marginal distributions. As the
reference marginal, we take the logistic, which is the horizontal
line in the figure, i.e., assigning uniform weights. We see that for
a uniform marginal, the weight goes to zero in the tails. The
weights for a normal marginal are for most $u$ intermediate between
the weights for the uniform and logistic marginals. In contrast to
uniform and normal marginals, the Laplace distribution gives more
weight to the tails than a logistic marginal. An exponential
marginal gives little weight to the lower tail, but much weight to
the upper tail. Finally, the chi-square distribution gives very
large weight to the lower tail and very small weight to the upper
tail. Among the distributions considered, the biggest difference is
between the chi-square and the exponential distribution.

\section{Data analysis: investigating the nature of the association}\label{sec6}

Gaining an understanding of the nature of an association between two
random variables is probably best viewed as an art rather than a
science, and in this section we present some visual tools based on
$\rho^*$ and its components which may be helpful in reaching this
objective. For an iid bivariate sample $\{(X_i,Y_i)\}$ we propose
the following two procedures.

Firstly, we calculate $\hat\rho^*$ from the sample and test its
significance. If found to be significant, then for each data point
$(X_i,Y_i)$ we calculate the weight
\begin{eqnarray*}
   W_{i} = \frac{\frac1n\sum_{j=1}^n h_{\hat F_1}(X_i,X_j)\,h_{\hat F_2}(Y_i,Y_j)}{\sqrt{\hat\kappa(X,X)\,\hat\kappa(Y,Y)}}
\end{eqnarray*}
Since
\begin{eqnarray*}
   \hat\rho^*(X,Y) = \frac{1}{n}\sum_{i=1}^nW_{i}
\end{eqnarray*}
the weights $W_{i}$ give an indication of how much the sample
element $(X_i,Y_i)$ contributes to $\hat\rho^*$, and so can be used
to discover the nature of a possible association between $X$ and
$Y$.

Secondly, we calculate the component correlations $\hat\rho_{kl}$ of
$\hat\rho^*$ and test their significance using the Bonferroni
corrections described in Section~\ref{bonf}. Then for each
significant component correlation $\hat\rho_{kl}$, we compute the
weights
\begin{eqnarray*}
   W_{i}^{(k,l)} = \hat g_{1k}(X_i)\hat g_{2l}(Y_i)
\end{eqnarray*}
where $g_{1k}$ and $g_{2l}$ are the eigenfunctions belonging to
$h_{F_1}$ and $h_{F_2}$. Since
\begin{eqnarray*}
   \hat\rho_{kl} = \frac{1}{n}\sum_{i=1}^nW_{i}^{(k,l)}
\end{eqnarray*}
the weight $W_{i}^{(k,l)}$ is the amount the sample element
$(X_i,Y_i)$ contributes to $\hat\rho_{kl}$ (conditionally on the
marginals), and so, like $W_i$, can be used to investigate the
association between $X$ and $Y$.

In this section we show how to visualize the weights $W_i$ and
$W_i^{(k,l)}$, both for continuous and categorical data, and show
how this can be used to gain an understanding of the association.
Some artificial continuous data sets are considered in
Section~\ref{sec art}, a real categorical data set is considered in
Section~\ref{sec cat} and a real time series data set is considered
in Section~\ref{sec time}

\subsection{Some artificial data sets}\label{sec art}

\begin{figure}[tbp]
 \subfigure[Increasing location ($\hat\rho^*=.36$)]{\includegraphics[width=60mm]{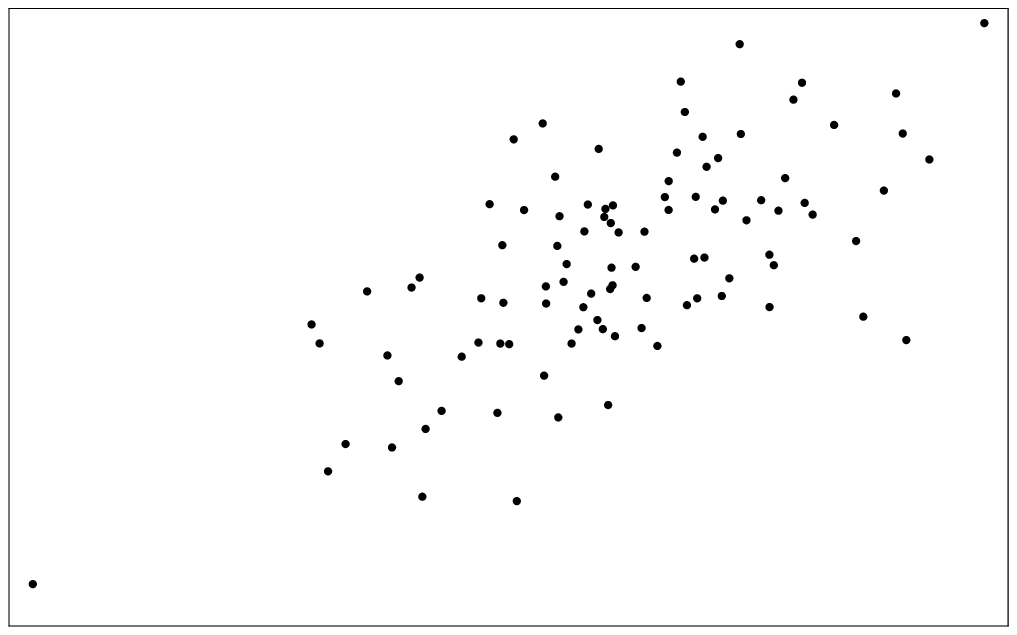}}
\hfill
 \subfigure[Decreasing then increasing location ($\hat\rho^*=.17$)]{\includegraphics[width=60mm]{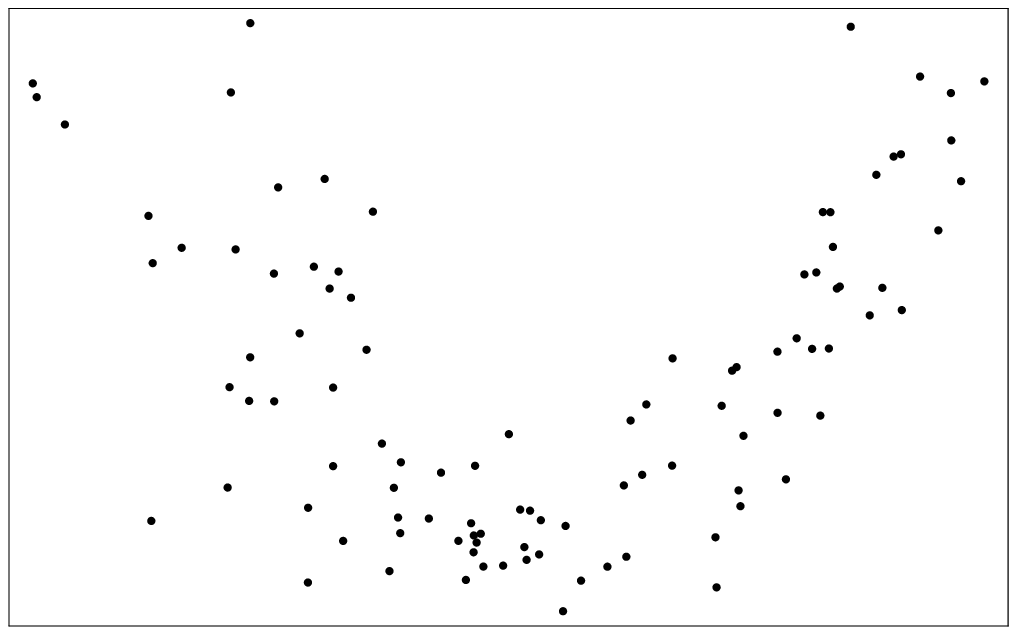}}
\\
 \subfigure[Increasing dispersion ($\hat\rho^*=.11$)]{\includegraphics[width=60mm]{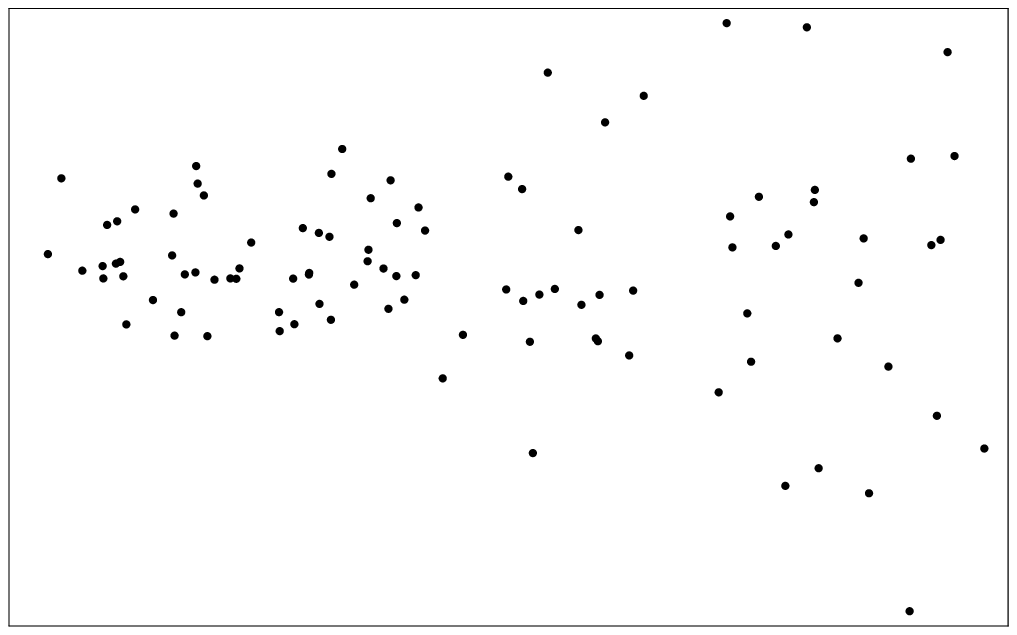}}
\hfill
 \subfigure[Increasing then decreasing dispersion ($\hat\rho^*=.02$)]{\includegraphics[width=60mm]{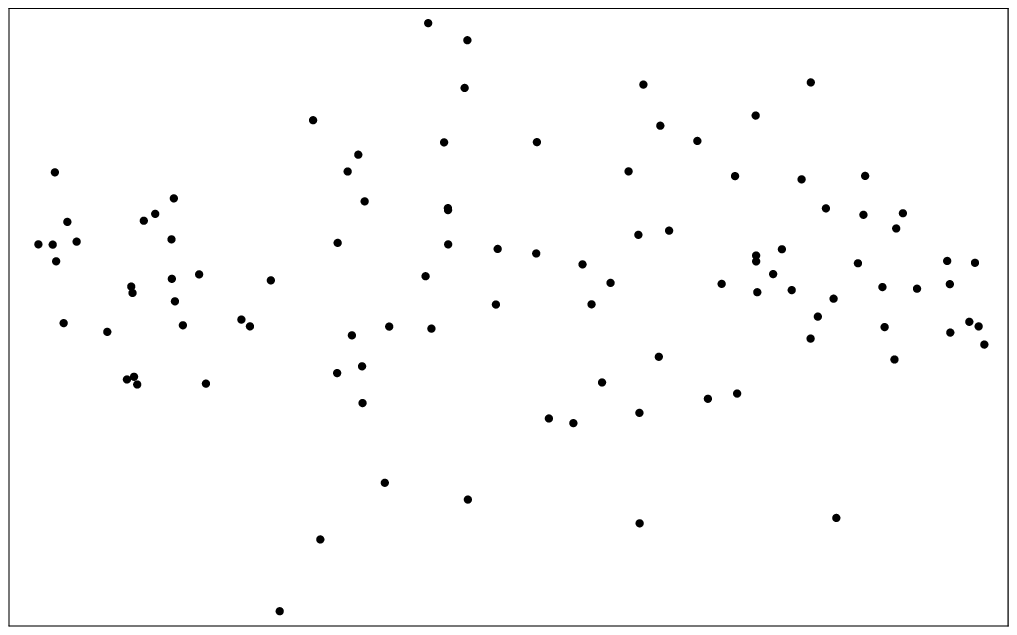}}
\caption{Scatterplots of artificial data sets. The captions denote
what happens to the conditional $Y$ distribution as $X$
increases.}\label{4scats}
\end{figure}

In Figure~\ref{4scats}, four artificial data sets are plotted, each
consisting of 100 iid points. For completeness, we explain how the
data were generated. In the following, $U$ is uniformly distributed
on $[0,1]$ and $Z_1$ and $Z_2$ are iid standard normal random
variables, and $Z(u)$ has a normal distribution with mean zero and
standard deviation $u$. The data in Figure~\ref{4scats}(a) are from
a bivariate normal distribution with $\rho=2/3$. The data in
Figure~\ref{4scats}(b) are of the form
$(X,Y)=(U,(U-1/2)^2)+(Z_1/10,Z_2/10)$. The data in
Figure~\ref{4scats}(c) are of the form $(X,Y)=(U,Z(1/5+U))$.
Finally, the data in Figure~\ref{4scats}(d) are of the form
$(X,Y)=(U,Z(1/5+\min\{U,1-U\}))$.

For all four data sets, we performed permutation tests for the
significance of $\hat\rho^*$ and its component correlations
$\hat\rho_{ij}$ based on 10,000 random permutations. This took us
about 48 minutes for each data set. For the significant component
correlations we took 1 million random permutations to get a more
accurate $p$-value, and this took about 110 seconds per component
correlation. We also computed the ordinary correlation, and, not
surprisingly, only for the data in Figure~\ref{4scats}(a) it is
significantly different from zero. There we found that
$\hat\rho=.44$ ($p=.000$).

\begin{figure}[tbp]
 \subfigure[Increasing location ($\hat\rho^*=.36$)]{\includegraphics[width=60mm]{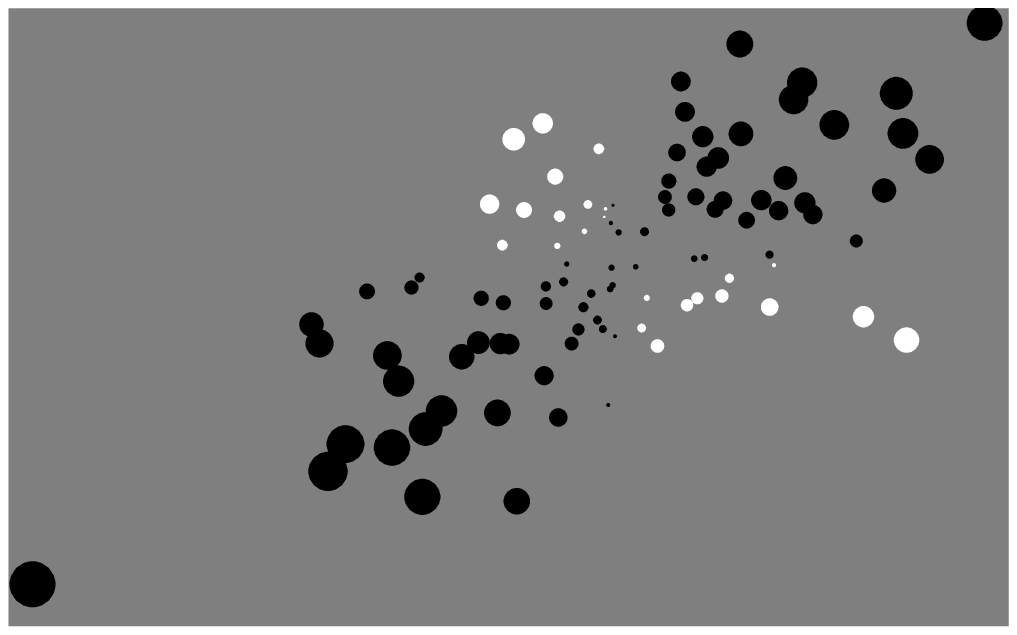}}
\hfill
 \subfigure[Decreasing then increasing location ($\hat\rho^*=.17$)]{\includegraphics[width=60mm]{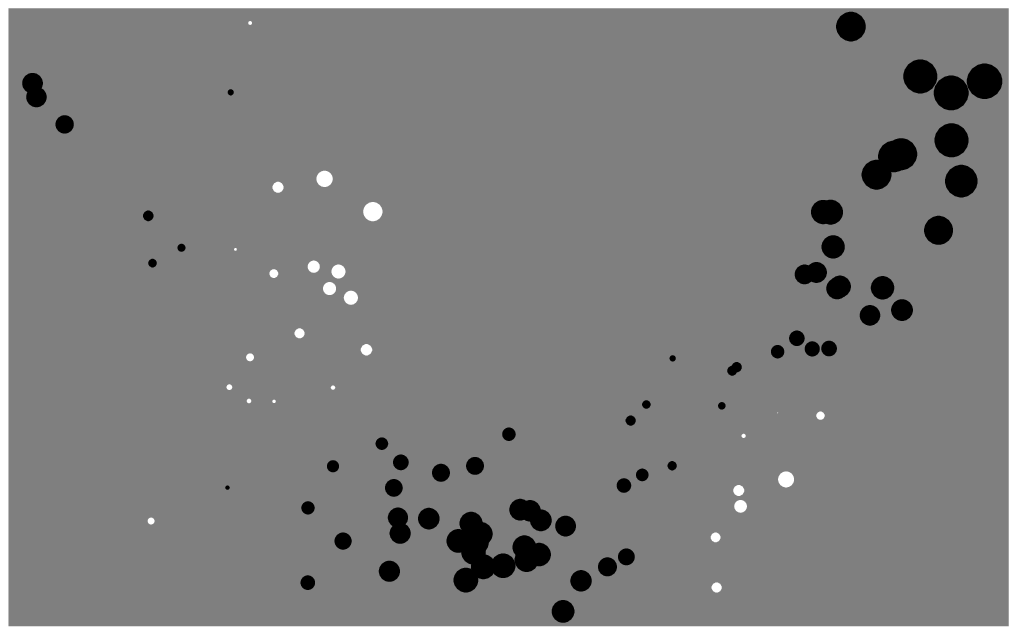}}
\\
 \subfigure[Increasing dispersion ($\hat\rho^*=.11$)]{\includegraphics[width=60mm]{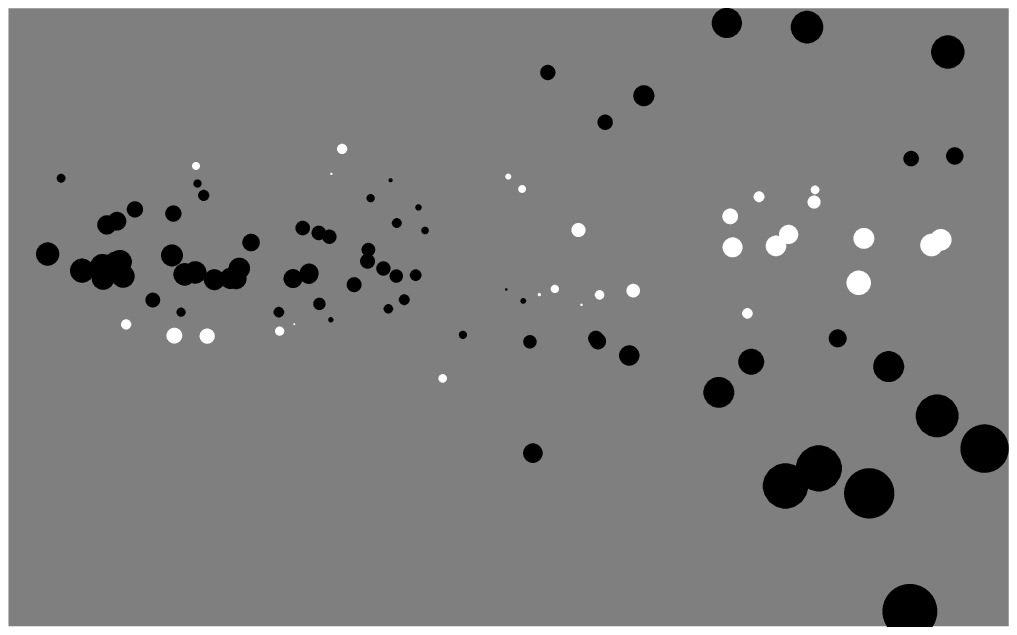}}
\hfill
 \subfigure[Increasing then decreasing dispersion ($\hat\rho^*=.02$)]{\includegraphics[width=60mm]{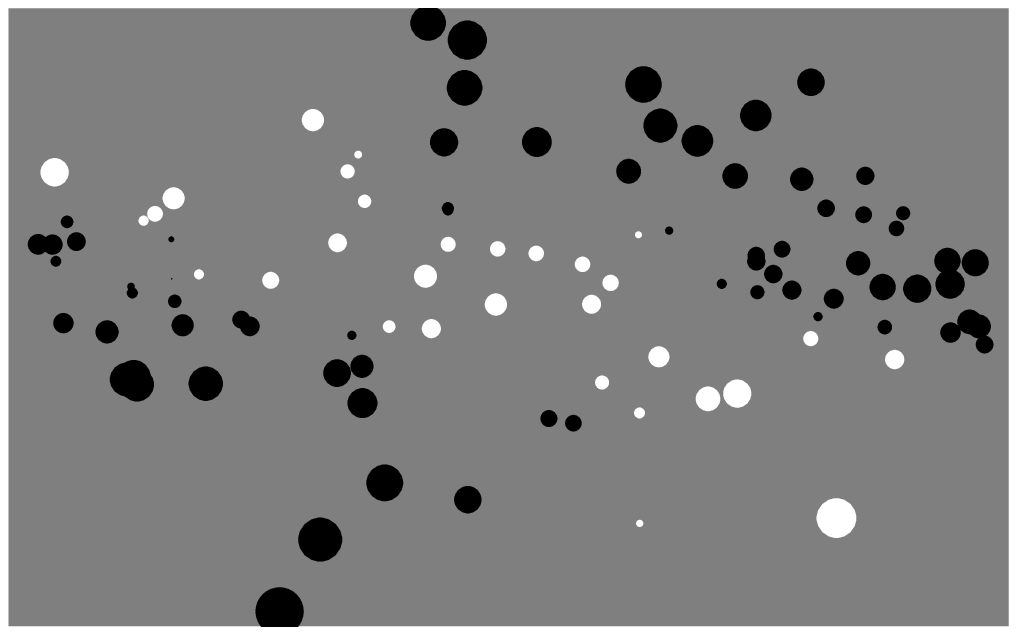}}
\caption{Representation of weights $W_i$ for data in
Figure~\ref{4scats}. The size of each dot is proportional to
$|W_i|$; black dots represent positive $W_i$ and white represent
negative $W_i$. Total area of dots is scaled to a fixed constant for
each plot.}\label{4scats2}
\end{figure}

\begin{figure}[tbp]
 \subfigure[$\hat\rho_{11}=.61$ (data from Figure~\ref{4scats}(a))]
 {\includegraphics[width=60mm]{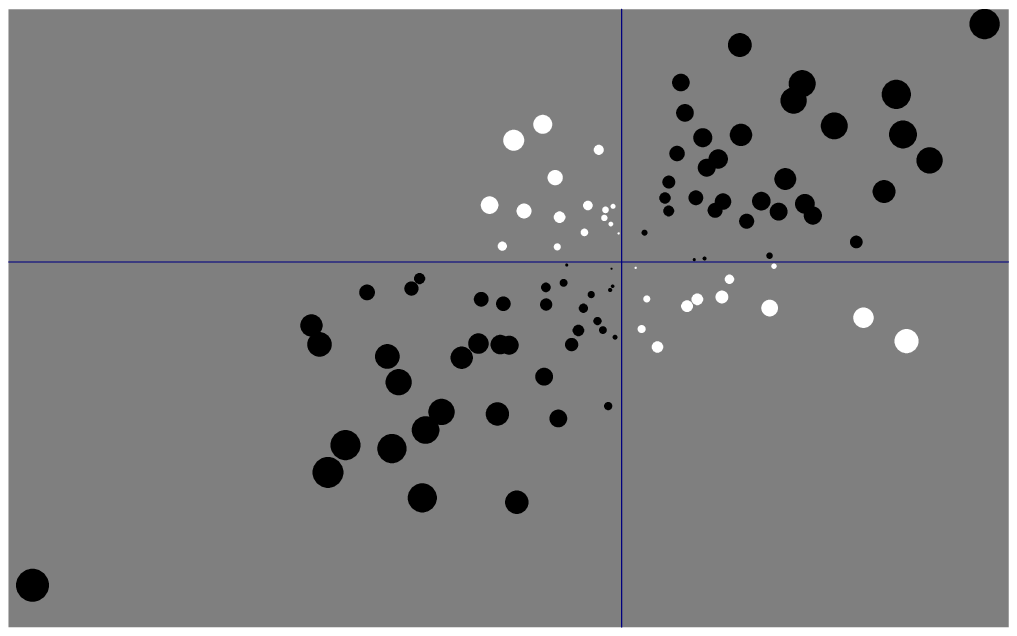}}
\hfill
 \subfigure[$\hat\rho_{22}=.38$  (data from Figure~\ref{4scats}(a))]
 {\includegraphics[width=60mm]{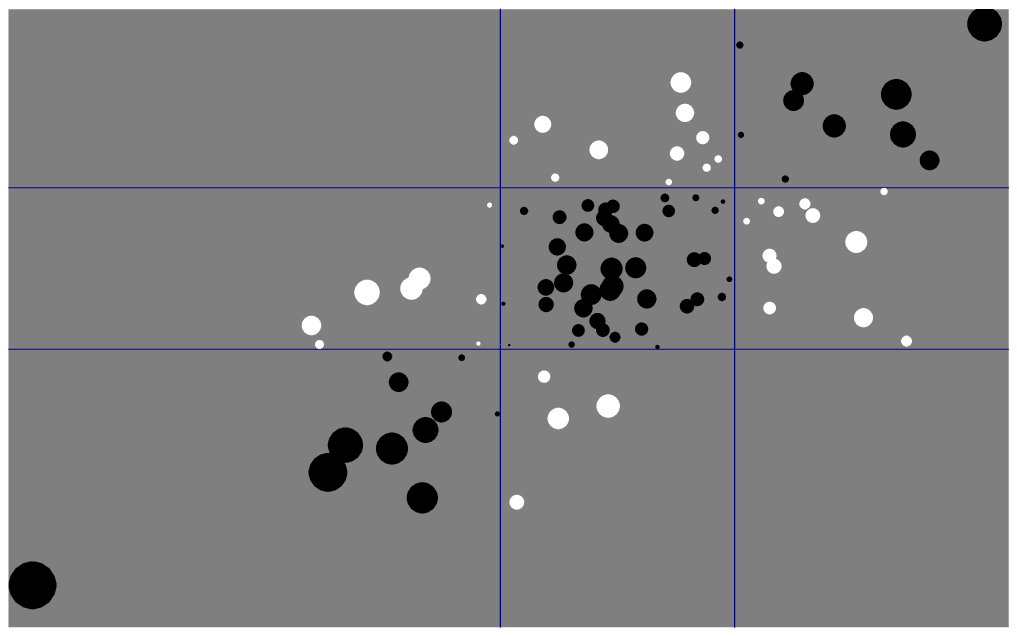}}
\\
 \subfigure[$\hat\rho_{21}=-.78$  (data from Figure~\ref{4scats}(b))]
 {\includegraphics[width=60mm]{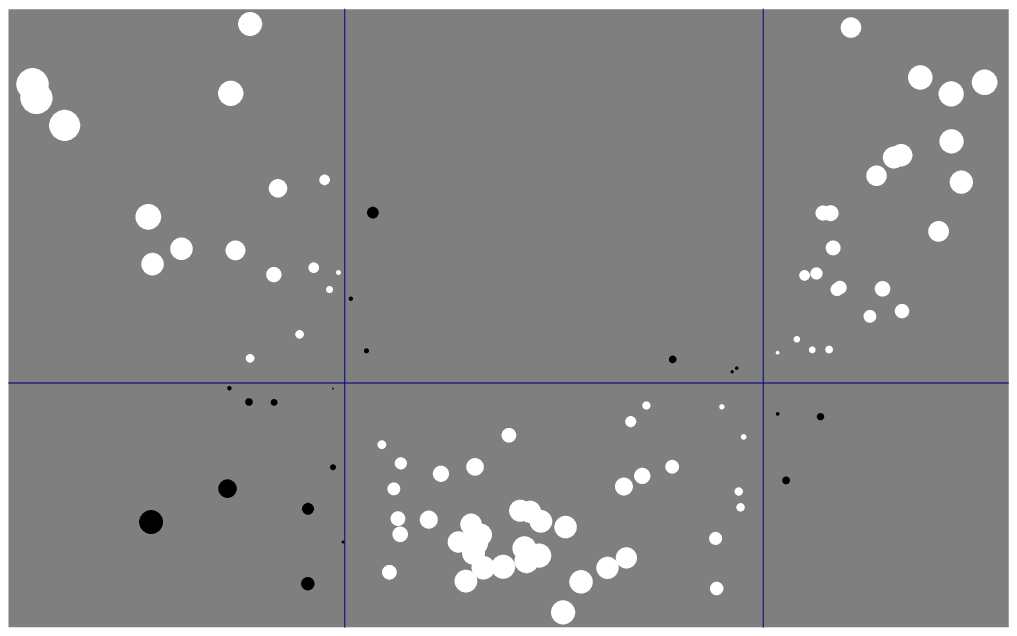}}
\hfill
 \subfigure[$\hat\rho_{42}=.54$  (data from Figure~\ref{4scats}(b))]
 {\includegraphics[width=60mm]{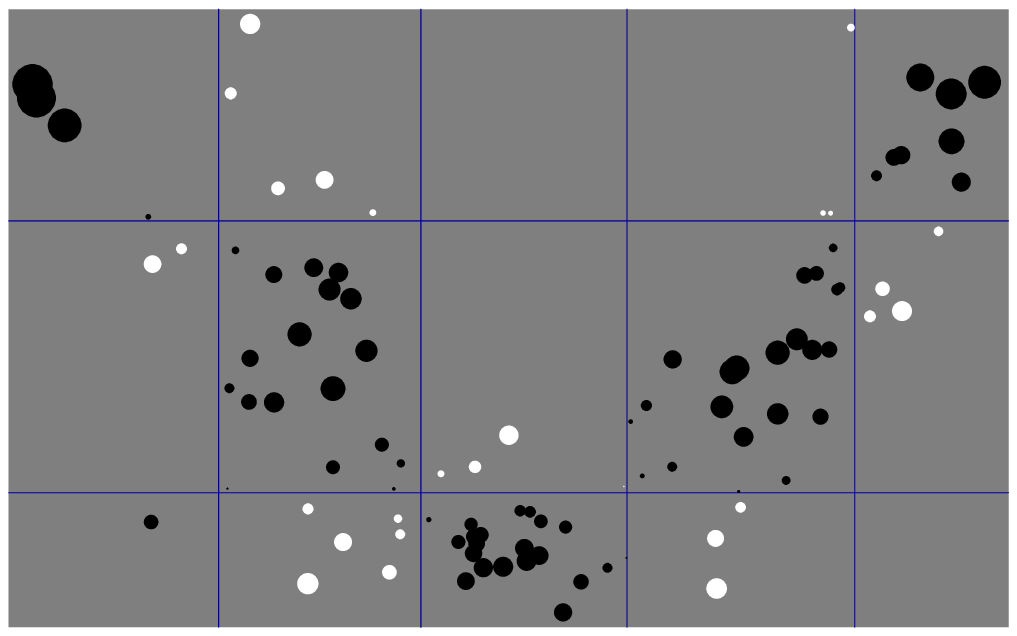}}
\\
 \subfigure[$\hat\rho_{12}=-.51$ (data from Figure~\ref{4scats}(c))]
 {\includegraphics[width=60mm]{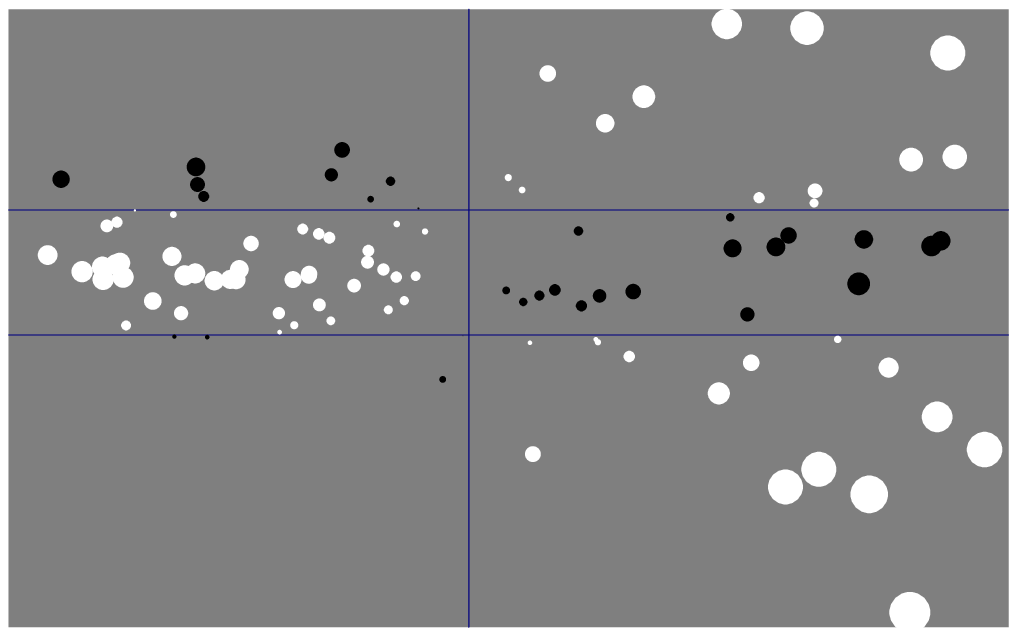}}
\hfill
 \subfigure[$\hat\rho_{22}=-.38$  (data from Figure~\ref{4scats}(d))]
 {\includegraphics[width=60mm]{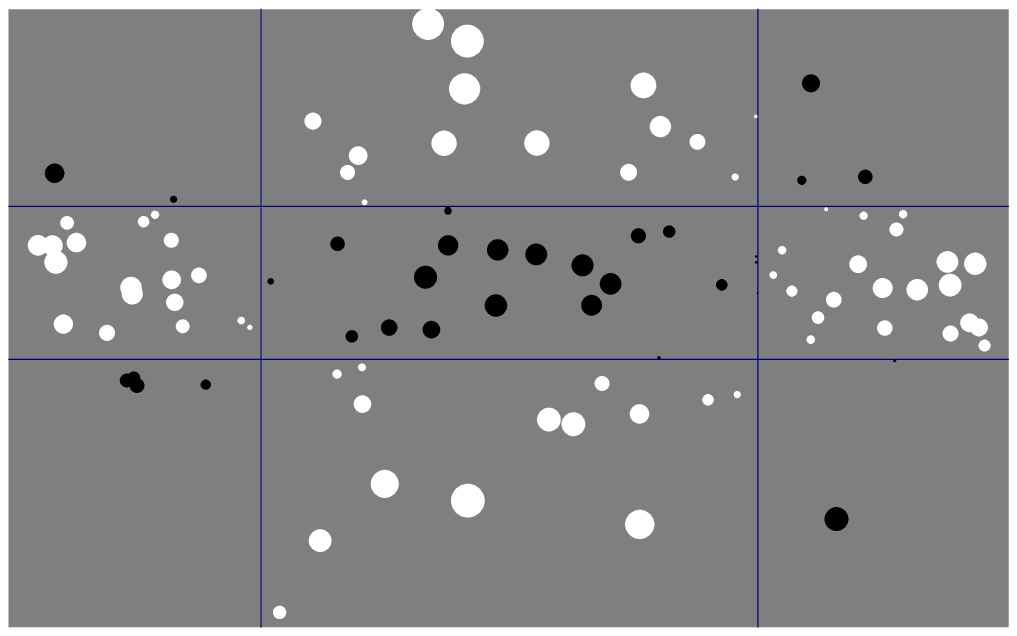}}
\caption{Representation of weights $W_i^{(k,l)}$ contributing to
$\hat\rho_{kl}$ for data in Figure~\ref{4scats}. The meaning of the
dots is otherwise the same as in
Figure~\ref{4scats2}.}\label{4scats3}
\end{figure}

For the data in Figure~\ref{4scats}(a), we found that
$\hat\rho^*=.36$ ($p=.000$), i.e., there is significant association.
In Figure~\ref{4scats2}(a), the data are plotted again, this time
each sample element $(X_i,Y_i)$ is represented with size
proportional to $|W_{i}|$; black dots represent positive $W_i$ and
white represent negative $W_i$. In all these plots, the total area
of the dots is scaled to a constant to make it easier to study them.
From Figure~\ref{4scats2}(a),we see that the association seems to
consist of a linearity in the data. It may be worthwhile however to
check if there are other forms of association present by looking at
the component correlations of $\hat\rho^*$. We found two significant
components: $\hat\rho_{11}=.61$ ($p=.000$) and $\hat\rho_{22}=.38$
($p=.004$). In Figures~\ref{4scats3}(a) and~\ref{4scats3}(b) the
weights $W_{i}^{(11)}$ and $W_{i}^{(22)}$ are visualized. The
interpetation of the black and white dots is the same as above. The
gridlines correspond to the zeroes of the marginal eigenfunctions.
Therefore, within any rectangle the dots have the same color. Also
in this case, the plots point to a linearity in the data.

For the data in Figure~\ref{4scats}(b), we found $\hat\rho^*=.17$
($p=.000$) and we found two significant component correlations:
$\hat\rho_{21}=-.78$ ($p=.000$) and $\hat\rho_{42}=.54$ ($p=.000$).
The plots in Figures~\ref{4scats2}(b),~\ref{4scats3}(c)
and~\ref{4scats3}(d) all point to a curved relationship.

For the data in Figure~\ref{4scats}(c), we found that
$\hat\rho^*=.11$ ($p=.001$). There is significant association at the
5\% level, but the evidence is not as overwhelming as in the
previous two cases. We only found one significant component
correlation: $\hat\rho_{12}=-.51$ ($p=.000$).
Figures~\ref{4scats2}(c) and~\ref{4scats3}(e) both point towards and
increase in dispersion of the $Y$ variable as $X$ increases.

For the data in Figure~\ref{4scats}(d), we found that
$\hat\rho^*=.02$ ($p=.522$). This time, the test based on $\rho^*$
does not yield a significant association. However, there is one
significant component correlation: $\hat\rho_{22}=-.38$ ($p=.004$).
Figure~\ref{4scats3}(f) indicates that the association is due to an
increase and then decrease in dispersion. Since $\hat\rho^*$ is not
significant, we should refrain from giving an interpretation based
on Figure~\ref{4scats2}(d).

\begin{figure}[tbp]
 \subfigure[Most common pattern]{\includegraphics[width=60mm]{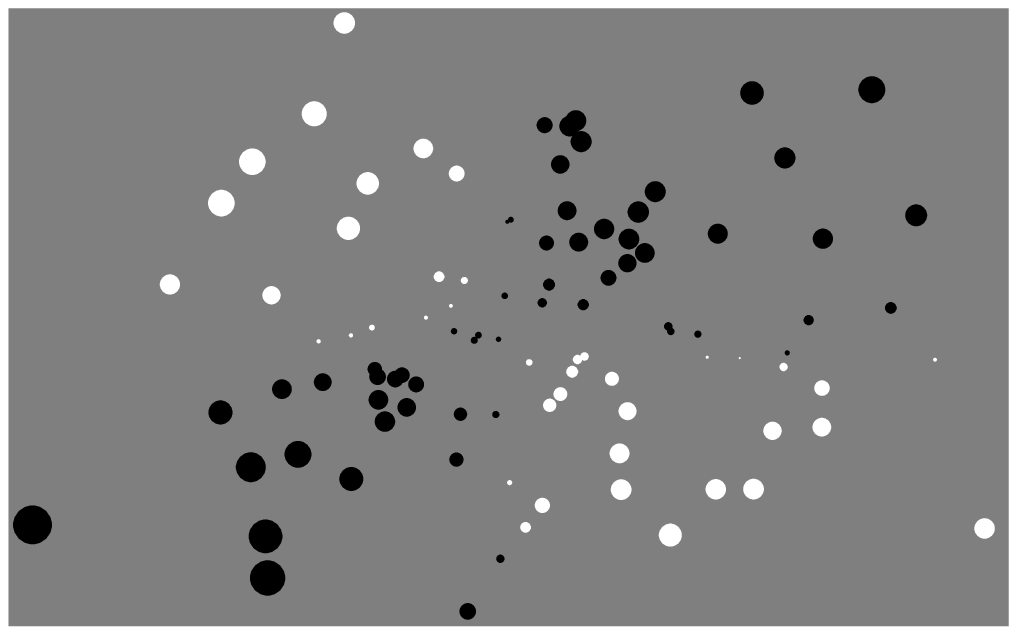}}
\hfill
 \subfigure[Less common pattern]{\includegraphics[width=60mm]{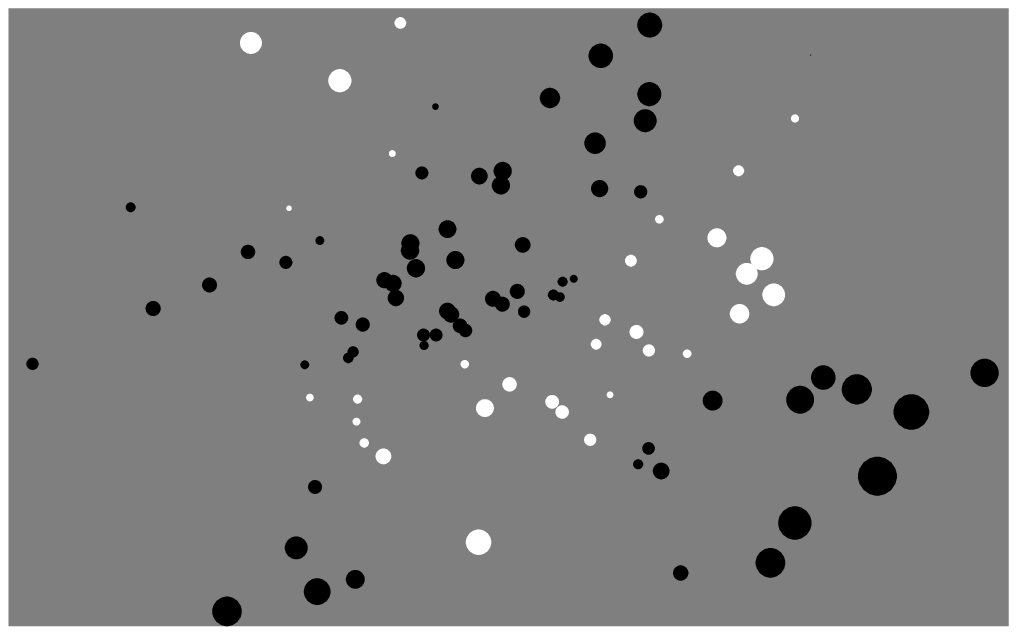}}
\caption{Representation of weights $W_i$ for data drawn from
distributions satisfying independence.}\label{4scats4}
\end{figure}

For comparative purposes, we plotted the weights $W_i$ for data
drawn from a distribution satisfying independence in
Figure~\ref{4scats4}. Both sets consist of 100 sample elements. The
most common pattern is that of Figure~\ref{4scats4}(a), with two
diagonally opposing clusters of black dots and two diagonally
opposing clusters of white dots. In a very limited investigation,
this type of pattern occured about half of the time. Otherwise more
complex patterns were obtained, such as the one in
Figure~\ref{4scats4}(b). These figures indicate that it doesn't seem
to make sense to interpret this kind of plot if $\hat\rho^*$ is not
significant, such as Figure~\ref{4scats2}(d).

Concluding, we see that significance tests combined with an
inspection of the two types of plots in Figures~\ref{4scats2}
and~\ref{4scats3} can give us insight into the kind of association
present between two random variables.

\subsection{Mental health data}\label{sec cat}

\begin{table}
\begin{center}
\begin{tabular}{|l|c|c|c|c|}\hline
  &\multicolumn{4}{|l|}{Mental Health Status}\\  \hline
Parents'      &      & Mild      & Moderate  &  \\
Socioeconomic &      & Symptom   & Symptom   &      \\
Status        & Well & Formation & Formation & Impaired\\ \hline %
A (high) & 64 & 94& 58& 46 \\  %
B & 57 & 94 & 64 & 40\\  %
C & 57 & 105 & 65 & 60\\ %
D & 72 & 141 & 77 & 94\\ %
F & 36 & 97 &54 & 78\\  %
G (low) & 21 & 71 & 54 & 71\\ \hline
\end{tabular}
\end{center}\caption{Cross-classification of Mental Health Status and Socioeconomic Status}
\label{table ment}
\end{table}

\begin{figure}[htbp]
\begin{center}
 \subfigure[$\hat\rho^*=.02$]
 {\includegraphics[width=36mm]{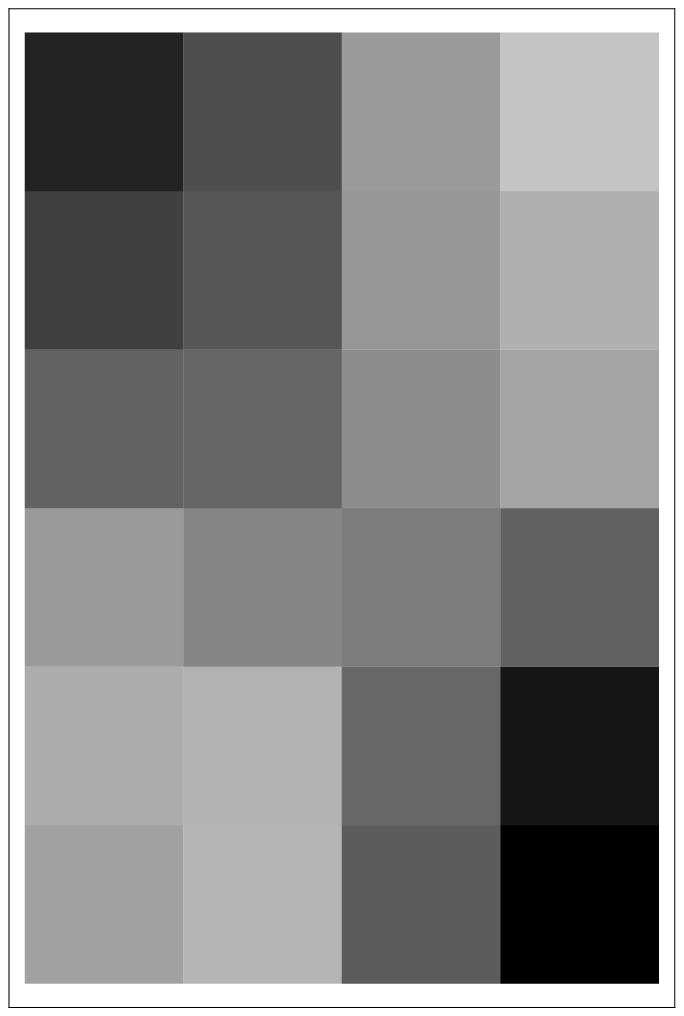}}
 \hspace{6mm}
 \subfigure[$\hat\rho_{11}=.13$]
 {\includegraphics[width=36mm]{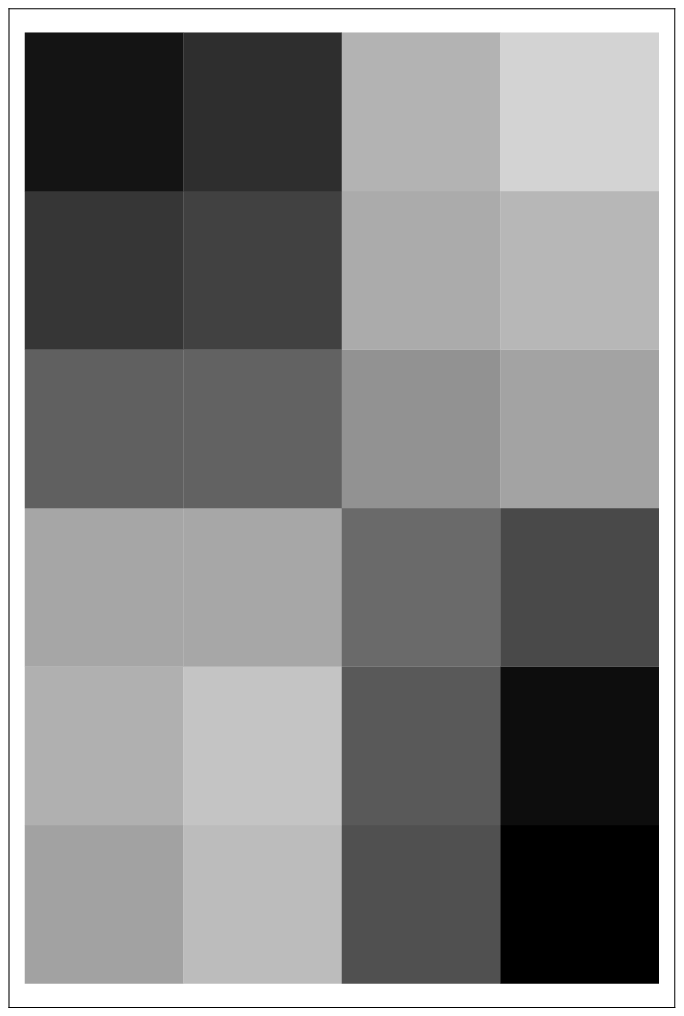}}
 \hspace{6mm}
 \subfigure[$\hat\rho_{13}=.08$]
 {\includegraphics[width=36mm]{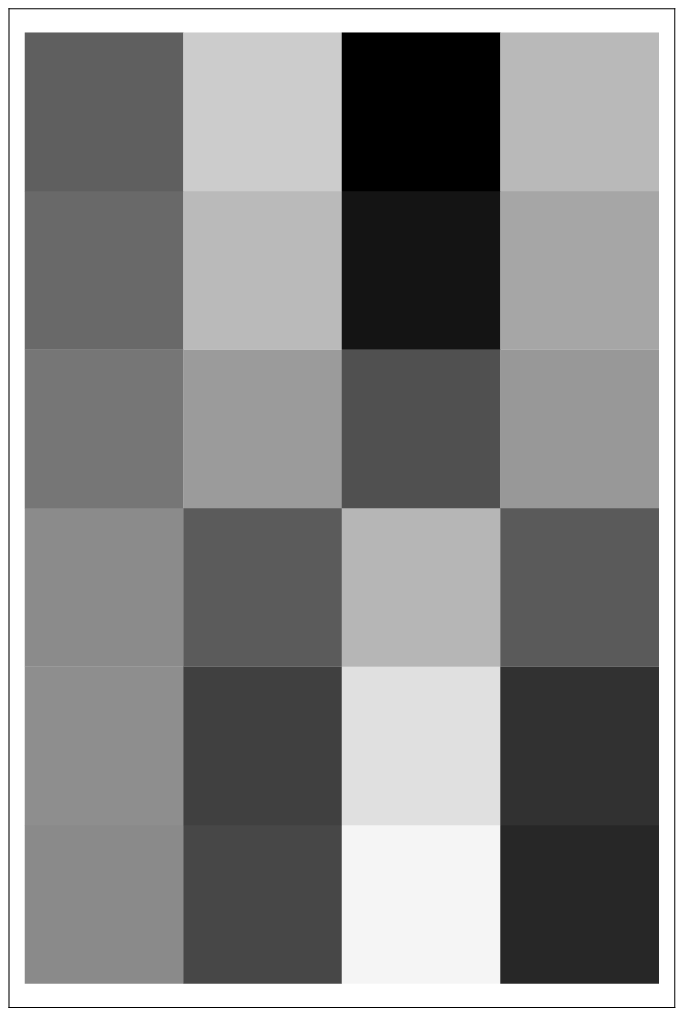}}
 \hfill
\end{center}
\caption{Representation of weights $W_{ab}$, $W_{ab}^{(11)}$ and
$W_{ab}^{(13)}$, respectively, for mental health data of
Table~\ref{table ment}. The darker the shade of gray in cell $(a,b)$
the larger $W_{ab}$; $W_{ab}=0$ is represented by an intermediate
shade of gray.} \label{table ment2}
\end{figure}

Table~\ref{table ment} describes the relationship between child's
mental impairment and parents' socioeconomic status  for a sample of
residents of Manhattan \cite[and references
therein]{goodman85,agresti02}. Goodman used this table to illustrate
various association models for categorical data, including the
so-called linear by linear association model, the row and columns
effects model, and correspondence analysis based on canonical
correlations. Here we illustrate the use of $\hat\rho^*$ and its
components as yet an alternative method for analyzing these data. We
relied on asymptotic $p$-values because with 1670 observations
approximate evaluation of the permutation tests would have been too
time consuming using our implementation.

In the categorical case, for an $I\times J$ contingency table, it
suffices to calculate weights for the cells, i.e., it is not
necessary to calculate separately a weight for each individual
observation. For an observation $(X_i,Y_i)$ in cell $(a,b)$, the
weight $W_i$ reduces to
\begin{eqnarray*}
   W_{ab} = p_{ab}\frac{\sum_{i=1}^I\sum_{j=1}^Jp_{ij}h_{\hat F_1}(i,a)h_{\hat F_2}(j,b)}{\sqrt{\hat\kappa(X,X)\hat\kappa(Y,Y)}}
\end{eqnarray*}
where $p_{ab}$ is the proportion of observations in cell $(a,b)$.
Similarly, the weights belonging to component correlation
$\rho_{kl}$ are
\begin{eqnarray*}
   W_{ab}^{(k,l)} = p_{ab}g_{1k}(a)g_{2l}(b)
\end{eqnarray*}
for $a=1,\ldots,I$, $b=1,\ldots,J$, $k=1,\ldots,I$ and
$l=1,\ldots,J$. Note that
\begin{eqnarray*}
   \rho^* = \sum_{a=1}^I\sum_{b=1}^J W_{ab}
\end{eqnarray*}
and
\begin{eqnarray*}
   \rho_{kl} = \sum_{a=1}^I\sum_{b=1}^J W_{ab}^{(k,l)}
\end{eqnarray*}

We found that $\hat\rho^*=.02$ ($p=.000$), i.e., there is
significant association in the data. The weights $W_{ab}$ for the
cells are represented in Figure~\ref{table ment2}(a). Here, the
grayscale represents the size of $W_{ab}$: the darker the cell, the
larger $W_{ab}$; $W_{ab}=0$ is represented by a fixed intermediate
shade of gray. From Figure~\ref{table ment2}(a), it can be seen that
most of the association is of a monotone nature: the higher the
parents' socioeconomic status, the better the mental health status
of their children. We also investigated the component correlations
and found two components to be significant at the 5\% level:
$\hat\rho_{11}=.13$ ($p=.000$) and $\hat\rho_{13}=.08$ ($p=.026$).
In Figures~\ref{table ment2}(b) and~\ref{table ment2}(c) we
represented the $W_{ab}^{(11)}$ and $W_{ab}^{(13)}$ using grayscales
as above. From Figure~\ref{table ment2}(b), we see that
$\hat\rho_{11}$ indicates linearity again. However, in
Figure~\ref{table ment2}(c) we see evidence of some nonlinearity in
the data, namely an apparent reversal of the association if only the
middle categories `Mild Symptom Formation' and `Moderate Symptom
Formation' are considered. Hence, it appears that the association
which is present in the data cannot be fully explained by linearity.

\subsection{Norwegian stock exchange}\label{sec time}

\begin{figure}[htbp]
\begin{center}
 \subfigure[Time series of stock price index]
 {\includegraphics[width=7cm]{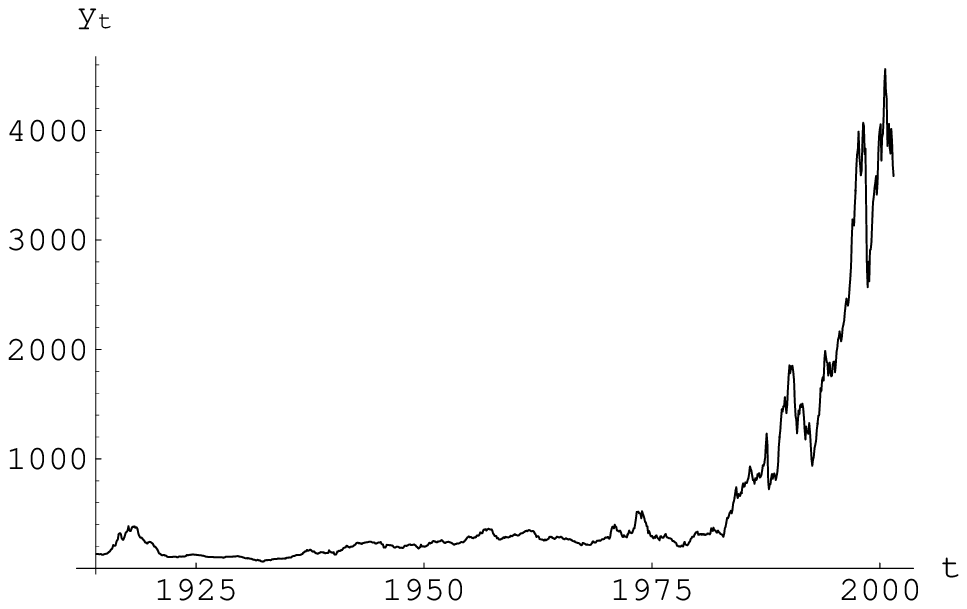}}
\end{center}
\hfill \\
\begin{center}
 \subfigure[Scatterplot of $(Z_t,W_t)$ showing association between successive jumps in Figure~\ref{fig norway}(a).
 $\hat\rho^*=.09$]
 {\includegraphics[width=7cm]{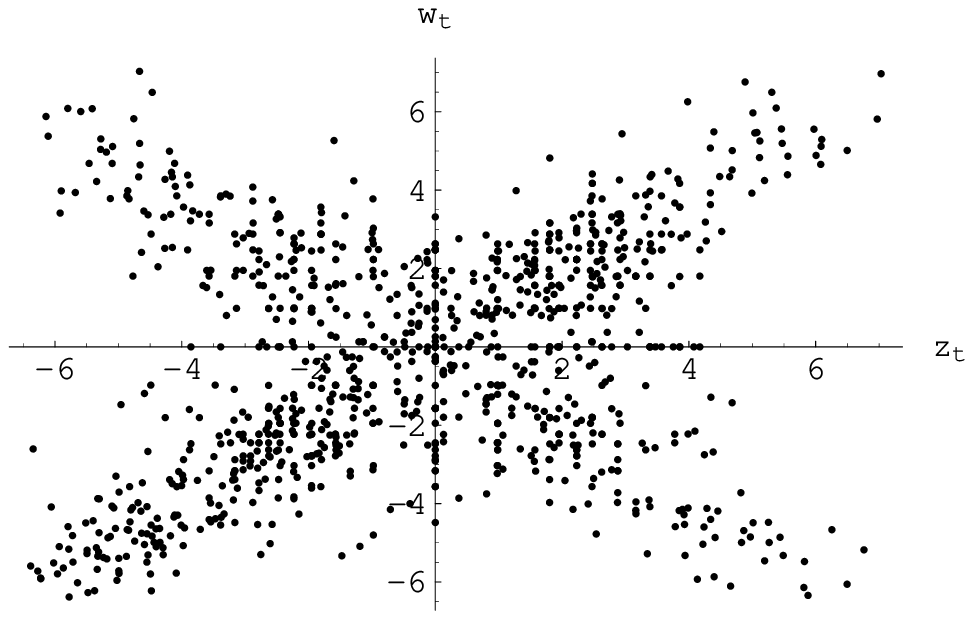}}
\end{center}
\hfill \\
\begin{center}
 \subfigure[Weights $W_i$ for data in Figure~\ref{fig norway}(b)]
 {\includegraphics[width=7cm]{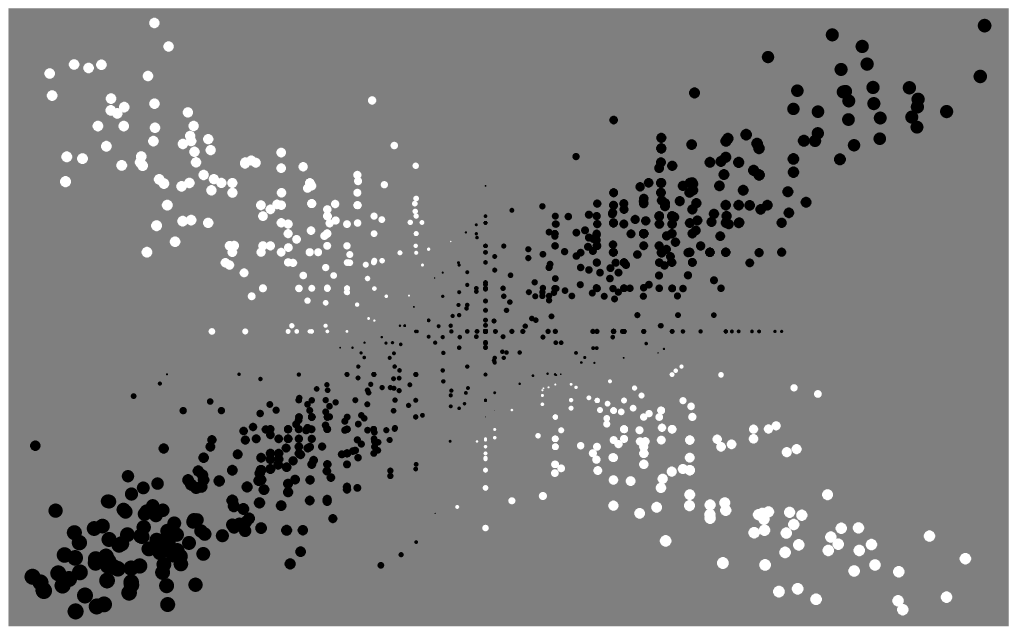}}
\end{center}
\caption{Monthly Norwegian stock price indices, 1914-2001}
\label{fig norway}
\end{figure}

The Norwegian stock exchange data represented in Figure~\ref{fig
norway}(a) yield an example with an especially interesting form of
association. In Figure~\ref{fig norway}(a) the original time series
data $Y_t$ are plotted. In Figure~\ref{fig norway}(b) we plotted
\begin{eqnarray*}
   (Z_t,W_t)=(\mbox{arcsinh}(Y_{t}-Y_{t-1}),\mbox{arcsinh}(Y_{t+1}-Y_t))
\end{eqnarray*}
The arcsinh transformation was done to make the marginal
distributions less heavy tailed. We found a highly significant
association with $\hat\rho^*(Z,W)=.09$ ($p=.000$). From
Figure~\ref{fig norway}(b) we can already interpret the association:
large jumps (up or down) in stock prices tend to be followed by
large jumps, and small jumps by small jumps, indicating periods of
volatility. In Figure~\ref{fig norway}(c), the weights $W_{i}$ are
represented by the size and color of the dots (see Section~\ref{sec
art} for further explanation). This plot points to the same
conclusion that large jumps are followed by large jumps and small
jumps by small jumps. Note that in this case, not only the positive
weights (the black dots) but also the negative weights (the white
dots) are highly indicative of association. The plot indicates that
the up-arm (with the black dots) is `heavier' than the down-arm
(with the white dots), that is, there is evidence that in the data
generating process a jump tends to be of the same sign as the
previous jump.

\begin{figure}[htbp]
\begin{center}
 \subfigure[$\hat\rho_{11}=.25$]
 {\includegraphics[width=4cm]{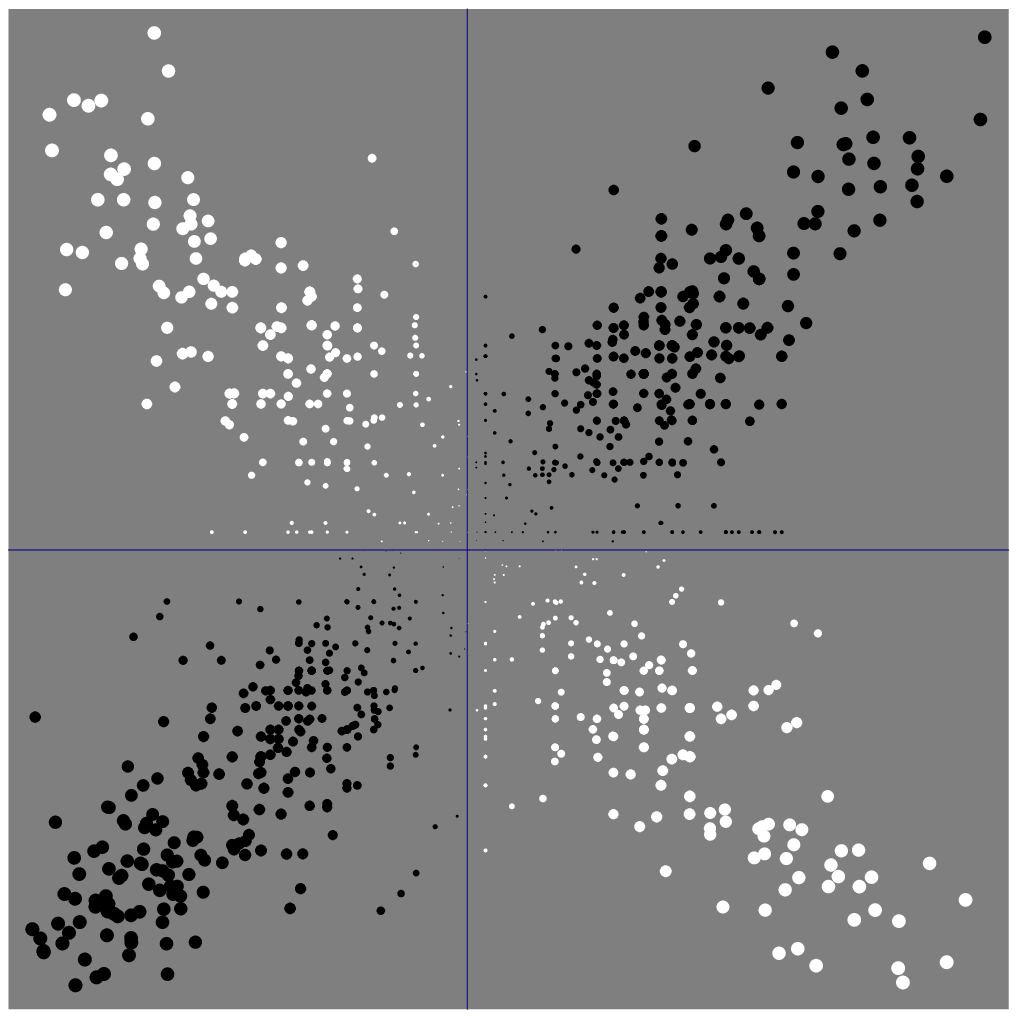}}
\hfill
 \subfigure[$\hat\rho_{12}=.13$]
 {\includegraphics[width=4cm]{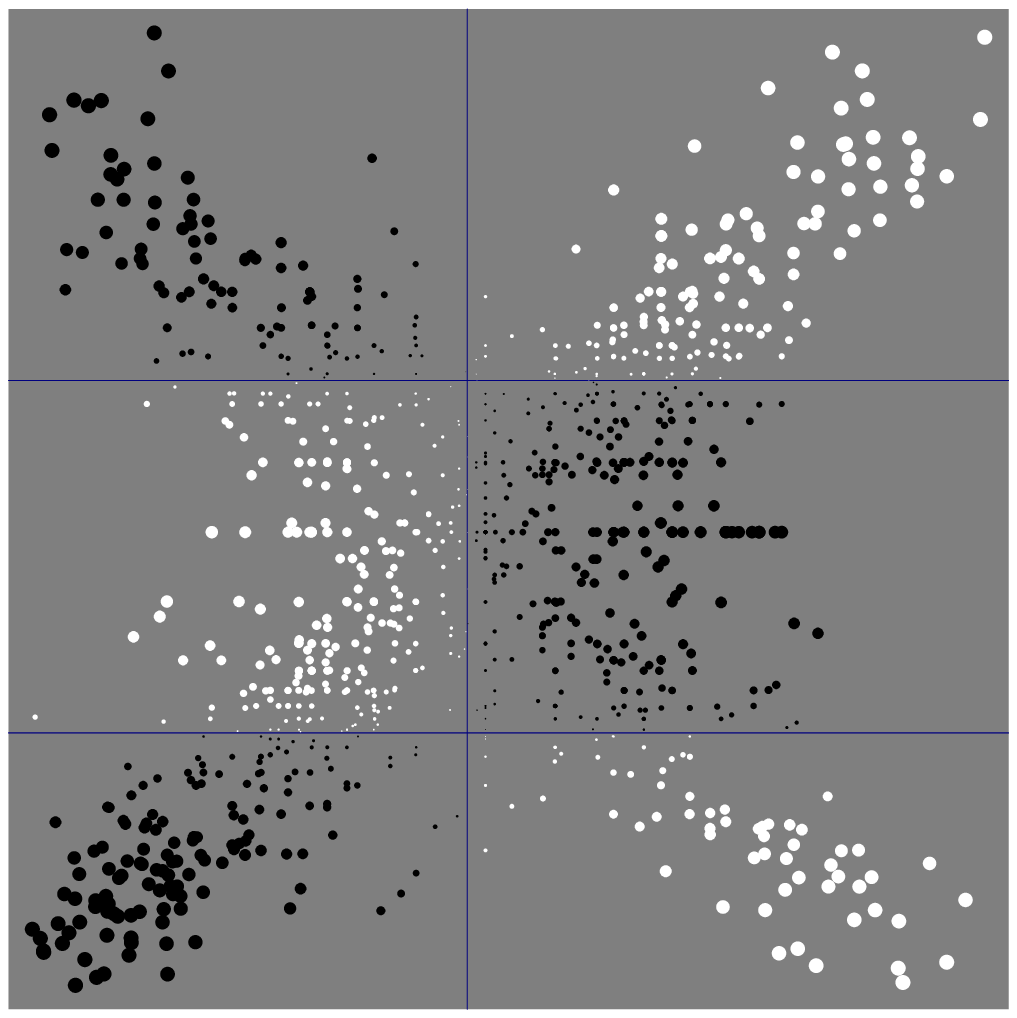}}
\hfill
 \subfigure[$\hat\rho_{22}=.64$]
 {\includegraphics[width=4cm]{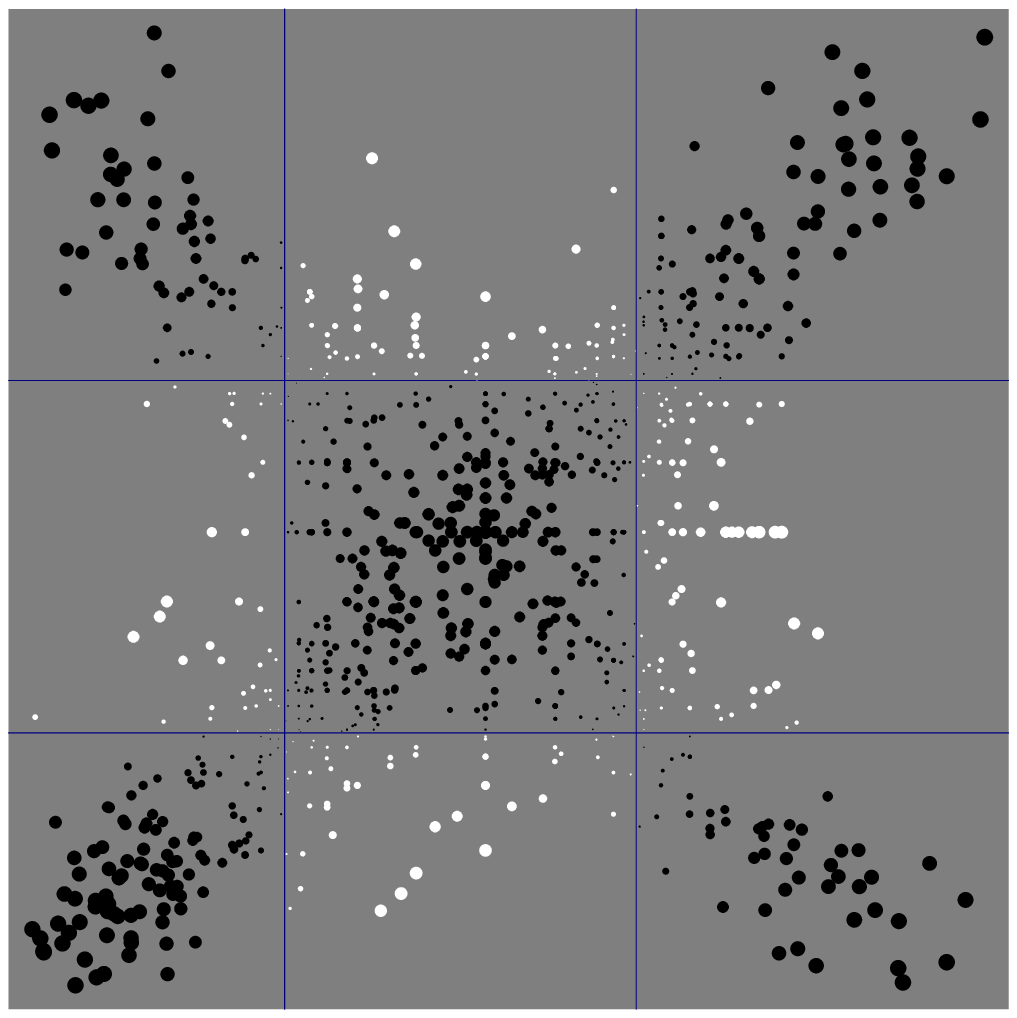}}
\end{center}\hfill\\
\begin{center}
 \subfigure[$\hat\rho_{24}=.19$]
 {\includegraphics[width=4cm]{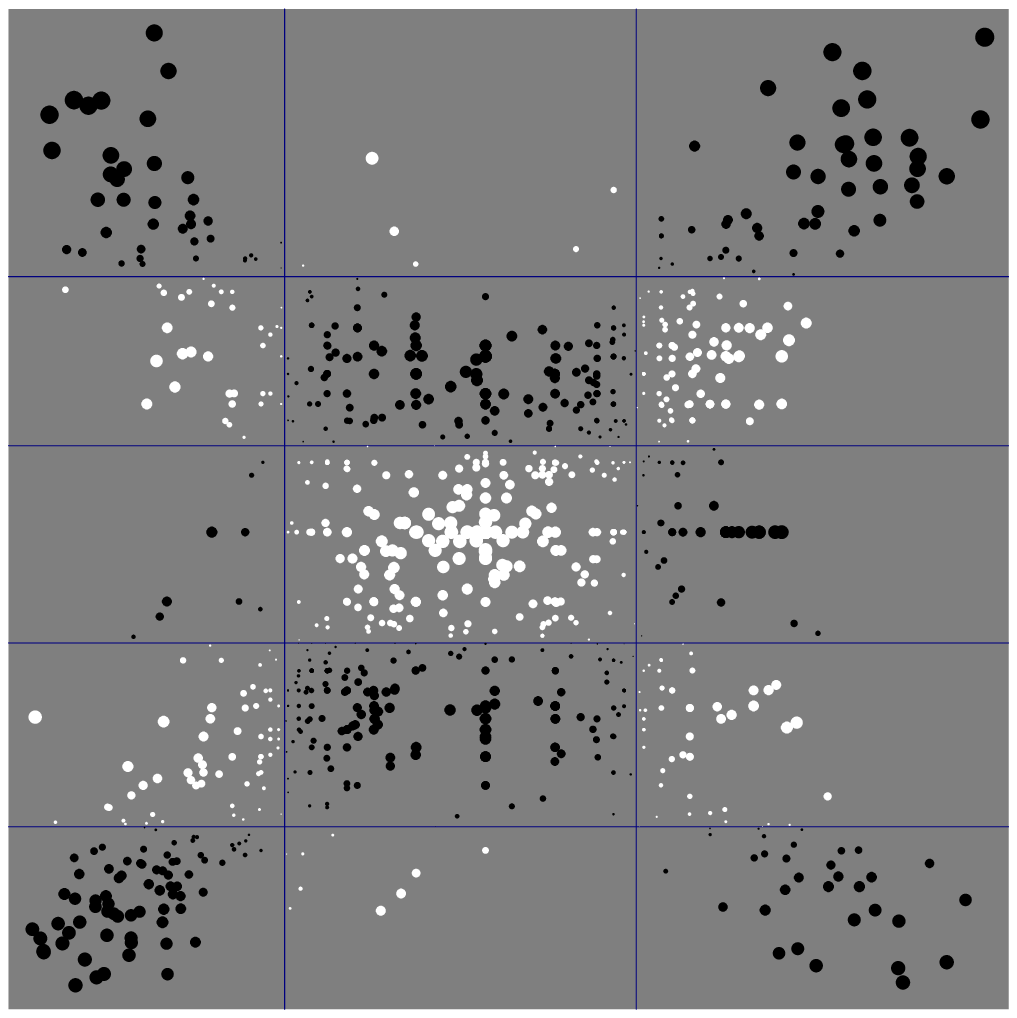}}
\hfill
 \subfigure[$\hat\rho_{33}=.19$]
 {\includegraphics[width=4cm]{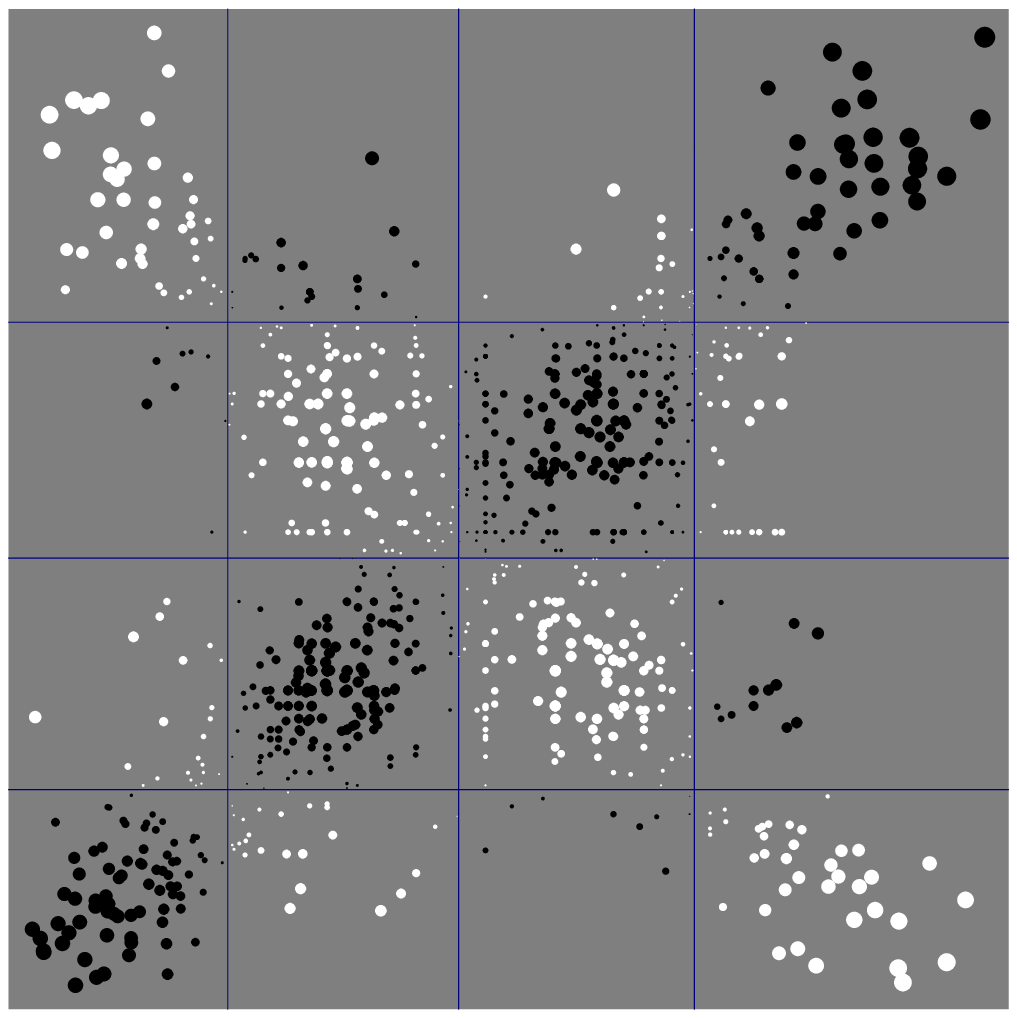}}
\hfill
 \subfigure[$\hat\rho_{44}=.38$]
 {\includegraphics[width=4cm]{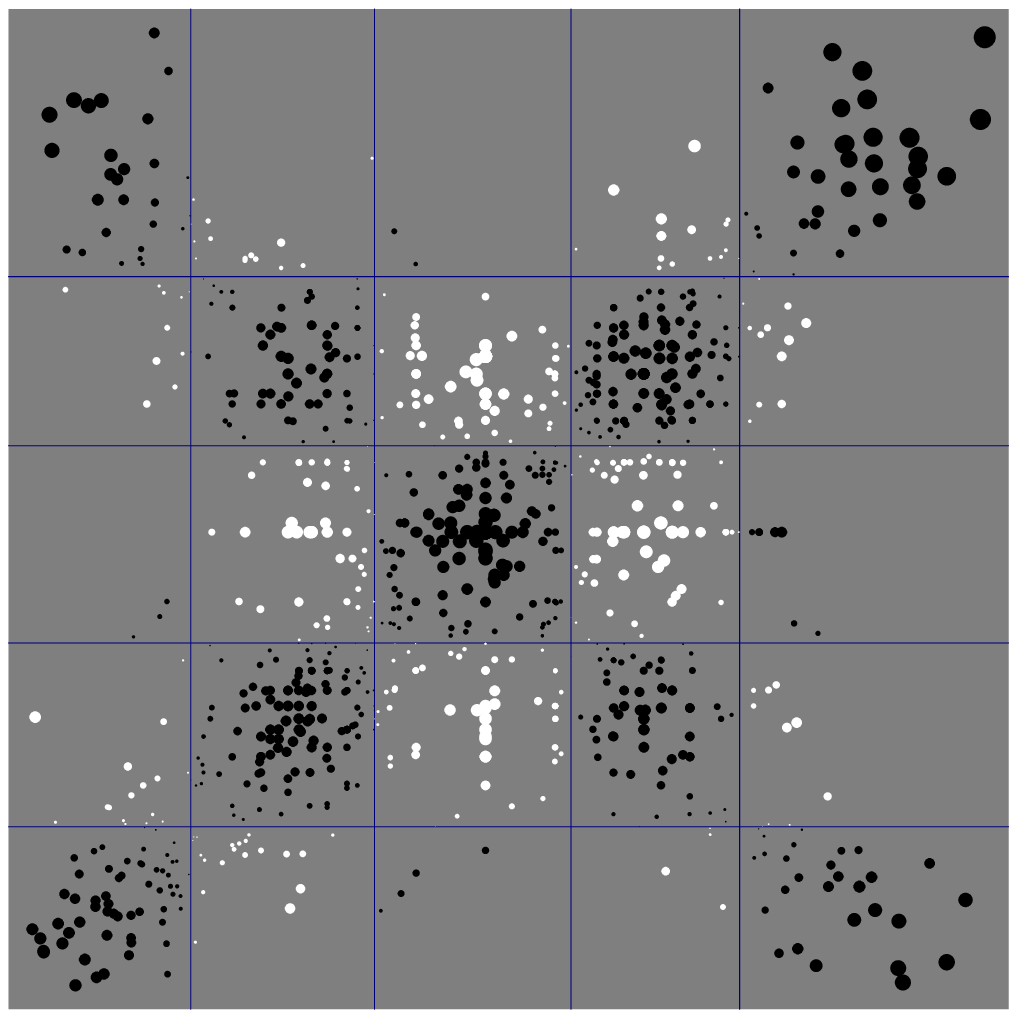}}
\end{center}\hfill\\
\begin{center}
 \subfigure[$\hat\rho_{66}=.25$]
 {\includegraphics[width=4cm]{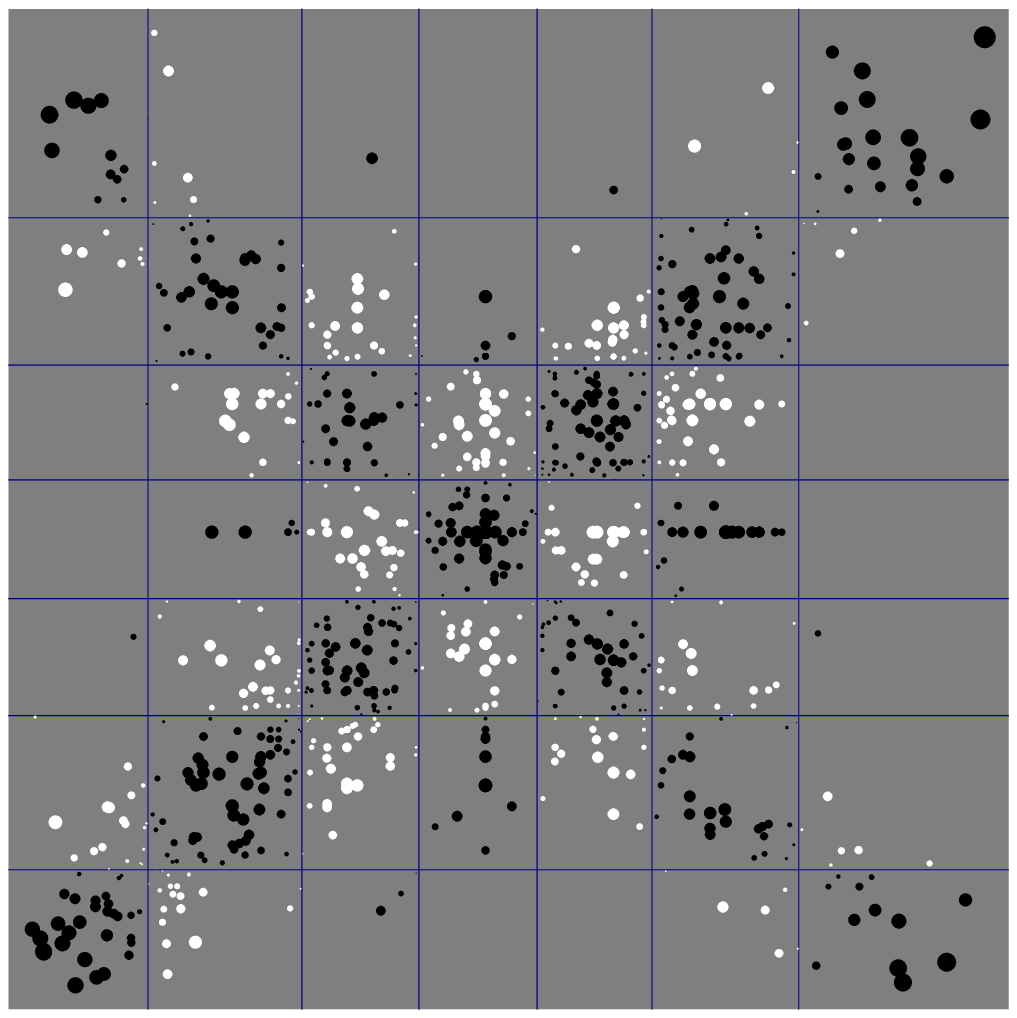}}
\end{center}\hfill\\
\caption{Representation of weights contributing to $\hat\rho_{ij}$
for data in Figure~\ref{fig norway}(b)} \label{fig norway 2}
\end{figure}

Seven component correlations were found to be significant at the 5\%
level after applying the Bonferroni correction: $\hat\rho_{11}=.25$,
$\hat\rho_{12}=.13$, $\hat\rho_{22}=.64$, $\hat\rho_{24}=.19$,
$\hat\rho_{33}=.19$, $\hat\rho_{44}=.38$ and $\hat\rho_{66}=.25$,
all with $p=.000$. In Figure~\ref{fig norway 2}, the weights
corresponding to the component correlations are represented.
Figures~\ref{fig norway 2}(c), (d), (f) and (g) point to the
cross-like nature of the data. Figures~\ref{fig norway 2}(a) and (e)
indicate that the up-arm is heavier than the down-arm. We did not
find a meaningful explanation for Figure~\ref{fig norway 2}(b).

\subsection{Discussion}

If a researcher investigating the association between two variables
decides on the use of $\rho^*$, we recommend the following approach.
First a test of the significance of $\hat\rho^*$ should be done,
and, if found to be significant, the weights $W_i$ should be
visualized as described above in order to determine the nature of
the association. If this does not yield the desired insight, it can
be worthwhile to investigate the component correlations and
visualize the corresponding weights $W_i^{(k,l)}$. These components
form a (unique) orthogonal decomposition of the
`infinite-dimensional' object $\rho^*$ into `one-dimensional'
objects $\rho_{kl}$, and the orthogonality ensures, in a limited
sense, that the different components measure different things; by
the latter we mean that for large samples and close to independence,
the sample component correlations are approximately independent. The
question may arise: why not investigate correlations between other
sets of orthogonal functions? A sketch of an answer is as follows.
Because of the various optimality properties of the eigenfunctions
in describing the marginal kernels, these component correlations are
likely to be a better choice for investigating the deviation of
$\rho^*$ from zero than correlations between arbitrarily chosen
functions for the marginal distributions. One way to make this
intuitive is as follows: in those regions where the marginal
distributions are sparse, the eigenfunctions vary relatively slowly
(in second derivative sense, see remark after Lemma~\ref{cor2}), and
therefore power of a test based on a component correlation will be
concentrated in those regions where there are many observations.
Thus, if we believe $\rho^*$ to be a good measure of deviation from
independence, then good `one-dimensional' objects to look at are the
correlations between the marginal eigenfunctions of $h_{F_1}$ and
$h_{F_2}$.



\section*{Acknowledgements}

The author would like to thank the following institutes where he has
been employed for providing a stimulating environment in which to
carry out this research: The Methodology Department of the Faculty
of Social Sciences at Tilburg University (Tilburg, the Netherlands)
and EURANDOM (Eindhoven, the Netherlands).


\bibliographystyle{apacite}
\bibliography{stats}

\end{document}